\documentclass[11pt,amsfonts]{amsart}

\input fullpage.sty

\def\endproof{\qed\smallskip}
\def\blacksquare{\hbox to .60em{\vrule width .60em height .60em}}

\newtheorem{theorem}{Theorem}[section]

\newtheorem{corollary}[theorem]{Corollary}

\newtheorem{definition}[theorem]{Definition}

\newtheorem{lemma}[theorem]{Lemma}

\newtheorem{proposition}[theorem]{Proposition}

\theoremstyle{remark}
\newtheorem{remark}[theorem]{\bf Remark}

\begin{document}

\title[]{Scalar Curvature and the Existence of Geometric Structures on 
3-Manifolds, II}

\author[]{Michael T. Anderson}

\thanks{Partially supported by NSF Grants DMS 9802722 and 0072591}
\thanks{2000 {\it Mathematics Subject Classification}. Primary 58J60, 
57M50; Secondary 58E11, 53C20}

\maketitle

\tableofcontents

\abstract 
This paper studies the degeneration, i.e curvature blow-up, of 
sequences of metrics approaching the Sigma constant $\sigma(M)$ on 
3-manifolds $M$ with $\sigma(M) \leq 0$. The degeneration is related to 
the sphere decomposition of $M$ in case $M$ is $\sigma$-tame.
\endabstract

\setcounter{section}{-1}

\section{Introduction.}
\setcounter{equation}{0}

 This paper is a continuation of [1], and is mainly concerned with the 
Sphere conjecture of [1,\S 0]. We recall the basic issues and set-up.

 Let $M$ be a closed oriented 3-manifold and $\sigma (M)$ the Sigma 
constant of $M$, i.e. the supremum of the scalar curvatures of unit 
volume Yamabe metrics on $M$. Throughout the paper, it is assumed that 
\begin{equation} \label{e0.1}
\sigma (M) \leq  0; 
\end{equation}
this is the case if for instance $M$ has a $K(\pi ,1)$ factor in its 
sphere decomposition, c.f. [11]. Consider the $L^{2}$ norm of the 
negative part of the scalar curvature as a functional on the space of 
metrics on $M$, i.e.
\begin{equation} \label{e0.2}
\mathcal{S}_{-}^{2}(g) = \bigl( v^{1/3}\int_{M}(s^{-})^{2}dV_{g} 
\bigr)^{1/2}, 
\end{equation}
where $s^{-} = \min(s, 0)$, $s = s_{g}$ is the scalar curvature of the 
metric $g$, $v$ is the volume of $(M, g)$. The volume power in (0.2) is 
chosen so that $\mathcal{S}_{-}^{2}$ is scale invariant. By [2, 
Prop.3.1], 
\begin{equation} \label{e0.3}
\inf \mathcal{S}_{-}^{2} = |\sigma (M)|, 
\end{equation}
under the assumption (0.1), so that a minimizing sequence for 
$\mathcal{S}_{-}^{2}$ has similar characteristics to a maximizing 
sequence of Yamabe metrics on $M$.

 Conjectures I and II of [1] concern the limiting behavior of certain 
minimizing sequences $\{g_{i}\}$ for $\mathcal{S}_{-}^{2}$ on 
irreducible 3-manifolds. These conjectures imply the geometrization 
conjecture for 3-manifolds satisfying (0.1). In [1, Thm.0.2], the 
validity of these two conjectures was proved for {\it  tame}  
3-manifolds, i.e. 3-manifolds $M$ which admit a minimizing sequence of 
unit volume metrics $\{g_{i}\}$ for $\mathcal{S}_{-}^{2}$ such that, 
for some $K <  \infty ,$
\begin{equation} \label{e0.4}
\mathcal{Z}^{2}(g_{i}) = \int_{M}|z_{g_{i}}|^{2}dV_{g_{i}} \leq  K, \ \ 
{\rm as} \ \ i \rightarrow  \infty , 
\end{equation}
where $z = r -  \frac{s}{3}g$ is the trace-free Ricci curvature. c.f. 
also [5] for an outline of the proof.

 The assumption (0.4) is a rather strong topological assumption and it 
is not clear which 3-manifolds actually satisfy (0.4). For instance if 
$N_{i}$, $i = 1,2$, are hyperbolic 3-manifolds then $M = N_{1}\#N_{2}$ 
does not satisfy (0.4). The fundamental issue concerning Conjectures I 
and II above is in fact whether the converse holds, i.e. if $M$ is 
irreducible, is then $M$ necessarily tame.

\medskip

 In this paper, building on the work in [1], substantial progress is 
made on this problem. To explain this, let
\begin{equation} \label{e0.5}
\widetilde \sigma(g) = \mathcal{S}_{-}^{2}(g) -  |\sigma (M)|, 
\end{equation}
so that $\widetilde \sigma(g) \geq  0$, and a sequence $\{g_{i}\}$ is 
minimizing for $\mathcal{S}_{-}^{2}$ if and only if $\widetilde 
\sigma(g_{i}) \rightarrow 0$. 

\begin{definition} \label{d 0.1.}
 A 3-manifold $M$ satisfying (0.1) is {\sf $\sigma$-tame}, if there 
exist constants $k <  \infty $ and $K <  \infty$, and some unit volume 
minimizing sequence $\{g_{i}\}$ for $\mathcal{S}_{-}^{2},$ such that
\begin{equation} \label{e0.6}
\widetilde \sigma(g_{i})^{k}\cdot \mathcal{Z}^{2}(g_{i}) \leq  K, \ \ 
{\rm as} \ \ i \rightarrow  \infty . 
\end{equation}
\end{definition}

 Thus, $M$ is $\sigma$-tame if there exists a sequence of unit volume 
metrics $\{g_{i}\}$ with $\widetilde \sigma(g_{i}) \rightarrow 0$ such 
that the (trace-free) curvature of $g_{i}$ blows up in $L^{2}$ at most 
as fast as some polynomial in $\widetilde \sigma^{-1}.$ This is of 
course a topological condition on $M$, which is clearly much weaker 
than the tameness condition (0.4), (which corresponds to the case $k = 
0$).

 In fact we believe that {\it  any}  closed 3-manifold $M$ with $\sigma 
(M) \leq  0$ is $\sigma$-tame, regardless of whether $M$ is irreducible 
or reducible; this will be pursued in subsequent work since the methods 
dealing with this question are not related to the main issues here. In 
[2,3], numerous examples or, more precisely, models for potentially 
minimizing sequences $\{g_{i}\}$ for $\mathcal{S}_{-}^{2}$ were 
constructed. All of these models satisfy (0.6) with $k = 1$. Using 
these models it is not difficult to show that if Conjectures I and II 
are valid, then any 3-manifold $M$, $\sigma (M) \leq $ 0, is 
$\sigma$-tame with $k=1$. Conversely, a much simpler proof of the Main 
Theorem below can be given for $\sigma$-tame 3-manifolds with $k=1$, 
(based on the work in [2, Thm.C]).

 The purpose of this paper is to prove the following result:

\smallskip
\noindent
{\bf Main Theorem.}
 {\it If $M$ is $\sigma$-tame, then Conjectures I and II of [1] are 
valid.}

\medskip

 As noted above, this implies the geometrization conjecture for 
$\sigma$-tame irreducible 3-manifolds $M$ with $\sigma (M) \leq  0$, as 
well as the evaluation of their Sigma constant $\sigma (M).$ We refer 
to [1,5] for background discussion on these Conjectures. Exact 
statements of the Main Theorem in the cases $\sigma (M) < 0$ and 
$\sigma (M) = 0$ are given in Theorems 7.6 and 7.7 respectively. 

\medskip

 The overall strategy of proof for the Main Theorem is briefly as 
follows; further details are given at the beginning of each section, 
c.f. also Remarks 1.3, 1.6 and 2.2. Basically, a modified version of 
the Sphere conjecture in [1] is proved for $\sigma$-tame 3-manifolds. 
Thus, suppose $M$ is $\sigma$-tame but not tame, so that the curvature 
of any minimizing sequence $\{g_{i}\}$ for $\mathcal{S}_{-}^{2}$ blows 
up in $L^{2},$ i.e. (0.4) fails. By far the most crucial issue of the 
paper is to locate ``natural'' 2-spheres $S^{2}$ from the geometry of 
the sequence $\{g_{i}\},$ as $i \rightarrow  \infty$. Under two 
relatively weak hypotheses, namely a degeneration hypothesis and a 
non-collapse hypothesis, it is shown that such 2-spheres arise as a 
certain locus where the curvature of $\{g_{i}\}$ blows-up, i.e. becomes 
very large, for $\widetilde \sigma(g_{i})$ small. This core issue of 
the paper requires most of the technical work.

 Given the existence of such 2-spheres, one then needs to prove that 
they are essential in $M$, i.e. do not bound 3-balls. This is done by 
means of the cut and paste or comparison methods already developed and 
used in [1, \S 2,\S 4] just for this purpose. The upshot is that $M$ is 
then necessarily a reducible 3-manifold.

 Recall that Conjectures I and II apply only to irreducible 
3-manifolds. Hence, if $M$ is irreducible, one of the two 
aforementioned hypotheses must fail. In the last section of the paper, 
\S 7, this is shown to imply that there exists a minimizing sequence 
which does not degenerate at all, so that $M$ is in fact tame. This 
proves the Main Theorem, as discussed preceding (0.4). In sum, if $M$ 
is $\sigma$-tame, but not tame, then $M$ is necessarily reducible.

\medskip

 The crucial issue is thus to locate geometrically natural 2-spheres 
from the geometry of the degeneration of a suitable minimizing sequence 
$\{g_{i}\}.$ The approach taken in the previous papers [1,2] was to 
find such spheres near regions where the curvature is blowing up at a 
{\it  maximal}  rate, since these are geometrically the most obvious 
places. Blow-up limits of the metric sequence based at such points are 
solutions of a certain elliptic system of equations, the 
$\mathcal{Z}_{c}^{2}$ equations, c.f. \S 1, studied in detail in [2]. 
The Sphere conjecture states that such limit solutions or metrics have 
an asymptotically flat end, c.f. also Remark 1.6; such an end thus 
carries a natural 2-sphere.

 While we believe the Sphere conjecture above is true, in this paper a 
somewhat different approach is taken. Namely, 2-spheres in $M$ are 
detected as a locus where the curvature of $\{g_{i}\}$ is blowing up, 
but blowing up at a, comparatively speaking, almost {\it  minimal}  
rate; (similar issues regarding the rate of curvature blow-up play a 
central role in [3]). The metric blow-up limits based at such points 
are {\it  flat}, and the near limiting metric behavior is governed by 
solutions of linearized equations closely analogous to the linearized 
$\mathcal{Z}_{c}^{2}$ equations. These linearized equations are derived 
in \S 2 and their implications for the near limiting geometry are 
analysed further in \S 3-\S 4. The main point, developed in \S 5, is 
the identification of 2-spheres in almost flat regions of rescalings 
$(M, g_{i}')$ of $(M, g_{i}),$ which enclose regions where the 
curvature is much larger. Such 2-spheres serve as a substitute for the 
spheres at infinity in the Sphere conjecture. Remark 3.6 explains 
briefly one important reason why this approach is advantageous to the 
analysis of general $\mathcal{Z}_{c}^{2}$ solutions.

\smallskip

 Thus, the most significant part of the paper is the work in \S 5, c.f. 
in particular Theorems 5.4, 5.9 and 5.12, given the preparatory work in 
\S 2-\S 4 leading up to \S 5. More detailed descriptions of the 
contents and the strategy of the arguments are given at the beginning 
of each section as well as in the Remarks.

\section{Background Setting and Results.}
\setcounter{equation}{0}

 In this section, we set the stage and discuss background material 
needed for the work to follow; most of this is taken from [2,4]. 
Remarks 1.3 and 1.6 give some preliminary perspective on the overall 
strategy of the paper.

 For a given $\varepsilon  >  0$, consider the scale-invariant 
functional $I_{\varepsilon}^{~-} = v^{1/3}\varepsilon\mathcal{Z}^{2} + 
\mathcal{S}_{-}^{2}$ on the space of metrics on the closed, oriented 
3-manifold $M$, $\sigma(M) \leq 0$, i.e.
\begin{equation} \label{e1.1}
I_{\varepsilon}^{~-} = \varepsilon v^{1/3}\int_{M}|z|^{2}dV + \bigl( 
v^{1/3}\int_{M}(s^{-})^{2}dV \bigr)^{1/2}. 
\end{equation}
Although seemingly complicated at first sight, it is explained in 
detail in [2, \S 1], (c.f. also [1]), why this is, in a natural sense, 
an optimal functional to consider from the point of view of 
geometrization of $M$ via direct methods in the calculus of variations.

 It is proved in [2,Thm.3.9] that, for any $\varepsilon  >  0$, there 
exists a domain $\Omega_{\varepsilon}$ and a complete $C^{2,\alpha}\cap 
L^{3,p}$ smooth Riemannian metric $g_{\varepsilon}$ on 
$\Omega_{\varepsilon}$ with $vol_{g_{\varepsilon}}\Omega_{\varepsilon} 
=$ 1, which realizes the infimum of $I_{\varepsilon}^{~-}.$ The domain 
$\Omega_{\varepsilon}$ has a finite number $q = q(\varepsilon , M)$ of 
components and is weakly embedded in $M$ in the sense that any domain 
with smooth compact closure in $\Omega_{\varepsilon}$ embeds in $M$. 
Further, there exists an exhaustion $\{K_{j}\}$ of 
$\Omega_{\varepsilon}$ by compact subsets such that $M\setminus K_{j}$ 
is a graph manifold. The pair $(\Omega_{\varepsilon}, g_{\varepsilon})$ 
is called a {\it minimizing pair} for $I_{\varepsilon}^{~-}.$ In brief, 
$(\Omega_{\varepsilon}, g_{\varepsilon})$ gives a geometric 
decomposition of $M$ w.r.t. the functional $I_{\varepsilon}^{~-}.$

 As $\varepsilon  \rightarrow $ 0, the metrics $g_{\varepsilon}$ form a 
minimizing sequence or family for $\mathcal{S}_{-}^{2},$ in that
\begin{equation} \label{e1.2}
\mathcal{S}_{-}^{2}(g_{\varepsilon}) \rightarrow  |\sigma (M)| \ \ {\rm 
and} \ \ \varepsilon\mathcal{Z}^{2}(g_{\varepsilon}) \rightarrow  0, \ 
\ {\rm as} \ \ \varepsilon  \rightarrow  0. 
\end{equation}
The essential issue is to understand the behavior of sequences of 
minimizing pairs $(\Omega_{\varepsilon}, g_{\varepsilon}),$ where 
$\varepsilon  = \varepsilon_{i}$ is some sequence with $\varepsilon_{i} 
\rightarrow $ 0. 

 The domain $\Omega_{\varepsilon}$ is empty if and only if $M$ is a 
graph manifold which does not admit a flat metric; in fact, if 
$\Omega_{\varepsilon}$ is non-empty, then no component of 
$\Omega_{\varepsilon}$ is a graph manifold. Similarly, it is apriori 
possible that a sequence $(\Omega_{\varepsilon_{i}}, 
g_{\varepsilon_{i}})$ is constant, i.e. independent of 
$\varepsilon_{i}$, as $\varepsilon_{i} \rightarrow 0$. In this case of 
course (0.4) holds, so the Main Theorem is already proved; thus [1, 
Thm.0.2] implies that $(\Omega_{\varepsilon}, g_{\varepsilon}) = (H, 
g_{o})$, where $g_{o}$ is a constant curvature metric with scalar 
curvature $\sigma (M) < $ 0, with $\partial H$ incompressible in $M$, 
or $\Omega_{\varepsilon} = M$ is flat if $\sigma (M) =$ 0. In these 
situations, the geometric decomposition of $M$ w.r.t. 
$I_{\varepsilon}^{~-}$ above already gives the geometrization of $M$, 
for any $\varepsilon  > $ 0, c.f. [1, \S 2] for further details.

 Thus, in effect, only the situation where $(\Omega_{\varepsilon}, 
g_{\varepsilon})$ is not empty and not constant in $\varepsilon $ needs 
to be considered. All of the work in \S 2-\S 7 is to prove that this 
situation can occur only if $M$ is reducible. Hence, if $M$ is 
irreducible, then the point is to prove that either 
$\Omega_{\varepsilon}$ is empty, (so $M$ is a graph manifold), or 
$(\Omega_{\varepsilon}, g_{\varepsilon})$ is constant in $\varepsilon $ 
as above, c.f. Theorems 7.2, 7.6, and 7.7.

 Abusing terminology slightly, a metric of constant negative curvature 
will be called hyperbolic, even if the sectional curvature is not $- 
1$. The topology of $\Omega_{\varepsilon}$ may, apriori, change with 
$\varepsilon$. 

\medskip

 The metrics $g_{\varepsilon}$ satisfy the system of Euler-Lagrange 
equations for $I_{\varepsilon}^{~-},$ c.f. [2, Thm.3.3]:
\begin{equation} \label{e1.3}
L^{*}w  = \varepsilon\nabla\mathcal{Z}^{2} + \phi\cdot  g,
\end{equation}
\begin{equation} \label{e1.4}
2\Delta (w-\tfrac{\varepsilon s}{12}) + \tfrac{1}{4}sw = 
\tfrac{1}{2}\varepsilon|z|^{2} -  3c. 
\end{equation}
Here $\nabla\mathcal{Z}^{2}$ is the gradient of $\mathcal{Z}^{2} = 
\int|z|^{2}dV,$ c.f. (B.4), $\Delta = tr D^{2}$ where $D^{2}$ is the 
Hessian, and 
\begin{equation} \label{e1.5}
L^{*}u = D^{2}u -  \Delta u\cdot  g -  u\cdot  r, 
\end{equation}
where $r$ is the Ricci curvature, (all w.r.t. $(\Omega_{\varepsilon}, 
g_{\varepsilon})).$ The subscript $\varepsilon $ has been dropped from 
the notation in (1.3)-(1.4). The function $\phi $ is given by $\phi  = 
-\frac{1}{4}sw + c$, and the constant $c$ is given by
\begin{equation} \label{e1.6}
c = \frac{1}{12\sigma}\int 
(s^{-})^{2}+\frac{\varepsilon}{6}\int|z|^{2}, 
\end{equation}
where $\sigma  = (v^{1/3}\int (s^{-})^{2})^{1/2}.$ Observe that $\sigma 
$ is scale invariant and converges to $|\sigma (M)| = -\sigma (M)$ as 
$\varepsilon  \rightarrow $ 0, by (1.2). Most importantly, the function 
$w$ is given by
\begin{equation} \label{e1.7}
w = -  \frac{s^{-}}{\sigma} \geq  0, 
\end{equation}
so that the $L^{2}$ norm of $w$ over $\Omega_{\varepsilon}$ equals 1. 
The equation (1.4) is the trace of (1.3). The metric $g_{\varepsilon}$ 
is $C^{2,\alpha}\cap L^{3,p},$ for any $p <  \infty ,$ while the scalar 
curvature $s$ is Lipschitz. In any region where $s \neq $ 0, the metric 
is $C^{\infty}$, in fact real-analytic. The function $-w$ was denoted 
by $\tau$ in [1].

 Again the equations (1.3)-(1.4) are a complicated nonlinear system of 
elliptic PDE (in the metric $g$), but are explained and analysed in 
detail in [2, \S 1,3]. Essentially all of the work in this paper 
depends crucially on the exact form of the terms in these equations; 
the analysis developed is thus specific to these equations and may not 
be relevant to general variational problems on the space of metrics.

 Let $T = T_{\varepsilon} = \sup_{\Omega_{\varepsilon}}w.$ Then one has 
the estimate, (c.f. [2, (3.43)]),
\begin{equation} \label{e1.8}
1 \leq  T \leq  (1 + 2\frac{\varepsilon\mathcal{Z}^{2}}{\sigma})^{1/2}. 
\end{equation}
Thus, if $\sigma (M) < $ 0, then
\begin{equation} \label{e1.9}
T_{\varepsilon} \rightarrow  1, \ \ {\rm as} \ \ \varepsilon  
\rightarrow  0, 
\end{equation}
and so in particular
\begin{equation} \label{e1.10}
w_{\varepsilon} \rightarrow  1, a.e. 
\end{equation}
in the sense that for any $\mu  >  0$, the measure of the set where 
$|1- w_{\varepsilon}| \geq  \mu $ converges to 0, as $\varepsilon  
\rightarrow $ 0. Moreover, in general for $\sigma (M) \leq  0$, one has
\begin{equation} \label{e1.11}
\int_{\Omega_{\varepsilon}}|dw_{\varepsilon}|^{2}dV_{g_{\varepsilon}} 
\rightarrow  0, \ \ {\rm as} \ \  \varepsilon  \rightarrow  0, 
\end{equation}
c.f. [3, (3.39),(3.46)]. However, if $\sigma (M) = 0$ and $M$ is 
$\sigma$-tame, then (1.9) does not hold. Although it can be shown that 
$T_{\varepsilon}$ remains bounded as $\varepsilon  \rightarrow  0$, 
this will not be used here. Thus to unify these two situations, 
regardless of the behavior of $T$, the potential $w$ will {\it always}  
be renormalized as
\begin{equation} \label{e1.12}
u = \frac{w}{T}, 
\end{equation}
so that, by fiat,
\begin{equation} \label{e1.13}
\sup_{\Omega_{\varepsilon}}u = 1. 
\end{equation}
The points $x_{\varepsilon}\in\Omega_{\varepsilon}$ where 
$u(x_{\varepsilon}) \rightarrow 1$ play the central role throughout the 
paper.

 The equations (1.3)-(1.4) then read
\begin{equation} \label{e1.14}
L^{*}u  = \bar \varepsilon \nabla\mathcal{Z}^{2} + \bar \phi \cdot  g,
\end{equation}
\begin{equation} \label{e1.15} 
2\Delta (u-\tfrac{\bar \varepsilon s}{12}) + \tfrac{1}{4}su = 
\tfrac{1}{2}\bar \varepsilon |z|^{2} -  3\bar c, 
\end{equation}
where $\bar \varepsilon = \varepsilon /T, \bar c = c/T, \bar \phi = 
\phi /T.$ The equations (1.14)-(1.15) will be used throughout the paper.

\smallskip

 The non-negative function $u$ in (1.12) is of course just a 
renormalization of the non-positive part of the scalar curvature of 
$g_{\varepsilon}.$ It is viewed as a potential function on 
$(\Omega_{\varepsilon}, g_{\varepsilon}),$ since it satisfies the 
potential type equation (1.15). The main point of view of the paper is 
to understand the degenerating geometry of $(\Omega_{\varepsilon}, 
g_{\varepsilon})$ as $\varepsilon \rightarrow 0$ from analysis on the 
scalar curvature, or more precisely $u$, as a solution of the equations 
(1.14)-(1.15).

 While all the terms in (1.14)-(1.15) are important, the two most 
important are the potential $u$ and the constant term $\bar c.$ Most of 
the paper is concerned with the distribution of the values of $u$ over 
$\Omega_{\varepsilon},$ and its relation to the behavior of the 
curvature blow-up of $g_{\varepsilon}.$ It will be seen in \S 3, c.f. 
Theorem 3.5, that $c$ in (1.4) or $\bar c$ in (1.15) gives some 
important global control on the value distribution of $u$. Observe that 
$\bar c$ is the {\it  only}  global term in the equations 
(1.14)-(1.15). Lemma 1.7 below relates the two terms of $c$ in (1.4) to 
the condition that $M$ is $\sigma$-tame.

\medskip

 Let $\rho (x) = \rho_{\varepsilon}(x)$ be the $L^{2}$ curvature radius 
of $(\Omega_{\varepsilon}, g_{\varepsilon})$ at $x$; $\rho (x)$ is the 
largest radius such that, for all geodesic balls $B_{y}(s) \subset  
B_{x}(\rho (x)) \subset  (\Omega_{\varepsilon}, g_{\varepsilon}),$
\begin{equation} \label{e1.16}
\frac{s^{4}}{vol B_{y}(s)}\int_{B_{y}(s)}|r|^{2} \leq  c_{o}, 
\end{equation}
where $c_{o}$ is a fixed constant, c.f. [2,(2.1)], [4,\S 3]. The radius 
$\rho (x)$  measures the degree of curvature concentration near $x$; 
$\rho (x)$ is small if and only if the $L^{2}$ average of the curvature 
of $g_{\varepsilon}$ is large near $x$. The function $\rho $ is a 
Lipschitz function, with $|\nabla\rho| \leq  1$. 

 We recall the following result from [2,Rmk.7.3], (compare also with 
[3,Thm.3.3]), which will be used repeatedly throughout the paper.

\begin{proposition} \label{p 1.1.}
  There are constants $\upsilon_{o} > $ 0 and $\rho_{o} > \infty$, 
independent of $\varepsilon ,$ such that if $x\in U^{\upsilon_{o}} = 
\{x\in\Omega_{\varepsilon}: u(x) \geq  1-\upsilon_{o}\}$, then
\begin{equation} \label{e1.17}
\rho (x) \geq  \rho_{o}\cdot  \min(t_{\upsilon_{o}}(x), 1), 
\end{equation}
where $t_{\upsilon_{o}}(x) = dist(x, L_{\upsilon_{o}})$, and 
$L_{\upsilon_{o}} =\{x\in\Omega_{\varepsilon}: u(x) = 1-\upsilon_{o}\}$ 
is the $1-\upsilon_{o}$ level of u. Further, analogous higher order 
estimates also hold, i.e, for any $k \geq  0$,
\begin{equation} \label{e1.18}
|\nabla^{k}r|(x) \leq  c(k)\cdot  \max(t_{\upsilon_{o}}(x)^{-2- k}, 1), 
\ \  |\nabla^{k}du|(x) \leq  c(k)\cdot  \max(t_{\upsilon_{o}}(x)^{-1- 
k}, 1).  
\end{equation}
\end{proposition}

 The constant 1 on the right in (1.17)-(1.18) arises from the scalar 
curvature of $(\Omega_{\varepsilon}, g_{\varepsilon}).$ The estimate 
(1.17) can be re-expressed scale-invariantly as
\begin{equation} \label{e1.19}
\rho (x) \geq  \rho_{o}\cdot  \min(t_{\upsilon_{o}}(x), \rho_{s}(x)), 
\end{equation}
where $\rho_{s}$ is the $L^{2}$ scalar curvature radius, defined as in 
(1.16) with $s$ in place of $r$. Similarly, (1.18) holds when $1$ is 
replaced by the appropriate power of $\rho_{s}(x)$.

\medskip

 As throughout the previous papers, we use the $L^{2}$ Cheeger-Gromov 
theory, c.f. [7,8], [4, \S 3] to understand the limiting behavior of 
the metrics $(\Omega_{\varepsilon}, g_{\varepsilon})$ as $\varepsilon  
\rightarrow $ 0. Thus, suppose $x_{\varepsilon}\in\Omega_{\varepsilon}$ 
are base points such that 
\begin{equation} \label{e1.20}
u(x_{\varepsilon}) \geq  u_{o}, 
\end{equation}
for some $u_{o} > $ 0, and
\begin{equation} \label{e1.21}
\rho (x_{\varepsilon}) \sim  1, 
\end{equation}
i.e. $\rho (x_{\varepsilon})$ is bounded away from 0 and $\infty $ as 
$\varepsilon  = \varepsilon_{i} \rightarrow $ 0. Then there two 
possible behaviors.

\smallskip

{\bf Convergence (Non-Collapse).}
 Suppose $vol_{g_{\varepsilon}}B(\rho (x_{\varepsilon})) \geq  
\nu_{o},$ for some $\nu_{o} > $ 0. Let $U_{j} \equiv U_{j, \varepsilon} 
\subset  \Omega_{\varepsilon}$ be the component of 
$\{q_{\varepsilon}\in  \Omega_{\varepsilon}: 2^{-j}\leq  \rho 
(q_{\varepsilon}) \leq  2^{j}\}$ containing $x_{\varepsilon},$ and 
choose a sequence $j = j(\varepsilon ) \rightarrow  \infty $ as 
$\varepsilon  \rightarrow $ 0. Then for any sequence $\varepsilon  = 
\varepsilon_{i} \rightarrow $ 0, there is a subsequence such that 
$(U_{j}, g_{\varepsilon}, x_{\varepsilon})$ converges, (modulo 
diffeomorphisms), in the $C^{\alpha}$ topology, $\alpha  <  
\frac{1}{2}$, and weak $L^{2,2}$ topology, uniformly on compact 
subsets, to a maximal limit $(U, g_{o}, x)$. In particular, any compact 
domain in $U$ embeds in $U_{j},$ for $j$ sufficiently large, (in the 
subsequence).

\smallskip

{\bf Collapse.}
 Suppose $vol_{g_{\varepsilon}}B(\rho (x_{\varepsilon})) \rightarrow $ 
0 as $\varepsilon  = \varepsilon_{i} \rightarrow 0$. Then, for any 
given $R <  \infty $ and $\varepsilon $ sufficiently small, the domains 
$U_{j}(R) = U_{j}\cap B_{x_{\varepsilon}}(R)$ above are either 
Seifert-fibered spaces over a surface $V$, or are torus bundles over an 
interval. In both cases, the diameter of the fiber 
$F_{y_{\varepsilon}}$ over $y_{\varepsilon}\in U_{j}(R)$ satisfies 
$diam_{g_{\varepsilon}}F_{y_{\varepsilon}} \rightarrow $ 0 as 
$\varepsilon  \rightarrow $ 0 and the inclusion map of the fiber in 
$U_{j}(R)$ induces an injection on $\pi_{1}.$ In particular, there are 
finite covers $\widetilde U_{j}(R)$ of $U_{j}(R),$ of degree $d_{j} 
\rightarrow  \infty $ as $\varepsilon  \rightarrow $ 0, such that 
$(\widetilde U_{j}(R), g_{\varepsilon})$ does not collapse, and so has 
a $C^{\alpha}$ convergent subsequence as above. Letting $R = R_{j} 
\rightarrow  \infty , j = j(\varepsilon )$ as $\varepsilon  \rightarrow 
$ 0 gives rise to a maximal limit $(\widetilde U, g_{o}, x)$. The limit 
has a free isometric $S^{1}$ or $S^{1}\times S^{1} = T^{2}$ action, 
arising from the unwrapping of the collapse.

\smallskip

 There is of course some freedom in the choice of coverings unwrapping 
the collapse. It will be assumed throughout the paper that such 
coverings are chosen so that the diameter of the limiting $S^{1}$ 
orbit, or diameter and area of the limiting $T^{2}$ orbit, equal 1 at 
the limit base point $x$.

 Analogous to (1.16) define the volume radius $\omega (x)$ of 
$(\Omega_{\varepsilon}, g_{\varepsilon})$ by
\begin{equation} \label{e1.22}
\omega (x) = \sup\{r: \frac{vol B_{y}(s)}{s^{3}} \geq  \omega_{o}, 
\forall B_{y}(s)\subset B_{x}(r)\}, 
\end{equation}
where $\omega_{o}$ is a fixed (small) positive constant. Both $\rho $ 
and $\omega $ scale as distances. The two possibilities 
convergence/collapse above correspond to the situations $\omega 
(x_{\varepsilon}) \geq  v_{o}\cdot \rho (x_{\varepsilon})$ for some 
$v_{o} > $ 0, or $\omega (x_{\varepsilon}) <<  \rho (x_{\varepsilon})$ 
respectively. A similar collapse result holds if 
$B_{x_{\varepsilon}}(\rho (x_{\varepsilon}))$ is sufficiently 
collapsed, for any fixed $\varepsilon ,$ i.e. if $\omega 
(x_{\varepsilon}) \sim  \mu_{o}\cdot \rho (x_{\varepsilon}),$ for 
$\mu_{o}$ small. In this case, there are finite covers as above, of 
order about $\mu_{o}^{-1},$ such that in the covering space, $\omega 
(x_{\varepsilon}) \sim  \rho (x_{\varepsilon}).$

\smallskip

 Although the collapse case corresponds to a degeneration of the 
sequence of metrics, this situation is in fact somewhat easier to 
handle than the convergence situation, since the collapse may be 
(locally) unwrapped to obtain convergence to a limit which has special 
properties, namely a free isometric $S^{1}$ action. This extra 
structure is of importance and leads to simplifications not present in 
the convergence situation in general.

 In both cases, one has in fact $C^{\infty}$ convergence to the limit 
within regions satisfying (1.20)-(1.21), due to regularity estimates as 
in (1.18), c.f. [2, Rmk.4.3]. It is elementary but important to note 
that the Euler-Lagrange equations (1.14)-(1.15) are invariant under 
passing to covering spaces, as is, (modulo bounded factors), the 
curvature radius $\rho .$

\smallskip

 The metric boundary $\partial U$ of $U$ is the limiting locus where 
the curvature of $(U_{j}, g_{\varepsilon})$ blows-up. Thus 
$p\in\partial U$ if and only if there is a Cauchy sequence $p_{k}\in 
U,$ with $p_{k} \rightarrow  p$ in the metric completion $\bar U$ of 
$U$ and such that if $q_{\varepsilon}(k)\in U_{j}$ converges to $p_{k}$ 
as $\varepsilon  \rightarrow $ 0 then 
$\rho_{\varepsilon}(q_{\varepsilon}(k)) \leq  \mu (k),$ where $\mu (k) 
\rightarrow $ 0 as $k \rightarrow  \infty .$ Thus, $\partial U$ is 
typically, although not necessarily, singular; it is possible that $U$ 
may extend through (parts of) $\partial U$ to a larger smooth manifold.

\smallskip

 Since each $g_{\varepsilon}$ satisfies the Euler-Lagrange equations 
(1.14)-(1.15), the limit metric $g_{o}$ (of a subsequence) satisfies 
the same equation with $\varepsilon $ set to 0, i.e.
\begin{eqnarray}
L^{*}(u) + (\tfrac{1}{4}su + \tfrac{1}{12}\frac{\sigma (M)}{T})g & = & 
0, \\ 
2\Delta u + \tfrac{1}{4}(su - \frac{\sigma (M)}{T}) & = & 0. 
\end{eqnarray}
These are the Euler-Lagrange equations for $\mathcal{S}_{-}^{2}$ at the 
value $\sigma (M).$ When $\sigma (M) =$ 0, the equations (1.23)-(1.24) 
are the static vacuum Einstein equations with potential $u$: 
\begin{equation} \label{e1.25}
L^{*}u = 0, \ \ \Delta u = 0. 
\end{equation}
(Here one uses the fact from [2,(3.45)] that $s \rightarrow 0$ wherever 
$u \geq 0$).

 The following result characterizes these solutions at points 
$x_{\varepsilon}$ satisfying $u(x_{\varepsilon}) \rightarrow $ 1 as 
$\varepsilon  \rightarrow $ 0.

\begin{proposition} \label{p 1.2.}
  Let $x_{\varepsilon}$ be base points in $(\Omega_{\varepsilon}, 
g_{\varepsilon})$ such that $u(x_{\varepsilon}) \rightarrow $ 1, and 
suppose (1.21) holds. Then any limit (U, $g_{o},$ x) constructed above 
is hyperbolic, with constant scalar curvature $\sigma (M)$ if $\sigma 
(M) < $ 0, or flat if $\sigma (M) =$ 0.
\end{proposition}

\noindent
{\bf Proof:}
 If $\sigma (M) = 0$, then the maximum principle applied to trace 
equation in (1.25) implies that $u \equiv  1$ in $U$, and hence the 
full equation $L^{*}u = 0$ implies $g_{o}$ is flat, c.f. (1.5).

 The argument is essentially the same when $\sigma (M) < $ 0. Thus, 
(1.9) implies that in the limit one has $s = - u \cdot |\sigma (M)|$ 
and so, setting $\omega  = u- 1 \leq $ 0, the trace equation (1.24) can 
be written as
\begin{equation} \label{e1.26}
2\Delta\omega  -  \tfrac{1}{4}|\sigma (M)|(u+1)\omega  = 0. 
\end{equation}
Since $\omega (x) = 0 = \max \omega$, the strong maximum principle, 
c.f. [9, Thm. 3.5], implies that $\omega  \equiv $ 0 in $U$ and hence 
from (1.23), the metric $g_{o}$ is of constant negative curvature.
{\endproof}

 The situation where
\begin{equation} \label{e1.27}
\rho (x_{\varepsilon}) \geq  \rho_{o} >  0, 
\end{equation}
for all $x_{\varepsilon}\in (\Omega_{\varepsilon}, g_{\varepsilon})$ as 
$\varepsilon  \rightarrow $ 0, i.e. when the metrics $g_{\varepsilon}$ 
do not degenerate anywhere in $\Omega_{\varepsilon},$ is analysed in 
detail in [1]; in particular (1.27) holds as $\varepsilon  \rightarrow 
$ 0 if and only if $M$ is tame, i.e. (0.4) holds. As mentioned in \S 0, 
the Main Theorem has been proved in this case. 

\medskip

 Observe that Proposition 1.1 shows that the metrics $g_{\varepsilon}$ 
do not degenerate in the region $U^{\upsilon_{o}},$ provided one stays 
a fixed distance away from the boundary $L^{\upsilon_{o}}.$ However, 
while the estimates (1.9)-(1.11) give some useful information on the 
structure of $U^{\upsilon_{o}},$ for instance 
$vol_{g_{\varepsilon}}U^{\upsilon_{o}} \rightarrow $ 1 as $\varepsilon  
\rightarrow $ 0 in case $\sigma (M) < $ 0, it is possible apriori that 
the level $L^{\upsilon_{o}}$ is $\delta$-dense in 
$(\Omega_{\varepsilon}, g_{\varepsilon}),$ with $\delta  \rightarrow $ 
0 as $\varepsilon  \rightarrow $ 0.

 By far the most important issue is to understand the situation where 
\begin{equation} \label{e1.28}
\rho (x_{\varepsilon}) \rightarrow  0 \ \ {\rm as} \ \ \varepsilon  
\rightarrow  0, 
\end{equation}
for at least some $x_{\varepsilon}$ in the region where $u \geq  
u_{o},$ for some $u_{o} > $ 0, so that the curvature blows up at or 
very near $x_{\varepsilon}.$ We will not attempt to understand the 
global geometry of the degeneration in all of $U^{u_{o}},$ but it will 
be important to understand the local behavior of the degeneration, say 
in unit balls $B_{x_{\varepsilon}}(1) \subset (\Omega_{\varepsilon}, 
g_{\varepsilon})$.

\begin{remark} \label{r 1.3.}
  To give some perspective, here is a simple model situation. Suppose 
there are base points $x_{\varepsilon}$, $u(x_{\varepsilon}) 
\rightarrow $ 1 and satisfying (1.28), but which are isolated in the 
sense that there exists $d_{o} > $ 0 such that within 
$B_{x_{\varepsilon}}(d_{o})$, $\rho (y_{\varepsilon}) \geq  
dist_{g_{\varepsilon}}(y_{\varepsilon}, x_{\varepsilon}).$ Thus, if 
there is no collapse in $B(d_{o}) = B_{x_{\varepsilon}}(d_{o})$, and $u 
\rightarrow 1$ in $B(d_{o}/2)$, the metrics converge in $B(d_{o}/2)$ to 
a constant curvature metric $g_{o}$ outside $x = \lim x_{\varepsilon}$, 
but the curvature blows-up at $x_{\varepsilon}.$ The limit metric 
$g_{o}$ extends smoothly through $x$, giving a constant curvature 
metric on a 3-ball $B^{3}.$ In this situation, one has a natural 
2-sphere in $B(d_{o}/2)$ that surrounds $x_{\varepsilon},$ and it is 
easy to see, from the arguments in [1,\S 3], that this 2-sphere is 
essential, (at least on the "inside"). Unfortunately, it seems very 
difficult to prove that base points $x_{\varepsilon}$ with this good 
behavior necessarily exist.
\end{remark}

\medskip

 In the situation where (1.28) holds, we always consider the blow-up 
metrics
\begin{equation} \label{e1.29}
g_{\varepsilon}'  = \rho (x_{\varepsilon})^{-2}\cdot  g_{\varepsilon}, 
\end{equation}
so that $\rho' (x_{\varepsilon}) =$ 1, where $\rho' $ is the $L^{2}$ 
curvature radius w.r.t. $g_{\varepsilon}' .$ (The rescaling (1.28) is 
also used when $\rho (x_{\varepsilon}) \rightarrow  \infty ;$ this 
situation will arise only at the end of the paper, in \S 7.3).

 In the rescaling (1.29), the equations (1.14)-(1.15) become
\begin{equation} \label{e1.30}
L^{*}u  =  \bar \alpha\nabla\mathcal{Z}^{2} + \bar \phi \cdot  g,
\end{equation}
\begin{equation} \label{e1.31} 
2\Delta (u-\tfrac{\bar \alpha s}{12}) + \tfrac{1}{4}su = 
\tfrac{1}{2}\bar \alpha |z|^{2} -  3\bar c' . 
\end{equation}
All quantities in (1.30)-(1.31) are w.r.t. the $g_{\varepsilon}' $ 
metric. The constant $\bar \alpha = \bar \alpha(\varepsilon , 
x_{\varepsilon})$ is given by $\bar \alpha = \bar \varepsilon/\rho^{2}, 
\rho  = \rho (x_{\varepsilon})$, while $c'  = \rho^{2}c$. As emphasized 
in [2], the function $u$ is scale-invariant, i.e. it is considered as a 
function, and not as the normalized scalar curvature in the scale 
$g_{\varepsilon}'$.

\smallskip

 The following result from [2, Cor.5.2] implies there is an apriori 
maximal rate or scale at which the curvature can blow up on 
$(\Omega_{\varepsilon}, g_{\varepsilon})$ as $\varepsilon  \rightarrow 
$ 0.

\begin{proposition} \label{p 1.4.}
   There is a constant $\kappa  >$ 0, independent of $M$ and 
$(\Omega_{\varepsilon}, g_{\varepsilon}),$ such that
\begin{equation} \label{e1.32}
\rho^{2}(x_{\varepsilon}) \geq  \kappa\frac{\varepsilon}{T} = 
\kappa\bar \varepsilon, 
\end{equation}
for all $x_{\varepsilon}\in\Omega_{\varepsilon}.$ Equivalently, $\bar 
\alpha \leq  \kappa^{-1}.$
\end{proposition}

 The convergence/collapse possibilities, as well as Proposition 1.2, 
hold for the sequence $g_{\varepsilon}' $ exactly as in the case (1.21) 
holds; as before, the convergence is modulo diffeomorphisms but this 
will usually be ignored in the notation. The possible structure of 
blow-up limits, i.e. limits $(N, g', x)$ of the Riemannian manifolds 
$(\Omega_{\varepsilon}, g_{\varepsilon}' , x_{\varepsilon})$ has been 
analysed in detail in [2]; specifically from [2, \S 4,\S 5], one has:

\begin{proposition} \label{p 1.5.}
  If $u(x_{\varepsilon}) \geq  u_{o},$ for some constant $u_{o} > $ 0, 
then any blow-up limit (N, $g' ,$ x) of the metrics 
$(\Omega_{\varepsilon}, g_{\varepsilon}, x_{\varepsilon})$ is of the 
following form:

 (i) a complete non-flat solution of the $\mathcal{Z}_{c}^{2}$ 
equations, with uniformly bounded curvature,

 (ii) an incomplete non-flat solution of the static vacuum Einstein 
equations,

 (iii) an incomplete flat solution of the static vacuum Einstein 
equations.
\end{proposition}

\smallskip

 The static vacuum equations are given by (1.25) while the 
$\mathcal{Z}_{c}^{2}$ equations with potential $u$ are 
\begin{equation} \label{e1.33}
L^{*}u  =  \bar \alpha\nabla\mathcal{Z}^{2},
\end{equation}
\begin{equation} \label{e1.34}
2\Delta (u-\tfrac{\bar \alpha}{12}s) = \tfrac{1}{2}\bar \alpha|z|^{2}, 
\end{equation}
where $\bar \alpha$ is a positive constant, (the limit of $\bar 
\varepsilon/\rho^{2}).$ These equations are obtained from (1.30)-(1.31) 
by replacing $\bar \varepsilon = \varepsilon /T$ by $\bar \alpha,$ 
replacing $\bar \phi$ and $\bar c$ by 0 and setting $su = 0$. Thus, on 
a $\mathcal{Z}_{c}^{2}$ solution, the support of $s$ and of $u$ are 
disjoint in their interiors. Note that when $\bar \alpha =$ 0, the 
equations (1.33)-(1.34) become the static vacuum equations (1.25). 
Further, to leading order, the Euler-Lagrange equations (1.14)-(1.15) 
are identical to the $\mathcal{Z}_{c}^{2}$ equations (1.33)-(1.34).

\smallskip

 The three cases above are characterized as follows: Case (i) occurs 
exactly when $\bar \alpha > $ 0, Case (ii) occurs exactly when $\bar 
\alpha =$ 0 and $u$ is non-constant on the limit, while Case (iii) 
occurs exactly when $\bar \alpha =$ 0 and $u = const.$ on the limit. 

 The limit $(N, g', x)$ is a $\mathcal{Z}_{c}^{2}$ solution whenever 
$x_{\varepsilon}$ is a point of locally {\it  maximal}  curvature 
blow-up, i.e. $\rho (x_{\varepsilon}) \leq  \rho (y_{\varepsilon}),$ 
for all $y_{\varepsilon}\in B_{x_{\varepsilon}}(\delta ) \subset  
(\Omega_{\varepsilon}, g_{\varepsilon}),$ for a fixed $\delta  > $ 0. 
By [2, Thm.B], whenever $M$ is not tame, i.e. $\rho  \rightarrow $ 0 
somewhere on $(\Omega_{\varepsilon}, g_{\varepsilon})$ as $\varepsilon  
\rightarrow $ 0, there exist base points $x_{\varepsilon}$ whose 
blow-up limit is a $\mathcal{Z}_{c}^{2}$ solution. The limit is a 
possibly flat static vacuum solution when the curvature blows up at a 
rate much slower than the maximal rate. Observe that if the limit $(N, 
g' , x)$ based at $x_{\varepsilon}$ is flat, then at any base points 
$y_{\varepsilon}$ converging to points in $\partial N,$ the curvature 
blows up much faster than at $x_{\varepsilon},$ i.e. $\rho 
(y_{\varepsilon}) <<  \rho (x_{\varepsilon}).$

\smallskip

\begin{remark} \label{r 1.6.}
  Recall the Sphere conjecture from [1]: if $(N, g', y)$ is a complete 
$\mathcal{Z}_{c}^{2}$ solution arising as a blow-up limit of 
$(\Omega_{\varepsilon}, g_{\varepsilon}, y_{\varepsilon}),$ then $(N, 
g')$ has an asymptotically flat end. The Sphere conjecture implies 
Conjectures I-II.

 Any asymptotically flat end of $(N, g')$ of course carries natural 
2-spheres, namely the spheres $S^{2}$ approximating the spheres 
$S^{2}(R)$ of large radius in ${\mathbb R}^{3}.$ Since the convergence 
to the limit is smooth, such spheres also lie in 
$(\Omega_{\varepsilon}, g_{\varepsilon}' ),$ and hence also in $M$. 
Theorem C of [2] proves the Sphere conjecture in case the potential 
function $u$ on $(N, g')$ is bounded away from 0 outside a compact set, 
and the level sets of $u$ are compact. The essential remaining 
difficulty in proving the conjecture lies in dealing with the situation 
where the levels of $u$ are non-compact.

 In this paper, in place of attempting to locate natural 2-spheres in 
ends of a $\mathcal{Z}_{c}^{2}$ solution as in the Sphere conjecture 
above, we instead locate 2-spheres directly in {\it  flat}  blow-up 
limits, as in Case (iii) of Proposition 1.5. Of course these two issues 
may, (and in fact should be), closely related since if a 
$\mathcal{Z}_{c}^{2}$ solution has an asymptotically flat end, then by 
``blowing it down'', i.e. rescaling the metric by factors converging to 
0, one obtains a limit which is flat, with an isolated singular point; 
the 2-sphere near $\infty $ thus becomes a standard 2-sphere, of size 
on the order of 1, surrounding the origin in ${\mathbb R}^{3}.$ This 
blow-down of the limit $\mathcal{Z}_{c}^{2}$ solution can be obtained 
directly as a flat ``blow-up'' limit of $(\Omega_{\varepsilon}, 
g_{\varepsilon})$ at base points $x_{\varepsilon}$ nearby but distinct 
from $y_{\varepsilon}.$ Note that such a 2-sphere in 
$(\Omega_{\varepsilon}, g_{\varepsilon})$ is very small; it has size on 
the order of $\rho (x_{\varepsilon}) \rightarrow $ 0 and is being 
rescaled to approximately unit size in $g_{\varepsilon}' = \rho 
(x_{\varepsilon})^{-2}\cdot  g_{\varepsilon}$; compare with Remark 1.3.

 The main concern of the paper is to locate such 2-spheres by carefully 
analysing the structure of level sets of the potential $u$, i.e. 
studying the value distribution theory of $u$ on 
$(\Omega_{\varepsilon}, g_{\varepsilon}).$ In effect, the main point is 
to show that there exist base points $x_{\varepsilon}$ satisfying 
$u(x_{\varepsilon}) \rightarrow  1$ and (1.28) such that in rescalings 
(1.29) converging to a flat limit, the level sets of $u$ have at least 
one uniformly compact and separated component $\mathcal{C}$, c.f. 
Remark 2.2. It is then not difficult to show that such components are 
surrounded in a natural way by 2-spheres.
\end{remark}

\smallskip

 We conclude this section with the following elementary consequence of 
the $\sigma$-tame condition.

\begin{lemma} \label{l 1.7.}
  Suppose $M$ is $\sigma$-tame, so that there exist constants $k < 
\infty $ and $K < \infty$ such that, on the space of unit volume 
metrics on $M$,
\begin{equation} \label{e1.35}
\liminf_{\widetilde \sigma \rightarrow  0}\widetilde 
\sigma^{k}\mathcal{Z}^{2} \leq  K, 
\end{equation}
where $\widetilde \sigma$ is as in (0.5). Then there is a sequence 
$\varepsilon  = \varepsilon_{i} \rightarrow $ 0, such that on any 
minimizing pair $(\Omega_{\varepsilon}, g_{\varepsilon}), \varepsilon  
= \varepsilon_{i}$,
\begin{equation} \label{e1.36}
 \varepsilon\mathcal{Z}^{2} \leq  (1+K)\varepsilon^{1/(k+1)}, \ \ {\rm 
and} \ \ \widetilde \sigma \leq  (1+K)\varepsilon^{1/(k+1)}. 
\end{equation}
\end{lemma}

\noindent
{\bf Proof:}
 For any given $\varepsilon$ the metric $g_{\varepsilon}$ is a 
minimizer of $I_{\varepsilon}^{~-},$ so that
\begin{equation} \label{e1.37}
\varepsilon\mathcal{Z}^{2}(g_{\varepsilon}) + \widetilde 
\sigma(g_{\varepsilon}) \leq  \varepsilon\mathcal{Z}^{2}(g' ) + 
\widetilde \sigma(g' ), 
\end{equation}
for any unit volume metric $g'$ on $M$. By (1.35), there is a sequence 
of metrics $g_{i}$ on $M$ with $\widetilde \sigma_{i} = \widetilde 
\sigma(g_{i}) \rightarrow $ 0 and
$$\mathcal{Z}^{2}(g_{i}) \leq  K\widetilde \sigma_{i}^{-k}. $$
Substituting this in (1.37) gives
$$\varepsilon\mathcal{Z}^{2}(g_{\varepsilon}) + \widetilde 
\sigma(g_{\varepsilon}) \leq  \varepsilon K\widetilde \sigma_{i}^{-k} + 
\widetilde \sigma_{i}. $$
Now this estimate holds for any $\varepsilon .$ Choose $\varepsilon  = 
\varepsilon_{i}$ by setting
$$\varepsilon_{i} = \widetilde \sigma_{i}^{k+1}, $$
so that $\varepsilon_{i} \rightarrow $ 0 as $i \rightarrow  \infty .$ 
This gives
$$\varepsilon_{i}\mathcal{Z}^{2}(g_{\varepsilon_{i}}) + \widetilde 
\sigma(g_{\varepsilon_{i}}) \leq  K\varepsilon_{i}\cdot 
\varepsilon_{i}^{-k/(k+1)} + \varepsilon_{i}^{1/(k+1)}\leq  
\varepsilon_{i}^{1/(k+1)}(1+K), $$
which gives the result.
{\endproof}

 The estimate (1.36) is the only consequence of the $\sigma$-tame 
condition that will be used in the paper. Observe that (1.36) gives an 
estimate for the constant term $\bar c$ in (1.15); this will be used to 
obtain the global control in Theorem 3.5 and Theorem 5.4.

\smallskip

 Finally, recall the following formulas for the change of curvature 
under conformal deformation of the metric, c.f. [6, 1.159ff]. Thus, if 
$\widetilde g = f^{2}\cdot  g,$ then the Ricci curvature $\widetilde r$ 
of $\widetilde g$ is given by
\begin{equation} \label{e1.38}
\widetilde r = r -  f^{-1}D^{2}f + 2(d\ln f)^{2} - f^{-1}\Delta f\cdot  
g. 
\end{equation}
Similarly, if $\widetilde g = \psi^{4}\cdot  g,$ then
\begin{equation} \label{e1.39}
\psi^{5}\widetilde s = - 8\Delta\psi  + s\psi . 
\end{equation}

 It should also be understood that throughout the paper, we always pass 
to subsequences where convergence issues are concerned, without 
explicitly mentioning this at each use. A sequence is said to 
sub-converge if a subsequence converges.

\section{The Linearized Equations.}
\setcounter{equation}{0}

 In this section, we study in detail the linearization of the metrics 
$g_{\varepsilon}$ and their rescalings $g_{\varepsilon}' $ in (1.29) 
about their constant curvature limits given by Proposition 1.2. This 
forms the foundation for the work to follow in later sections. 
Throughout the paper, we study the behavior of $(\Omega_{\varepsilon}, 
g_{\varepsilon}, x_{\varepsilon})$ at base points $x_{\varepsilon}$ 
satisfying
\begin{equation} \label{e2.1}
u_{\varepsilon}(x_{\varepsilon}) \rightarrow  1 ,
\end{equation}
on some arbitrary sequence $\varepsilon  = \varepsilon_{i} \rightarrow 
$ 0. From \S 4 on, the sequence $\varepsilon_{i}$ will no longer be 
arbitrary, but will be required to satisfy (1.36).

 A general introductory discussion on such linearizations is given in 
\S 2.1, while in \S 2.2 the three possible forms of the linearized 
equations and the associated form of the metrics are analysed in 
detail. These are then summarized and unified in \S 2.3, c.f. Theorem 
2.11.

\smallskip

{\bf 2.1.}
 The metrics $g_{\varepsilon}$ satisfy the Euler-Lagrange equations 
(1.14)-(1.15). After a little algebra on the trace equation, these may 
be rewritten in the form
\begin{equation} \label{e2.2}
L^{*}u  =  \bar \varepsilon\nabla\mathcal{Z}^{2} + \bar \phi\cdot  g,
\end{equation}
\begin{equation} \label{e2.3}
2\Delta (u - \tfrac{\bar \varepsilon s}{12}) = \psi_{\varepsilon} + 
\tfrac{1}{2}\bar \varepsilon|z|^{2}, 
\end{equation}
where
\begin{equation} \label{e2.4}
\psi_{\varepsilon} = \tfrac{1}{4}\sigma T(u+1)(u- 1) + \bar 
d_{\varepsilon}, 
\end{equation}
and $\bar d_{\varepsilon}$ is a constant, depending on the {\it  
global}  geometry of $(\Omega_{\varepsilon}, g_{\varepsilon}),$ given by
\begin{equation} \label{e2.5}
\bar d_{\varepsilon} = \tfrac{1}{4}\sigma T(1- T^{-2}) -  
\tfrac{1}{2}\bar \varepsilon\mathcal{Z}^{2}. 
\end{equation}
All quantities in (2.2)-(2.5) are w.r.t. $(\Omega_{\varepsilon}, 
g_{\varepsilon})$ and $\bar \varepsilon\mathcal{Z}^{2} \equiv \bar 
\varepsilon\mathcal{Z}^{2}(g_{\varepsilon})$. Observe that (1.8) 
implies that $\bar d_{\varepsilon} \leq $ 0, so that both summands of 
$\psi_{\varepsilon}$ are non-positive. The term $\sigma T \rightarrow  
|\sigma (M)|$ when $\sigma (M) < $ 0, while $\sigma T \rightarrow $ 0 
when $\sigma (M) =$ 0, again by (1.8). Thus, $\psi_{\varepsilon}$ is 
also uniformly bounded below as $\varepsilon  \rightarrow $ 0, i.e. 
there exists $K <  \infty $ such that
\begin{equation} \label{e2.6}
- K \leq  \psi_{\varepsilon} \leq  0. 
\end{equation}
Further, in regions where $u \rightarrow $ 1, one has
\begin{equation} \label{e2.7}
\psi_{\varepsilon} \rightarrow  0. 
\end{equation}

 Apriori, the $L^{2}$ curvature radius $\rho (x_{\varepsilon}) = 
\rho_{\varepsilon}(x_{\varepsilon})$ may converge to 0, remain bounded 
away from 0 and $\infty ,$ or diverge to $\infty ,$ (as usual in 
subsequences), as $\varepsilon  \rightarrow $ 0. Regardless of which of 
these behaviors occurs, we will always work in the scale 
$g_{\varepsilon}' $ given by
\begin{equation} \label{e2.8}
g_{\varepsilon}'  = \rho (x_{\varepsilon})^{-2}\cdot  g_{\varepsilon}, 
\end{equation}
so that $\rho' (x_{\varepsilon}) =$ 1. The metrics $g_{\varepsilon}' $ 
then satisfy the rescaled Euler-Lagrange equations (1.30)-(1.31).

 As discussed in \S 1, one may then pass, (in subsequences), to a 
maximal connected (blow-up) limit $(F, g_{o}', x_{o})$ of the metrics 
$(\Omega_{\varepsilon}, g_{\varepsilon}' )$ based at 
$\{x_{\varepsilon}\},$ passing to sufficiently large finite covers if 
the sequence collapses at $x_{\varepsilon}.$ As noted following 
Proposition 1.4, the assumption (2.1) and Proposition 1.2, (in the 
scale (2.8)), imply that the maximal limit $(F, g_{o}', x_{o})$ at 
$x_{o} = \lim x_{\varepsilon}$ is a constant curvature manifold, i.e. 
either flat or hyperbolic, with limit potential $u \equiv $ 1, provided 
$\rho (x_{\varepsilon}) \leq  R_{o},$ for some $R_{o} <  \infty .$ If 
$\rho (x_{\varepsilon}) \rightarrow  \infty ,$ so that 
$g_{\varepsilon}$ is being 'blown-down', this remains true but does not 
follow directly from Proposition 1.2; since the proof is out of place 
here, it is given in Appendix A. (The very special situation when $\rho 
(x_{\varepsilon}) \rightarrow  \infty $ will only arise in \S 7.3 and 
so may be ignored until then). Note also that $F$ so defined is not 
apriori maximal w.r.t. a smooth flat or hyperbolic structure, i.e. $F$ 
may possibly extend to a larger smooth flat or hyperbolic manifold, 
(c.f. the discussion preceding (1.23)).

 It is of course apriori possible that the metrics 
$(\Omega_{\varepsilon}, g_{\varepsilon})$ do not degenerate at all, 
i.e. $\rho (x_{\varepsilon}) \geq  \rho_{o} > $ 0 whenever 
$u(x_{\varepsilon})$ is sufficiently close to 1. The question of 
whether the metrics $g_{\varepsilon}$ degenerate in this region plays 
an important role in \S 4, but is not used until then. The main case of 
interest is when the metrics degenerate, so that $\rho 
(x_{\varepsilon}) \rightarrow $ 0. 
\medskip

 To study the linearization of the metrics $g_{\varepsilon}' $ at the 
limit constant curvature metric $(F, g_{o}'),$ write
\begin{equation} \label{e2.9}
g_{\varepsilon}'  = g_{o}'  + \delta\cdot  h + o(\delta ), 
\end{equation}
where the parameter $\delta  = \delta (\varepsilon) \rightarrow  0$ as 
$\varepsilon  = \varepsilon_{i} \rightarrow  0$ and $h$ is a symmetric 
bilinear form, formally corresponding to 
$\frac{dg_{\varepsilon}'}{d\delta}|_{\delta =0}$. The choices of the 
parameter $\delta $ and the form $h$ of course depend on each other. 
While the sequence $\{g_{\varepsilon}'\}, \varepsilon  = 
\varepsilon_{i},$ converges smoothly to the limit $g_{o}'$, it may not 
converge smoothly with respect to $\varepsilon$, i.e. 
$dg_{\varepsilon}' /d\varepsilon $ may not exist at $\varepsilon =  0$, 
so that one cannot necessarily choose $\delta  = \varepsilon$. (Of 
course, this is related to the fact that the choice of base points 
$x_{\varepsilon}$ is far from unique). 

 The natural value for the parameter $\delta $ is the local size of the 
curvature near the base point $x_{\varepsilon},$ which gives a measure 
as to how far away $g_{\varepsilon}' $ is from the limit metric $g_{o}' 
.$ For reasons which will only become clear later, we use for this the 
trace-free Ricci curvature $z$ in place of the full Ricci curvature 
$r$. Thus, set
\begin{equation} \label{e2.10}
\delta   = \delta_{z} = (\oint_{B_{x_{\varepsilon}}' 
(\frac{1}{2})}|z_{\varepsilon}'|^{2}dV_{g_{\varepsilon}'})^{1/2} = 
(\rho^{4}(x_{\varepsilon})  \oint_{B_{x_{\varepsilon}}(\frac{1}{2}\rho 
(x_{\varepsilon}))}|z_{\varepsilon}|^{2}dV_{g_{\varepsilon}})^{1/2} , 
\end{equation}
where $\oint$ denotes the average value, i.e. the integral divided by 
the volume of the domain. Observe that the right side of (2.10) is 
scale invariant, so that $\delta $ is scale invariant. The quantity 
$\delta $ depends only on the isometry class of $g_{\varepsilon}$ and 
so is independent of diffeomorphisms 'reparametrizing' 
$g_{\varepsilon}.$ Clearly $\delta  = \delta (x_{\varepsilon})$ depends 
strongly on the choice of the base point. The factor $\frac{1}{2}$ in 
(2.10) could be replaced by any other fixed constant $< $ 1. The 
quantities $\delta_{r}$ and $\delta_{s}$ are defined in the same way as 
$\delta ,$ by replacing $z$ by the Ricci curvature $r$ or scalar 
curvature $s$ respectively. 

 Note that from the definition of $\rho $ and $\delta ,$
$$\delta  \leq  1 $$
always. Further, since the convergence to the limit is smooth in 
$(B_{x}(1), g_{o}' )$ and the limit $g_{o}' $ is of constant curvature 
since (2.1) holds, one has
\begin{equation} \label{e2.11}
\delta (x_{\varepsilon}) \rightarrow  0, \ \ {\rm as} \ \  \varepsilon  
\rightarrow  0. 
\end{equation}
On the other hand, only $x_{\varepsilon}$ such that $\delta 
(x_{\varepsilon}) >  0$ for $\varepsilon$ sufficiently small need to be 
considered. For if $\delta (x_{\varepsilon}) =$ 0 on some sequence 
$\varepsilon  = \varepsilon_{i} \rightarrow $ 0, then $g_{\varepsilon}$ 
is of constant curvature, (hyperbolic or flat), and $u =$ const, on 
$B_{x_{\varepsilon}}(\rho (x_{\varepsilon})).$ Since $g_{\varepsilon}$ 
is real-analytic in the region where $u > $ 0, $g_{\varepsilon}$ is of 
constant curvature on the full component $C_{\varepsilon}$ of 
$\Omega_{\varepsilon}$ containing $x_{\varepsilon}.$ In this case, the 
metrics $(C_{\varepsilon}, g_{\varepsilon}) = (C, g_{o})$ are constant, 
i.e. independent of $\varepsilon ,$ for all $\varepsilon ,$ and no 
further analysis is necessary, c.f. also \S 1.

 There are three basic cases for the behavior of $\delta $ at base 
points $x_{\varepsilon},$ as $\varepsilon  = \varepsilon_{i} 
\rightarrow $ 0:
\begin{eqnarray}
& (i)   & \delta  >>  \rho^{2}, \nonumber \\
& (ii)  & \delta  \sim  \rho^{2}, \\
& (iii) & \delta  <<  \rho^{2}. \nonumber
\end{eqnarray}

By definition (2.10), $\delta /\rho^{2}$ is just the $L^{2}$ average of 
$z$ over $B(\frac{1}{2}\rho ).$ Observe also that $s_{g_{\varepsilon}'} 
= \rho^{2}s_{g_{\varepsilon}} = O(\rho^{2}),$ so that in all cases 
\begin{equation} \label{e2.13}
\delta_{s} \leq  c\cdot \rho^{2}. 
\end{equation}

 It is obvious from (2.11) that $\rho (x_{\varepsilon}) \rightarrow $ 0 
in Cases (i) and (ii). In Case (iii), either of the behaviors $\rho 
(x_{\varepsilon}) \rightarrow $ 0 or $\rho (x_{\varepsilon}) \geq  
\rho_{o} > $ 0 is possible. Thus, in Cases (i) and (ii), the maximal 
blow-up limit $(F, g_{o}', x)$ based at $\{x_{\varepsilon}\}$ with $x = 
\lim x_{\varepsilon}$, is flat, while in Case (iii), it is either flat 
or hyperbolic. 

 For $\delta $ as in (2.10), the local $L^{2}$ norm of $h$ in (2.9) is 
on the order of 1 in the Cases (i) and (ii) above. More precisely, 
define
\begin{equation} \label{e2.14}
h_{\varepsilon} = \frac{g_{\varepsilon}'  -  g_{o}'}{\delta}. 
\end{equation}
Then
\begin{equation} \label{e2.15}
\oint_{B_{x_{\varepsilon}}(\frac{1}{2})}|h_{\varepsilon}|^{2}
dV_{g_{\varepsilon}'} \sim  1 \ \ {\rm as} \ \  \delta  \rightarrow  0. 
\end{equation}
This is because the $L^{2}$ norm of the curvature $r_{\varepsilon}' $ 
of $g_{\varepsilon}' $ is bounded by $\delta $ in 
$(B_{x_{\varepsilon}}' (\frac{1}{2}), g_{\varepsilon}' )$ and the 
components of $g_{\varepsilon}' - g_{o}'$ in a harmonic coordinate 
chart on these balls are bounded in the $L^{2,2}$ topology by the 
$L^{2}$ norm of the curvature. Moreover, since $g_{\varepsilon}'$ 
satisfies the elliptic system (1.30)-(1.31), the covariant derivatives 
$\nabla^{k}z_{g_{\varepsilon}'}/\delta$, $k \geq  0$, are uniformly 
bounded in $L^{2}$, while $\alpha \nabla^{k}z_{g_{\varepsilon}'}/\delta 
\rightarrow 0$ as $\delta \rightarrow 0$, both on compact subsets of 
$(F, g_{o}')$. The proof of this (technical) result is given in 
Appendix B.

 The estimate (2.15) may not hold in Case (iii), since one might have 
$\delta_{z} <<  \delta_{s},$ so that $s_{g_{\varepsilon}}' /\delta $ 
may be unbounded as $\delta  \rightarrow $ 0. This case must then be 
treated somewhat more carefully.

 In analogy to (2.9), let
\begin{equation} \label{e2.16}
u_{\varepsilon} = 1 + \nu\cdot \delta  + o(\delta ), 
\end{equation}
where $\nu $ is (formally) $du_{\varepsilon}/d\delta $ at $\delta  =$ 
0. More precisely, as in (2.14), set
\begin{equation} \label{e2.17}
\nu_{\delta} = \frac{u_{\varepsilon}- 1}{\delta}, 
\end{equation}
and let $\nu  = \lim \nu_{\delta}$ if the limit exists. Special care 
must be taken here, since apriori $\nu $ might be infinite. 

 We complete this subsection with the following characterization of the 
boundary $\partial F$.
\begin{proposition} \label{p 2.1.}
  For any maximal limit $(F, g_{o}', x)$ with $x_{\varepsilon} 
\rightarrow  x$ satisfying (2.1), the boundary $\partial F$ is the 
Gromov-Hausdorff limit of the level set $L^{\upsilon_{o}}$ from (1.17).
\end{proposition}

\noindent
{\bf Proof:}
 Recall that $\partial F$ is formed by limits of sequences where $\rho' 
 \rightarrow $ 0. By the scale-invariant estimate (1.19), $\rho'  
\rightarrow $ 0 implies $t_{\upsilon_{o}}'  \rightarrow $ 0, where 
$t_{\upsilon_{o}}'  = dist_{g_{\varepsilon}'}(L^{\upsilon_{o}}, \cdot  
).$ Hence $\partial F$ is contained in the Gromov-Hausdorff limit of 
$L^{\upsilon_{o}}$, c.f. [10,Ch.3], or [14,Ch.10] for discussion of 
Gromov-Hausdorff limits.

 To prove the converse, it suffices to prove that the converse of 
Proposition 1.1 also holds, i.e. there is a constant $R_{o} < \infty$, 
independent of $\varepsilon ,$ such that
\begin{equation} \label{e2.18}
\rho (x_{\varepsilon}) \leq  R_{o}\cdot  
t_{\upsilon_{o}}(x_{\varepsilon}). 
\end{equation}
To see this, the estimate (2.18) is scale-invariant, and so it suffices 
to prove it in the scale $g_{\varepsilon}' $ where $\rho' 
(x_{\varepsilon}) =$ 1. The metrics $(\Omega_{\varepsilon}, 
g_{\varepsilon}' , x_{\varepsilon})$ converge smoothly to the maximal 
limit $(F, g_{o}', x)$, and the potential $u$ converges smoothly on 
compact subsets of $F$ to the limit function $u_{o}.$ Proposition 1.2, 
or more precisely its proof, shows that $u_{o} \equiv $ 1. It follows 
that the level set $L^{\upsilon_{o}}$ cannot converge to any interior 
region of $F$, and hence $t_{\upsilon_{o}}'(x_{\varepsilon})$ is at 
least 1 in the limit, which gives (2.18).

{\endproof}

  Note that for any flat manifold $(F, g_{o}')$, the curvature radius 
$\rho'$ satisfies
\begin{equation} \label{e2.19}
\rho'  = t', 
\end{equation}
where $t'(x) = dist_{F}(x, \partial F)$; in (1.16) the balls 
$B_{x}(\rho'(x))$ are required to be contained in $F$. Similarly, if 
the limit is hyperbolic, i.e. of constant negative curvature, then 
$\rho'  \leq  t'$. Note that equality, i.e. (2.19), holds when 
$(s')^{2} \leq  3c_{o},$ where $s'$ is the (constant) scalar curvature 
of $g_{o}'$ and $c_{o}$ is as in (1.16).

\medskip

 The structure of $\partial F$ plays an important role throughout the 
paper. The metric closure $\bar F = F\cup\partial F$ is a length space, 
c.f. [10, Ch.1], and is viewed as a singular flat (or hyperbolic) 
manifold. In situations where $F$ is embedded in ${\mathbb R}^{3},$ (or 
a quotient of ${\mathbb R}^{3}), \bar F$ is then just the closure of a 
domain in ${\mathbb R}^{3},$ (or the quotient).

 Some of the more technical issues in \S 3 and \S 4 arise because unit 
balls in $\bar F$ are not necessarily compact. Consider for instance 
$F$ of the form of a product $F = {\mathbb R}\times V,$ where $V$ is 
the universal cover of ${\mathbb R}^{2}\setminus \{0\}$; $V$ may 
naturally be identified with the upper half plane. In this setting, the 
topological boundary of ${\mathbb R}\times V$ is ${\mathbb R}^{2},$ 
while the metric boundary is a line, i.e. ${\mathbb R}\times \{0\}$. 
Hence, $\partial F,$ or even better $\partial V =$ \{0\}, acts as a 
"wormhole", in that points of $V$ that are far apart, when measured by 
the length of curves a bounded distance away from $\partial V,$ may be 
close when the curves are allowed to be arbitrarily close to or within 
$\partial V.$ On the other hand, it will be seen later following 
Proposition 4.8 that this is essentially the only way that $\partial F$ 
can act as a wormhole.

\begin{remark} \label{r 2.2.}
  Much of the work of the paper through \S 5 can be characterized as 
follows: find a flat blow-up limit $(F, g', y)$ of 
$(\Omega_{\varepsilon}, g_{\varepsilon})$ such that $\partial F$ is 
compact with $F$ non-compact, or at least there is a ball $B_{x}(R)$ of 
finite radius $R$ with center $x\in\bar F,$ and a constant $d_{o} > $ 0 
such that
\begin{equation} \label{e2.20}
B_{x}(R)\cap\partial F \neq  \emptyset  , \ \ {\rm but} \ \  A_{x}(R, 
(1+d_{o})R) \cap  \partial F = \emptyset  , 
\end{equation}
where $A_{x}(r,s)$ is the geodesic annulus about $x$ of radii $r, s$. 
Thus $\partial F$ has a compact region, separated by a definite amount 
from any other part of $\partial F$. The relation (2.20) is established 
in Theorem 5.4. It is in such gap regions that geometrically natural 
2-spheres in $(\Omega_{\varepsilon}, g_{\varepsilon})$ will be located.
\end{remark}

\medskip

{\bf  2.2.}
 The main results of this subsection, (Propositions 2.3-2.4 and 
Corollary 2.8), describe the form of the linearizations of 
$g_{\varepsilon}' $ and $u$ at their limits $(g_{o}', 1)$. For clarity, 
the results are separated into three cases according to the three 
possibilities in (2.12). These three distinct situations will then be 
unified in \S 2.3, c.f. Theorem 2.11. An analysis similar to the 
analysis here was previously carried out in [1, \S 4], and also in [2, 
\S 7] in the context of $\mathcal{Z}_{c}^{2}$ solutions; this may serve 
as an introduction to the discussion below.

 The first, and most important case, is the following.

\begin{proposition} \label{p 2.3.}
  Suppose the base points $x_{\varepsilon}$ satisfy (2.1) and suppose 
that
\begin{equation} \label{e2.21}
\delta  >>  \rho^{2}(x_{\varepsilon}), \ \ {\rm as} \ \ \varepsilon  = 
\varepsilon_{i} \rightarrow  0. 
\end{equation}
Then there is a sequence of affine functions $a_{\delta}$ such that 
$\nu_{\delta}- a_{\delta}$ subconverges to a smooth limit function $\nu 
$ on $F$. The linearizations $h = \lim_{\varepsilon\rightarrow 
0}h_{\varepsilon}$ and $\nu $ are a non-trivial solution of the 
linearized static vacuum Einstein equations 
\begin{equation} \label{e2.22}
r'  = \frac{dr}{d\delta}|_{\delta =0} = D^{2}\nu , \ \ \Delta\nu  = 0, 
\end{equation}
at the flat metric $g_{o}'$. Thus, the limit function $\nu$ is a 
non-affine, locally bounded and smooth harmonic function on $(F, 
g_{o}')$. 

 To $1^{st}$ order in $\delta ,$ the metric $g_{\varepsilon}' $ is 
conformally flat. In fact, for $\varepsilon $ small, and on any compact 
domain $K \subset  F$, the metric $g_{\varepsilon}' $ has the 
expansion, modulo diffeomorphisms,
\begin{equation} \label{e2.23}
g_{\varepsilon}'  = (1- 2\nu\delta )(g_{o}' +\delta\chi ) + o(\delta ), 
\end{equation}
where $\delta  \rightarrow $ 0 as $\varepsilon  \rightarrow $ 0, 
$|o(\delta )| <<  \delta $, and $\chi $ is an infinitesimal holonomy 
deformation, i.e. $g_{o}' +\delta\chi $ is a curve (in $\delta$) of 
flat metrics on $F$.
\end{proposition}

\noindent
{\bf Proof:}
 The metrics $g_{\varepsilon}' $ satisfy the Euler-Lagrange equations 
(2.2)-(2.3) in this scale, i.e.
\begin{equation} \label{e2.24}
 L^{*}u  =  \bar \alpha\nabla\mathcal{Z}^{2} + (-\tfrac{1}{4}su+\bar 
c_{\varepsilon})\cdot  g
\end{equation}
$$2\Delta (1+\tfrac{\varepsilon\sigma}{12})u = \psi_{\varepsilon}' + 
\tfrac{\bar \alpha}{2}|z|^{2}.$$
Here and below, the prime and subscript $\varepsilon $ are usually 
omitted from the notation, and $u > 0$. The term $\psi_{\varepsilon}' $ 
is given by $\psi_{\varepsilon}'  = \rho^{2}\cdot \psi_{\varepsilon}$, 
$\rho  = \rho (x_{\varepsilon})$. To obtain the linearization, divide 
the equations (2.24) by $\delta ,$ and consider the limit as $\delta  
\rightarrow  0$, (as always in subsequences).

 By Proposition 1.5 and the discussion following it, $\bar \alpha = 
\varepsilon /\rho^{2}(x_{\varepsilon}) \rightarrow $ 0 as $\varepsilon  
\rightarrow $ 0. Consider first the linearization of the $\bar 
\alpha|z|^{2}$ term, and write this in the form:
$$\frac{\bar \alpha|z|^{2}}{\delta} = \bar \alpha|z|\frac{|z|}{\delta}. 
$$
Proposition 1.1 implies that $|z| \rightarrow 0$ in $L^{\infty}$ on any 
compact subset of $F$, while $\bar \alpha|z|/\delta $ is bounded, by 
the discussion following (2.15), c.f. also Appendix B. It follows that, 
as $\delta  \rightarrow  0$,
\begin{equation} \label{e2.25}
\frac{\bar \alpha|z|^{2}}{\delta} = o(1), 
\end{equation}
on compact subsets of $F$. The same reasoning, again from the 
discussion following (2.15), implies
\begin{equation} \label{e2.26}
(\bar \alpha\nabla\mathcal{Z}^{2})/\delta  = o(1). 
\end{equation}

 Next consider the remaining terms on the right in (2.24). By definition
$$\frac{\psi_{\varepsilon}'}{\delta}= 
\frac{\rho^{2}}{\delta}\psi_{\varepsilon}. $$
From (2.7), $\psi_{\varepsilon} \rightarrow 0$ as $\varepsilon 
\rightarrow 0$. Since, by assumption, $\delta  >>  \rho^{2}$, this term 
tends to $0$ uniformly on compact subsets of $F$. The same analysis 
holds for the terms $su$ and $\bar c_{\varepsilon}$; note that $s = 
g_{g_{\varepsilon}'} = \rho^{2}s_{g_{\varepsilon}}$, and similarly for 
the term $\bar c_{\varepsilon}$ in this scale.

  Hence, these arguments imply that
$$(L^{*}u)/\delta  = o(1) \ \ {\rm and} \ \ (\Delta u/\delta ) = o(1), 
$$
in $L^{\infty}$ on compact subsets of $F$. In particular, by (1.5),
\begin{equation} \label{e2.27}
u\frac{r}{\delta} = D^{2}\nu_{\delta} + o(1), \ \  \Delta\nu_{\delta} = 
0 + o(1). 
\end{equation}

 The equations (2.27), with the $o(1)$ and $\delta $ terms removed, are 
the static vacuum Einstein equations (1.25). Since $u \rightarrow  1$, 
it follows that any limit $(h, \nu )$ of $(h_{\delta}, \nu_{\delta})$ 
is a solution of the linearization of the static vacuum equations 
(2.22) at the flat limit $g_{o}'$, with $u \equiv 1$. To prove the 
existence of such limits, the regularity discussion following (2.15), 
(c.f. Appendix B), implies that $r/\delta $ converges, in a 
subsequence, to a limit $r' $ and further $h_{\varepsilon}$ converges, 
modulo diffeomorphisms, to a limit $h$ on $(F, g_{o}')$. 

 Thus, (2.27) implies that the functions $\nu_{\delta}$ are bounded, 
modulo addition of functions in the kernel of $D^{2}.$ Since Ker 
$D^{2}$ consists of affine functions on $F$, there are affine functions 
$a_{\delta}$ such that $\nu_{\delta}- a_{\delta}$ converges to a limit 
harmonic function $\nu $ on $(F, g_{o}')$. This proves the first 
statement.

 To prove the second statement, write
\begin{equation} \label{e2.28}
\frac{u- 1}{\delta} = \nu  + a_{\delta} + o(1), 
\end{equation}
and set $\upsilon  = u- a_{\delta}\cdot \delta  = 1+\nu\delta  + 
o(\delta ).$ Clearly $a_{\delta}\cdot \delta  \rightarrow $ 0 as 
$\delta  \rightarrow $ 0. Now consider the conformally equivalent 
metric 
\begin{equation} \label{e2.29}
\widetilde g_{\varepsilon} = \upsilon^{2}g_{\varepsilon}' , 
\end{equation}
A simple computation using (2.27) and (1.38) shows that
\begin{equation} \label{e2.30}
\widetilde r = 2(d\ln \upsilon )^{2} -  
\frac{2}{\upsilon}\Delta\upsilon\cdot  g -  \frac{\bar 
\alpha}{\upsilon}\nabla\mathcal{Z}^{2} , 
\end{equation}
where all metric quantities on the right in (2.30) are w.r.t. $g = 
g_{\varepsilon}' .$ Hence,
$$\frac{\widetilde r}{\delta} = \frac{2(d\ln \upsilon )^{2}}{\delta} -  
\frac{2\Delta\upsilon}{\upsilon\delta}\cdot  g -  \frac{\bar 
\alpha\nabla\mathcal{Z}^{2}}{\upsilon\delta} . $$
By the arguments in (2.25)-(2.27), the last two terms on the right go 
to 0 as $\delta  \rightarrow $ 0. Similarly, since $d\ln \upsilon  
\rightarrow  0$ as $\delta  \rightarrow  0$, and $(d\ln \upsilon 
)/\delta  = \frac{1}{\upsilon}d\nu_{\delta}+o(1) \rightarrow  d\nu$ as 
$\delta  \rightarrow  0$, it follows that
\begin{equation} \label{e2.31}
\frac{\widetilde r}{\delta} \rightarrow  0, \ \ {\rm as} \ \ \delta  
\rightarrow  0. 
\end{equation}
Thus, to first order in $\delta$, $\widetilde g_{\varepsilon}$ is flat 
and so $g_{\varepsilon}' $ is conformally flat. More precisely, one may 
write
\begin{equation} \label{e2.32}
\widetilde g_{\varepsilon} =  \widetilde g_{o} + \delta\chi  + o(\delta 
) = g_{o}'  + \delta\chi  + o(\delta ), 
\end{equation}
where $\chi $ is an infinitesimal flat variation of the metric $g_{o}' 
.$ The variation $\chi $ is a sum of a trivial term, due to variations 
of $\widetilde g_{\varepsilon}$ by diffeomorphisms, and a possibly 
non-trivial term due to variations in the holonomy, i.e. periods, of 
$(F, g_{o}')$. Ignoring the contribution due to diffeomorphisms, it 
follows that
$$g_{\varepsilon}'  = \upsilon^{-2}(g_{o}' +\delta\chi ) + o(\delta ) = 
(1- 2\nu\cdot \delta )(g_{o}'  + \delta\chi ) + o(\delta ), $$
which proves the result.
{\endproof}

 As noted above in the proof, the affine indeterminacy of $\nu $ arises 
from the fact that Ker $D^{2}$ consists of affine functions. 

 Next we turn to the Case (ii) of (2.12).
\begin{proposition} \label{p 2.4.}
  Suppose the base points $x_{\varepsilon}$ satisfy (2.1) and suppose 
that
\begin{equation} \label{e2.33}
\delta  \sim  \rho^{2}(x_{\varepsilon}), \ \ {\rm as} \ \ \varepsilon  
\rightarrow  0. 
\end{equation}
Then there exists a sequence of affine functions $a_{\delta}$ such that 
the functions $\nu_{\delta} - a_{\delta}$ converge to a limit function 
$\nu$. The linearizations $h_{\varepsilon}$ from (2.14) converge to a 
solution of the linearized $\mathcal{S}_{-}^{2}$ equations,
\begin{equation} \label{e2.34}
z'  = D^{2}\nu , \ \ \Delta\nu  = 0, 
\end{equation}
on $(F, g_{o})$. This solution is non-trivial in the sense that 
$D^{2}\nu  \neq $ 0. On compact subsets of $F$, the metric 
$g_{\varepsilon}' $ has the expansion, modulo diffeomorphisms,
\begin{equation} \label{e2.35}
g_{\varepsilon}'  = (1- 2\nu\cdot \delta )(g_{-\delta}' +\delta\chi ) + 
o(\delta ), 
\end{equation}
where $g_{-\delta}' $ is the space form of constant curvature $(\sigma 
(M)r_{o}/6)\delta$, $r_{o} = \lim_{\varepsilon\rightarrow 
0}(\rho^{2}/\delta )(x_{\varepsilon})$ and $\chi$ is an infinitesimal 
holonomy deformation of $(F, g_{-\delta}')$.
\end{proposition}

\noindent
{\bf Proof:}
Following the proof of Proposition 2.3, the estimates (2.25)-(2.26) 
remain valid here. Further, since $\rho^{2}/\delta $ is bounded and 
$\psi_{\varepsilon} \rightarrow $ 0 on $F$ by (2.7), it follows from 
(2.24) as before that
\begin{equation} \label{e2.36}
\Delta\nu_{\delta} = o(1). 
\end{equation}
Similarly, as before from the full equation in (2.24), one obtains
\begin{equation} \label{e2.37}
L^{*}(\frac{u}{\delta}) = \frac{\sigma T}{3}\frac{\rho^{2}}{\delta}\cdot  
g_{\varepsilon}'  + o(1),
\end{equation}
using the fact that $u \rightarrow 1$ and $-s \sim \rho^{2} \sigma T$ 
in this scale. When $\sigma(M) < 0$, one does not obtain here a solution 
to the linearized static vacuum equations, since
\begin{equation} \label{e2.38}
\frac{s_{g_{\varepsilon}'}}{\delta} = -\frac{\rho^{2}}{\delta}u\sigma T 
= \frac{\rho^{2}}{\delta}\sigma (M) + o(1). 
\end{equation}
The equations (2.36) and (2.37) imply, in the limit $\delta  
\rightarrow $ 0, the equations
\begin{equation} \label{e2.39}
z'  = D^{2}\nu , \ \  \Delta\nu  = 0. 
\end{equation}
These equations are the linearization of the equations for critical 
points of $\mathcal{S}_{-}^{2},$ i.e. the linearized 
$\mathcal{S}_{-}^{2}$ equations.

 Exactly the same arguments as in Proposition 2.3 shows that the limit 
$h$ exists, and there exist affine functions $a_{\delta}$ such that 
$\nu_{\delta} - a_{\delta}$ converges to a harmonic function $\nu $ on 
$(F, g_{o}').$ Similarly, for $\widetilde g_{\varepsilon}$ as in 
(2.29), the arguments in (2.30)-(2.31) give here that
\begin{equation} \label{e2.40}
\widetilde r_{\varepsilon}/\delta  = \frac{\sigma (M)}{3}\cdot  
g_{\varepsilon}'  + o(1), 
\end{equation}
and so
\begin{equation} \label{e2.41}
\widetilde z_{\varepsilon}/\delta  = 0 + o(1). 
\end{equation}
Thus, to first order in $\delta$, $\widetilde g_{\varepsilon}$ is of 
constant curvature, and hence $g_{\varepsilon}$ to first order is again 
conformally flat. More precisely, the formula (1.39) for conformal 
deformation and the estimate (2.36) give
\begin{equation} \label{e2.42}
\frac{\widetilde s}{\delta} = 
\frac{s_{g_{\varepsilon}'}}{\delta}\upsilon^{-2} + o(1) = 
-\frac{\rho^{2}}{\delta}\sigma T(1 + o(1)). 
\end{equation}
Hence, to first order in $\delta$, $\widetilde g_{\varepsilon}$ is 
locally isometric to the space form of constant curvature $-\delta 
(\frac{\rho^{2}}{\delta}\sigma T)/6.$ As always, by passing to 
subsequences, one may assume that $\rho^{2}/\delta  \rightarrow  r_{o} 
>  0$, while by (1.8) $-\sigma T \rightarrow  \sigma (M).$ Thus, to 
first order in $\delta$, $\widetilde g_{\varepsilon}$ is locally 
isometric to the space form $g_{-\delta}' $ of constant curvature 
$\delta r_{o}\sigma (M)/6.$ The remainder of the proof, in particular 
(2.35), then follows as in Proposition 2.3.
{\endproof}

 Next, we turn to the situation where $\delta  <<  \rho^{2},$ which is 
somewhat more complicated. 
\begin{remark} \label{r 2.5.(i).}
  Actually, the situation below where $\delta  <<  \rho^{2}$ or the 
situation (2.33) where $\delta  \sim  \rho^{2},$ although logically 
needed for the work through \S 5, are needed only in a comparatively 
minor way, (in terms of the overall picture). Hence on a first reading, 
it is advisable to skip at this point to the summary Theorem 2.11, and 
then onto \S 3, returning to this case as needed.

{\bf (ii).}
 Similarly, given Proposition A, the work to follow in \S 2.2 - \S 2.3 
also holds when $\rho (x_{\varepsilon}) \rightarrow  \infty $ as 
$\varepsilon  \rightarrow $ 0; this will be used however only in \S 7.3.
\end{remark}

\medskip

 Some preliminary estimates are needed in this case before deriving the 
analogue of Propositions 2.3 and 2.4. First, return to the equations 
(2.24), at the blow-up scale $g_{\varepsilon}' .$ The estimate (2.26) 
also holds here, i.e. on compact subsets of $F$,
\begin{equation} \label{e2.43}
\frac{\bar \alpha}{\delta}\nabla\mathcal{Z}^{2} \rightarrow  0, \ \ 
{\rm as} \ \ \varepsilon \rightarrow 0,  
\end{equation}
c.f. again Appendix B for the details of the proof. Taking the 
trace-free part of (2.24) then gives
\begin{equation} \label{e2.44}
\frac{z}{\delta} = D^{2}_{o}\nu_{\delta} + o(1), 
\end{equation}
where $\nu_{\delta} = (u- 1)/\delta $ as before and $D^{2}_{o}$ is the 
trace-free Hessian. 

 Observe that the quantity $z/\delta$ in (2.44) is both bounded away 
from 0 and bounded away from $\infty ,$ in the blow-up scale. Thus it 
has a limit, in a subsequence, and the limit is not identically 0. 

\medskip

 We return to the estimate (2.44), the trace-free analogue of (2.27) or 
(2.37) in Corollary 2.8 below, but first need to understand in detail 
the behavior of the trace equation in (2.24) in this setting. 

 Referring to $\psi_{\varepsilon}$ in (2.4), let $\psi_{\delta} = 
\psi_{\varepsilon}/\delta ,$ so that $\psi_{\delta}'  = \rho^{2}\cdot 
\psi_{\delta}$ in the $g_{\varepsilon}' $ scale (2.8). Hence from 
(2.43), or more directly from the definition of $\delta ,$ one has
\begin{equation} \label{e2.45}
2\Delta\nu_{\delta} = \rho^{2}(x_{\varepsilon})(\psi_{\delta} + o(1)) 
\leq  o(1). 
\end{equation}
Thus, in any limit, the potential is superharmonic. Of course, since 
the behavior of $u$ is being examined near its maximal value, this is 
not unexpected.

 The following result gives control on the size of the right side of 
(2.45) from below.

\begin{lemma} \label{l 2.6.}
  Let $x_{\varepsilon}\in (\Omega_{\varepsilon}, g_{\varepsilon})$ be 
any base points satisfying (2.1) and suppose one of the following 
conditions holds: there are positive constants $\chi_{o}$ or $\rho_{o}$ 
such that
\begin{equation} \label{e2.46}
|u- 1|(x_{\varepsilon}) \geq  \chi_{o}\cdot \bar 
\varepsilon\mathcal{Z}^{2}, \ \ {\rm or} 
\end{equation}
\begin{equation} \label{e2.47}
\rho (x_{\varepsilon}) \geq  \rho_{o}. 
\end{equation}
Then there exists a constant $K <  \infty ,$ depending only on 
$\chi_{o}$ or $\rho_{o},$ such that
\begin{equation} \label{e2.48}
\xi (x_{\varepsilon}) \equiv  
\rho^{2}(x_{\varepsilon})\psi_{\delta}(x_{\varepsilon}) \geq  - K. 
\end{equation}
\end{lemma}

\noindent
{\bf Proof:}
 The proof is by contradiction, and so suppose
$$\xi (x_{\varepsilon}) = 
\rho^{2}(x_{\varepsilon})\psi_{\delta}(x_{\varepsilon}) \rightarrow  
-\infty , \ \ {\rm as} \ \ \varepsilon  \rightarrow  0. $$
Renormalize the equations (2.45) and (2.48) by dividing them by $|\xi 
(x_{\varepsilon})|,$ so that
\begin{equation} \label{e2.49}
2\Delta\bar \nu_{\delta} = \bar \xi + o(1), 
\end{equation}
where $\bar \nu_{\delta} = \nu_{\delta}/|\xi (x_{\varepsilon})|, \bar 
\xi = \xi /|\xi (x_{\varepsilon})|,$ so that 
\begin{equation} \label{e2.50}
\bar \xi(x_{\varepsilon}) = - 1. 
\end{equation}
It then follows from (2.44) and the uniform bound on 
$(z_{g_{\varepsilon}'}/\delta )$ that
\begin{equation} \label{e2.51}
D^{2}_{o}\bar \nu_{\delta} = o(1), 
\end{equation}
on compact subsets of $F$. Consider first functions $\phi $ such that
\begin{equation} \label{e2.52}
D^{2}_{o}\phi  = 0, 
\end{equation}
on the limit $(F, g_{o}', x)$, $g_{o}'  = \lim g_{\varepsilon}'$. If 
$\rho (x_{\varepsilon}) \rightarrow $ 0, the limit $F$ is flat, while 
if $\rho (x_{\varepsilon}) \geq  \rho_{o} > $ 0, then the limit is 
either flat or hyperbolic. We treat first the case where $(F, g_{o}')$ 
is flat. 

 On a flat manifold, the only functions $\phi $ satisfying (2.52) 
locally are, modulo affine functions,
\begin{equation} \label{e2.53}
\phi  = b_{o}\cdot  t_{o}^{2}, 
\end{equation}
where $b_{o}$ is a constant and $t_{o}(x) = dist(x, p_{o}),$ for some 
base point $p_{o}\in F$. It follows from (2.51)-(2.53) that $\bar 
\nu_{\delta}$ converges to a limit function $\bar \nu_{\infty}$ modulo 
addition of affine functions and functions $\phi $ of the form (2.53). 
By (2.51), the limit function $\bar \nu_{\infty}$ necessarily satisfies 
(2.52) and hence 
$$\Delta\bar \nu_{\infty} = 6b_{o}. $$
The normalization (2.50) gives $b_{o} = - 1/12$ and further $\bar \xi = 
\bar \xi_{\varepsilon}$ converges to the constant function $- 1$ as 
$\varepsilon  \rightarrow $ 0. Hence, provided $y_{\varepsilon} 
\rightarrow  y\in F,$ from (2.4) we have
$$\bar \xi(y_{\varepsilon}) = \frac{1}{|\psi 
(x_{\varepsilon})|}[\tfrac{1}{4}\sigma T(u+1)\nu_{\delta} + \bar 
d_{\delta}](y_{\varepsilon}) \rightarrow  - 1, $$
where $\bar d_{\delta} = \bar d_{\varepsilon}/\delta$ and $\bar 
d_{\varepsilon}$ is as in (2.5). Via (2.49), this implies that on 
compact subsets of $F$,
\begin{equation} \label{e2.54}
2\Delta\bar \nu_{\delta} = 
2\Delta\frac{\nu_{\delta}}{\rho^{2}(x_{\varepsilon})|\psi 
(x_{\varepsilon})|} = \frac{1}{|\psi 
(x_{\varepsilon})|}[\tfrac{1}{4}\sigma T(u+1)\nu_{\delta} + \bar 
d_{\delta}] \rightarrow  - 1. 
\end{equation}

 Now suppose first (2.46) holds, and consider the right side of (2.54). 
The term $\hat d_{\delta} \equiv  \bar d_{\delta}/|\psi 
(x_{\varepsilon})|$ is constant, and $u \rightarrow $ 1 uniformly. 
Hence the oscillation of $\hat \nu_{\delta} \equiv  \nu_{\delta}/|\psi 
(x_{\varepsilon})|$ converges to 0 as $\delta  \rightarrow $ 0, i.e. 
\begin{equation} \label{e2.55}
\sup_{B}\hat \nu_{\delta} -  \inf_{B}\hat \nu_{\delta} \rightarrow  0, 
\end{equation}
where $B$ is any compact subset of $F$. The hypothesis (2.46) and the 
definition of $\bar d_{\delta}$ imply that 
$|\nu_{\delta}|(x_{\varepsilon}) = -\nu_{\delta}(x_{\varepsilon}) \geq  
-\chi_{o}\cdot \bar d_{\delta} = \chi_{o}\cdot |\bar d_{\delta}|$ and 
so $|\hat \nu_{\delta}|(x_{\varepsilon}) \geq  \chi_{o}\cdot |\hat 
d_{\delta}|.$ Via (2.55), this now implies that 
\begin{equation} \label{e2.56}
\sup_{B}\hat \nu_{\delta}/\inf_{B}\hat \nu_{\delta} \rightarrow  1, 
\end{equation}
and so of course $\sup_{B}\bar \nu_{\delta}/\inf_{B}\bar \nu_{\delta} 
\rightarrow $ 1. In turn, this means that $\Delta\bar \nu_{\delta} 
\rightarrow $ 0 weakly, which contradicts (2.54). An essentially 
identical argument holds if (2.47) holds, since in this case, $|\bar 
\nu| \leq  \rho_{o}^{-1}|\hat \nu|.$ 

\smallskip

 Next suppose the limit $(F, g_{o}')$ is hyperbolic, i.e. of constant 
negative curvature, so that in particular (2.47) automatically holds. 
The argument in this case is essentially the same as the flat case, 
with some minor modifications. Thus, since $u(x_{\varepsilon}) 
\rightarrow $ 1, the scalar curvature of $g_{\varepsilon}' $ is 
approximately $-\sigma T\rho^{2}, \rho  = \rho (x_{\varepsilon}),$ so 
that $g_{\varepsilon}' $ approximates a metric of constant sectional 
curvature $-\kappa^{2}, \kappa  = \rho (\sigma T/6)^{1/2}.$

 The arguments (2.49)-(2.52) hold as before. An elementary computation 
shows that the only solutions of (2.52) on the limit, (mod constants), 
are
\begin{equation} \label{e2.57}
\phi  = b_{o}\cdot  \cosh \kappa t_{o}, 
\end{equation}
with $t_{o}$ again the distance from some point on $F$. A simple 
computation then gives
\begin{equation} \label{e2.58}
\Delta\phi  = 3\kappa^{2}\phi . 
\end{equation}
On the other hand, arguing exactly as in (2.54), one has
$$2\Delta\frac{\hat \nu_{\delta}}{\rho^{2}} \sim  \frac{\sigma 
T}{2}\hat \nu_{\delta} + \hat d_{\delta}, $$
or equivalently
\begin{equation} \label{e2.59}
\Delta\bar \nu_{\delta} \sim  \frac{\sigma T}{4}\rho^{2}\bar 
\nu_{\delta} + \frac{1}{2}\hat d_{\delta}. 
\end{equation}
Since $\bar \nu_{\delta}$ converges to a solution of (2.52), (mod 
constants), $\bar \nu_{\delta}$ approaches a solution of (2.58). 
However, since $\kappa^{2} = \rho^{2}(\sigma T / 6)$, this contradicts 
(2.59), (since $\kappa > 0$ in the limit).
{\endproof}

\smallskip

 Lemma 2.6 leads to the following definition:

\begin{definition} \label {d 2.7.}
 Base points $x_{\varepsilon}\in (\Omega_{\varepsilon}, 
g_{\varepsilon})$ are called {\sf allowable}  if $u(x_{\varepsilon}) 
\rightarrow $ 1 and $|\Delta\nu_{\delta}|(x_{\varepsilon})$ remains 
uniformly bounded in the scale $g_{\varepsilon}' = 
\rho(x_{\varepsilon})^{-2} \cdot g_{\varepsilon}$, as $\varepsilon  = 
\varepsilon_{i} \rightarrow  0$.
\end{definition}

 Thus, the results above imply that $x_{\varepsilon}$ is allowable if 
either $\delta  \geq  c\cdot \rho^{2}$ at $x_{\varepsilon},$ as 
$\varepsilon  \rightarrow $ 0, for some $c > $ 0, or if $\delta  <<  
\rho^{2},$ then either (2.46) or (2.47) holds. Observe that if 
$x_{\varepsilon}$ is allowable, then both terms in $\xi 
(x_{\varepsilon}),$ i.e. $\rho^{2}\bar d_{\delta}$ and 
$\rho^{2}\nu_{\delta}(x_{\varepsilon})$, are bounded as $\delta  
\rightarrow $ 0. 

 Note that $x_{\varepsilon}$ allowable does {\it  not}  imply that the 
linearized scalar curvature, i.e. $s_{g_{\varepsilon}'}/\delta $ 
remains uniformly bounded at $x_{\varepsilon}$ as $\varepsilon  
\rightarrow $ 0. On the other hand, by the statement following (2.44) 
and Lemma 2.6, the full Hessian $|D^{2}\nu_{\delta}|(x_{\varepsilon})$ 
does remain bounded on allowable sequences.

\smallskip

 The discussion above proves the following partial analogue of 
Propositions 2.3 and 2.4.
\begin{corollary} \label{c 2.8.}
  Suppose the base points $x_{\varepsilon}$ are allowable and suppose 
that
\begin{equation} \label{e2.60}
\delta  <<  \rho^{2}(x_{\varepsilon}), \ \ {\rm as} \ \  \varepsilon  
\rightarrow  0. 
\end{equation}
Then there exists a sequence of affine functions $a_{\delta}$ such that 
the functions $\nu_{\delta}- a_{\delta}$ converge to a non-affine limit 
function $\nu $ on the maximal limit $(F, g_{o}', x)$. The 
linearization of the metrics $g_{\varepsilon}' $ at $\varepsilon  =$ 0 
gives rise to a solution of the equations,
\begin{equation} \label{e2.61}
z'  = D^{2}_{o}\nu , \ \ 2\Delta\nu  = \xi  \leq  0, 
\end{equation}
where $\xi  = \lim \rho^{2}(x_{\varepsilon})(\frac{1}{2}\sigma 
T\nu_{\delta} + \bar d_{\delta})$ on $(F, g_{o}')$.  
\end{corollary}

\noindent
{\bf Proof:}
 The first equation in (2.61) follows from the limit of (2.44). For the 
second equation in (2.61), as noted above the constant term 
$\rho^{2}\bar d_{\delta}$ in $\xi $ is bounded as $\varepsilon  
\rightarrow $ 0. Further, the coefficient $\rho^{2}\sigma T$ of 
$\nu_{\delta}$ in $\xi $ is also bounded, since it approximates the 
scalar curvature of $g_{\varepsilon}' .$ Hence, elliptic regularity for 
the trace equation in (2.45) implies that $\xi  = \xi_{\varepsilon}$ is 
locally bounded on $F$ and so (2.61) follows as the limit of (2.45). 
The proof that the functions $\nu_{\delta}$ converge to $\nu $ modulo 
affine functions follows for the same reasons as before in Proposition 
2.3.
{\endproof}

\smallskip

 The situation in Corollary 2.8 of course differs somewhat from that in 
Propositions 2.3 and 2.4, where the limit potential $\nu $ is harmonic. 
It is not to be expected that $\xi  =$ 0 in (2.61) in general in 
situations where $\delta  <<  \rho^{2}.$ The form of the metric 
$g_{\varepsilon}' ,$ i.e. the analogue of (2.23) or (2.35) in this 
situation, is discussed below in Theorem 2.11.

\begin{remark} \label{r 2.9.}
  If the limit $(F, g_{o})$ in Corollary 2.8 is hyperbolic, i.e. of 
constant negative curvature, then the functions $\nu_{\delta}$ converge 
modulo constants to the limit $\nu ;$ this follows since hyperbolic 
manifolds do not carry, even locally, any non-constant affine functions.
\end{remark}

\begin{remark} \label{r 2.10.}
  Suppose the base points $x_{\varepsilon}$ are non-allowable, in that 
both (2.46) and (2.47) fail. The implication that (2.55) implies (2.56) 
no longer holds, since one may have $|\hat \nu_{\delta}| <<  |\hat 
d_{\delta}|.$ Hence $\xi  = \psi_{\delta}' $ may become unbounded as 
$\delta  \rightarrow $ 0. Observe that in this case the dominant term 
in $\xi $ is $\rho^{2}\bar d_{\delta}.$ Whether this dominant term 
$\rho^{2}\bar d_{\delta}$ stays bounded or not depends on the relation 
of $\rho^{2}/\delta ,$ which depends only on the local geometry at 
$x_{\varepsilon},$ with the global term $\bar 
\varepsilon\mathcal{Z}^{2}.$ There are other, much more fundamental 
reasons to work only with allowable base points which arise in \S 4 and 
\S 5 below.
\end{remark}

\smallskip

{\bf 2.3.}
 In this section, we complete and unify the three forms of the 
linearizations, according to the three cases of (2.12). To do this, 
given allowable base points $x_{\varepsilon}\in\Omega_{\varepsilon},$ 
let $g_{-\kappa}$ be the space form of constant curvature 
\begin{equation} \label{e2.62}
-\kappa^{2}  = 
-\frac{T\sigma_{g_{\varepsilon}}}{6}\rho^{2}(x_{\varepsilon}). 
\end{equation}
If $\delta (x_{\varepsilon}) >>  \rho^{2}(x_{\varepsilon})$ as 
$\varepsilon  \rightarrow $ 0, so $\kappa^{2} << \delta$, set
\begin{equation} \label{e2.63}
g_{-\kappa} '  \equiv  g_{o}' , 
\end{equation}
where $g_{o}' $ is the flat metric on $F$. If $\delta (x_{\varepsilon}) 
\sim  \rho^{2}(x_{\varepsilon})$ or $\delta (x_{\varepsilon}) <<  
\rho^{2}(x_{\varepsilon}),$ set
\begin{equation} \label{e2.64}
g_{-\kappa} '  \equiv  g_{-\kappa}. 
\end{equation}

\begin{theorem} \label{t 2.11.}
  Let $x_{\varepsilon}$ be any sequence of allowable base points. Then 
there is a sequence of affine functions $a_{\delta},$ either $0$ or 
diverging to $-\infty$, such that $\nu_{\delta} -  a_{\delta}$ 
sub-converges to a non-affine function $\nu ,$ defined on the maximal 
constant curvature limit $(F, g_{o}', x)$, $x = \lim x_{\varepsilon}$ 
of $(\Omega_{\varepsilon}, g_{\varepsilon}' , x_{\varepsilon}).$ The 
linearization $h = \lim_{\delta\rightarrow 0} h_{\varepsilon}$ from 
(2.14) exists and gives rise to a non-trivial solution of the equations
\begin{equation} \label{e2.65}
z'  = D^{2}_{o}\nu , \ \  \Delta\nu  = \lambda\nu  + a_{\infty} \leq  
0, 
\end{equation}
where $\lambda $ is a constant, either 0 or positive, and $a_{\infty}$ 
is a non-positive affine function, possibly constant. For $\varepsilon 
$ sufficiently small, the metric $g_{\varepsilon}' $ has the form, 
modulo diffeomorphisms,
\begin{equation} \label{e2.66}
g_{\varepsilon}'  = (1 -  2\nu\delta )(g_{-\kappa}' +\delta\chi ) + 
o(\delta ), 
\end{equation}
where $\chi $ is an infinitesimal holonomy deformation of $(F, 
g_{-\kappa}')$.
\end{theorem}

\noindent
{\bf Proof:}
 The first equation in (2.65) has already been verified in (2.22), 
(2.34) and (2.61). As in (2.29), consider the conformally equivalent 
metric $\widetilde g_{\varepsilon} = \upsilon^{2}g_{\varepsilon}' ,$ 
for $\upsilon $ as following (2.28). The trace-free part of (1.38) 
implies
\begin{equation} \label{e2.67}
\frac{\widetilde z_{\varepsilon}}{\delta} = 
\frac{z_{\varepsilon}}{\delta} -  
\frac{D^{2}_{o}\upsilon_{\delta}}{\upsilon} + 
2\frac{(d\upsilon_{\delta}d\upsilon )_{o}}{\upsilon^{2}}, 
\end{equation}
where the subscript $o$ denotes trace-free part and $\upsilon_{\delta} 
= (\upsilon -1)/\delta .$ Since $\upsilon_{\delta}$ converges to the 
limit $\nu$, $d\upsilon_{\delta}$ is bounded, while $d\upsilon  
\rightarrow $ 0 uniformly on compact subsets of $F$. Hence, by the 
first equation in (2.65), one obtains
\begin{equation} \label{e2.68}
\frac{\widetilde z_{\varepsilon}}{\delta} \rightarrow  0, \ \ {\rm as} 
\ \  \varepsilon  \rightarrow  0, 
\end{equation}
on $F$, so that $\widetilde g_{\varepsilon}$ is of constant curvature, 
and hence $g_{\varepsilon}' $ is conformally flat, both to first order 
in $\delta .$

 For the trace equation in (2.65), by (2.45) one has
\begin{equation} \label{e2.69}
\Delta\nu_{\delta} = \frac{\rho^{2}}{8}\sigma T(u+1)\nu_{\delta} + 
\frac{\rho^{2}\bar d_{\delta}}{2} + o(1), 
\end{equation}
where $\bar d_{\delta} = \bar d_{\varepsilon}/\delta .$ Recall that $u 
\rightarrow $ 1 uniformly on compact subsets of the limit $F$. Lemma 
2.6 implies that both terms on the right side of (2.69) are uniformly 
bounded as $\varepsilon  \rightarrow $ 0, while (2.5ff) implies both 
are non-positive. As in (2.28), $\nu_{\delta} = \nu  + a_{\delta} + 
o(1)$, where $a_{\delta}$ is either 0 or a divergent sequence of affine 
functions, so that (2.69) has the form
\begin{equation} \label{e2.70}
\Delta\nu_{\delta} = \frac{\rho^{2}}{4}\sigma T(\nu +a_{\delta}) + 
\frac{\rho^{2}\bar d_{\delta}}{2} + o(1). 
\end{equation}
Suppose first that there exists $c_{o} > 0$ such that, as $\varepsilon 
\rightarrow 0$,
\begin{equation} \label{e2.71}
\rho^{2}\sigma T \geq  c_{o}, 
\end{equation}
Then the term $a_{\delta}$ must be uniformly bounded as $\varepsilon  
\rightarrow $ 0, and so it may be ignored, i.e. absorbed into $\nu$. In 
this case, (2.70) limits to
\begin{equation} \label{e2.72}
\Delta\nu  = \lim(\frac{\rho^{2}}{4}\sigma T)\nu + 
\lim(\frac{\rho^{2}\bar d_{\delta}}{2}), 
\end{equation}
giving (2.65). On the other hand, if $\rho^{2}\sigma T \rightarrow $ 0, 
then $\rho^{2}\sigma T\nu  \rightarrow $ 0 also, so the limit is
\begin{equation} \label{e2.73}
\Delta\nu  = \lim(\frac{\rho^{2}}{4}\sigma Ta_{\delta})+ 
\lim(\frac{\rho^{2}\bar d_{\delta}}{2}), 
\end{equation}
giving (2.65) with $\lambda  = 0$.

 Next we relate the scalar curvatures $\widetilde s_{\varepsilon}$ and 
$s_{\varepsilon}' .$ A short calculation using (1.39) shows that
\begin{equation} \label{e2.74}
\frac{\widetilde s_{\varepsilon}}{\delta} = - 
4\upsilon^{-3}\Delta\nu_{\delta} + 2u^{-4}|d\nu_{\delta}||d\upsilon| + 
\upsilon^{-2}\frac{s_{\varepsilon}'}{\delta}, 
\end{equation}
where again the right side of (2.74) is w.r.t the blow-up metric 
$g_{\varepsilon}' .$ As noted above $d\nu_{\delta} \rightarrow  d\nu ,$ 
and so is bounded, while $d\upsilon  \rightarrow $ 0 as $\varepsilon  
\rightarrow $ 0. Hence the second term in (2.74) tends to 0 as $\delta  
\rightarrow $ 0. Further, 
\begin{equation} \label{e2.75}
\upsilon^{-2}\frac{s_{\varepsilon}'}{\delta} = -  
\upsilon^{-2}\frac{\rho^{2}}{\delta}T\sigma u, 
\end{equation}
so that (2.74) becomes
\begin{equation} \label{e2.76}
\frac{\widetilde s_{\varepsilon}}{\delta} = - 4u^{-3}\Delta\nu  -  
\upsilon^{-2}\frac{\rho^{2}}{\delta}T\sigma u + o(1). 
\end{equation}
Referring to (2.72) and (2.73), $\nu\cdot \delta$, $a_{\delta}\cdot 
\delta$ and $d_{\delta}\cdot \delta $ all converge to 0 as $\delta  
\rightarrow $ 0. Hence the second term in (2.76) dominates the first, 
so that
\begin{equation} \label{e2.77}
\frac{\widetilde s_{\varepsilon}}{\delta} = -  
\frac{\rho^{2}}{\delta}T\sigma (1 + o(1)). 
\end{equation}
In particular, one sees that $\widetilde s_{\varepsilon}$ is constant, 
to $1^{\rm st}$ order in $\delta$, which of course it must be by (2.68) 
and the Bianchi identity. The estimates (2.68) and (2.77) thus imply 
that $\widetilde g_{\varepsilon}$ is isometric, to first order in 
$\delta ,$ to the space form $g_{-\kappa}' ,$ as in (2.67)-(2.68). The 
remainder of the proof follows as in Proposition 2.3.
{\endproof}

 We point out that if (2.71) holds, i.e. $\lambda > 0$ in (2.65), then 
the proof shows that $\nu_{\delta}$ converges to $\nu ,$ i.e. there is 
no constant or affine indeterminacy; in particular, $\nu  \leq $ 0 in 
this situation.

\medskip

 Observe that Theorem 2.11 implies that {\it  any} allowable sequence 
$\{x_{\varepsilon}\}$ gives rise to a {\it non-trivial} and {\it 
uniformly controlled} solution of the corresponding linearized 
equations, modulo addition of affine functions. This fact plays a 
crucial role throughout the rest of the paper. 

 In particular, at any allowable base point sequence 
$\{x_{\varepsilon}\},$ there is a constant $\kappa_{o} > $ 0, 
independent of $\varepsilon $ and $x_{\varepsilon},$ such that, for 
most $y_{\varepsilon}\in B_{x_{\varepsilon}}(\rho (x_{\varepsilon})),$
\begin{equation} \label{e2.78}
|d\nu_{\delta (x_{\varepsilon})}|(y_{\varepsilon}) \geq  
\kappa_{o}/\rho (x_{\varepsilon}). 
\end{equation}
Note that the estimate (2.78) is scale-invariant. There is a minor 
technical issue here in that (2.78) need not hold for {\it all} 
$y_{\varepsilon}\in B_{x_{\varepsilon}}(\rho (x_{\varepsilon}));$ for 
instance, if $y_{\varepsilon}$ is a critical point of $u$, then 
$|d\nu_{\delta (x_{\varepsilon})}|(y_{\varepsilon}) =$ 0. However, 
(2.78) does hold for generic $y_{\varepsilon} \in 
B_{x_{\varepsilon}}(\rho(x_{\varepsilon}))$, since if not, the limit 
potential $\nu$ on the maximal limit $(F, g_{o}', x)$ would be 
constant, (modulo affine functions), contradicting Theorem 2.11. Thus, 
(2.78) holds for a large percentage of $y_{\varepsilon}\in 
B_{x_{\varepsilon}}(\rho (x_{\varepsilon})),$ where $\kappa_{o}$ 
depends only on the volume percentage. In particular, one can always 
slightly adjust the choice of $x_{\varepsilon},$ within $\rho 
(x_{\varepsilon})$, so that (2.78) holds at $y_{\varepsilon} = 
x_{\varepsilon}.$ Such a choice of $x_{\varepsilon}$, (of course with 
$u(x_{\varepsilon}) \rightarrow 1$), will be called {\it generic}, and it 
will always be assumed, (usually implicitly), in the work to follow 
that base points are generic. 

 The uniformity in Theorem 2.11 also implies that for generic base 
points there is a constant $a_{o} > 0$, depending only on the collapse 
ratio $\mu_{o} = \omega(x_{\varepsilon}) / \rho(x_{\varepsilon})$ such 
that
\begin{equation} \label{e2.79}
a_{o} \cdot \rho^{2}(x_{\varepsilon}) \leq 
area(B_{x_{\varepsilon}}(\rho (x_{\varepsilon})\cap L)) \leq  
a_{o}^{-1}\cdot \rho^{2}(x_{\varepsilon}), 
\end{equation}
where $L$ is the level set of $\nu_{\delta}$ through $x_{\varepsilon}.$

 Finally, it is worth pointing out that Theorem 2.11 is the reason for 
choosing $\delta  = \delta_{z}$ and not $\delta  = \delta_{r};$ 
linearizations w.r.t. $\delta_{r}$ in the situation $\delta << 
\rho^{2}$ always lead to trivial solutions of the linearized equations, 
at least when $\sigma (M) <  0$.

\section{Level Set Measures and Masses.}
\setcounter{equation}{0}

 In this section, we describe natural measures and corresponding mass 
functions associated to the potential functions $\nu_{\delta}$ and 
their limits $\nu $ from \S 2. These measures are supported on the 
level sets of $u$, and it is the distribution of these measures in 
relation to the geometry of the metrics $g_{\varepsilon},$ or 
$g_{\varepsilon}' ,$ which plays the central role in locating the 
2-spheres in $M$.

\medskip

 First, recall that the choice of $\delta $ in (2.10) led to the 
constant and affine indeterminacy in the construction of the limit 
potential $\nu ,$ at any given base point sequence $x_{\varepsilon}$, 
(as always satisfying (2.1)). At a later point in the paper in \S 5, it 
will be important to remove the affine indeterminacy of the potential 
$\nu_{\delta};$ however, until then this indeterminacy does not play a 
major role. 

 In analogy to the definition (2.10) of $\delta ,$ define then
\begin{equation} \label{e3.1}
\delta_{a} = \delta_{a}(x_{\varepsilon}) = 
\oint_{B_{x_{\varepsilon}}(\frac{1}{2})\cap L_{\varepsilon}}|\nabla 
u|dA_{g_{\varepsilon}'} = \rho (x_{\varepsilon})\cdot 
\oint_{B_{x_{\varepsilon}}(\frac{1}{2}\rho (x_{\varepsilon}))\cap 
L_{\varepsilon}}|\nabla u|dA_{g_{\varepsilon}}, 
\end{equation}
where the subscript stands for affine, $L_{\varepsilon}$ is the level 
set of $u$ through the base point $x_{\varepsilon}$ and $\oint$ is the 
average value, as before. The first integral is w.r.t $g_{\varepsilon}' 
$ while the second is w.r.t. $g_{\varepsilon},$ so that $\delta_{a}$ is 
scale-invariant, as $\delta $ is. Theorem 2.11 implies that the 
potentials $\nu_{\delta (x_{\varepsilon})},$ at any allowable base 
point sequence $x_{\varepsilon},$ (sub)-converge, modulo affine 
functions, to non-trivial limit functions. As mentioned in connection 
with (2.78), it is also assumed henceforth that $x_{\varepsilon}$ is 
generic. This implies immediately that
\begin{equation} \label{e3.2}
\delta_{a} \geq   c_{1}\cdot \delta , 
\end{equation}
for a fixed numerical constant $c_{1} > $ 0, independent of 
$\varepsilon $ and the generic base point $x_{\varepsilon}.$ The affine 
indeterminacy of $\nu_{\delta}$ corresponds to the possibility that
\begin{equation} \label{e3.3}
\delta_{a} >> \delta . 
\end{equation}
The reason for preferring $\delta_{a}$ to $\delta ,$ (at this stage), 
is that the potential function
\begin{equation} \label{e3.4}
\nu_{\delta_{a}(x_{\varepsilon})} = (u- 1)/\delta_{a}(x_{\varepsilon}) 
= \nu_{\delta (x_{\varepsilon})}\cdot  (\frac{\delta 
(x_{\varepsilon})}{\delta_{a}(x_{\varepsilon})}) 
\end{equation}
converges, modulo addition of constants, to a limit function $\nu  = 
\nu_{a},$ so that the affine indeterminacy is reduced to a constant 
indeterminacy. When (3.3) occurs at $x_{\varepsilon},$ the potential 
$\nu_{\delta_{a}}$ converges to a non-constant affine function 
$\nu_{a}$, with $|\nabla \nu_{a}|(x) = 1$, on the limit $(F, g_{o}', 
x_{\varepsilon}),$ while if 
\begin{equation} \label{e3.5}
\delta_{a} \leq  k_{o}\cdot \delta , 
\end{equation}
for some fixed constant $k_{o} <  \infty ,$ (on base points 
$x_{\varepsilon}),$ then the functions $\nu_{\delta_{a}}$ and 
$\nu_{\delta}$ have uniformly bounded ratios, and thus the limit 
functions $\nu $ and $\nu_{a}$ differ merely by a multiplicative 
factor. In particular, Theorem 2.11 holds also w.r.t. $\delta_{a}$ in 
place of $\delta $, with the understanding that the limit function 
$\nu_{a}$ is affine in case (3.3) holds.

\smallskip

 Note the following elementary result.
\begin{lemma} \label{l 3.1.}
  Suppose the base points $x_{\varepsilon}\in (\Omega_{\varepsilon}, 
g_{\varepsilon})$ are allowable. Then there is a constant $C_{o} <  
\infty$, (independent of $\varepsilon$), such that
\begin{equation} \label{e3.6}
|du|(x_{\varepsilon}) \leq  C_{o}\cdot 
\frac{\delta_{a}(x_{\varepsilon})}{\rho (x_{\varepsilon})}. 
\end{equation}
\end{lemma}

\noindent
{\bf Proof:}
 To see this, observe that (3.6) is equivalent to the scale invariant 
estimate 
$$\rho (x_{\varepsilon})|d\nu_{\delta_{a}}|(x_{\varepsilon}) \leq  
C_{o}. $$
Theorem 2.11 implies that the functions $\nu_{\delta_{a}}$ 
(sub)-converge smoothly, modulo constants, to a limit function 
$\nu_{a}$ on the maximal limit $(F, g_{o}', x)$, $x = \lim 
x_{\varepsilon}$, (unwrapping in the case of collapse as described in 
\S 1). By construction and elliptic regularity, the limit satisfies 
$|d\nu_{a}|(x) \leq  C_{o},$ with $C_{o}$ independent of the base 
points $x_{\varepsilon}.$ This gives (3.6), since the convergence to 
the limit is smooth.
{\endproof}

 For the same reasons, if $y_{\varepsilon}$ are points with 
$y_{\varepsilon} \rightarrow  y \in  (F, g_{o}')$, then there is a 
constant $C_{1} <  \infty ,$ depending only on $dist_{g_{o}'}(y, 
\partial F)$, and $dist_{g_{o}'}(x, y)$ such that
\begin{equation} \label{e3.7}
|du|(y_{\varepsilon}) \leq  C_{1}\cdot 
\frac{\delta_{a}(x_{\varepsilon})}{\rho (x_{\varepsilon})}. 
\end{equation}

 We now introduce Riesz-type measures and mass functions for the 
potentials $\nu_{\delta_{a}}$ and limit potentials $\nu_{a}$ 
constructed above. Let $L = L_{\varepsilon}$ be the level set of $u$ 
through the base point $x_{\varepsilon}$ and define
\begin{equation} \label{e3.8}
d\mu  = d\mu_{\delta_{a}}(L) = 
|\nabla\nu_{\delta_{a}}|dA_{g_{\varepsilon}}, 
\end{equation}
where $\delta_{a} = \delta_{a}(x_{\varepsilon})$ and 
$dA_{g_{\varepsilon}}$ is Lebesgue measure w.r.t. $g_{\varepsilon}$ on 
the level $L$. This measure depends on the base point, both with regard 
to $\delta_{a}$ and the level $L$. Observe that $d\mu$ is unchanged 
when adding constants to $\nu_{\delta_{a}}.$ For any compact set 
$E\subset F,$ one has the associated mass
\begin{equation} \label{e3.9}
m(E) = \int_{E \cap L}|\nabla\nu_{\delta_{a}}|dA_{g_{\varepsilon}}. 
\end{equation}
In particular, for any $q\in L,$ define the mass function $m_{q}$ of 
$d\mu $ by
\begin{equation} \label{e3.10}
m_{q}(s) = \int_{L\cap 
B_{q}(s)}|\nabla\nu_{\delta_{a}}|dA_{g_{\varepsilon}}. 
\end{equation}

 Similar measures and mass functions are defined on the other level 
sets $\bar L$ of $u$ as in (3.8)-(3.10), with $d\bar A$ the Lebesgue 
measure on $\bar L,$ and with the same potential $\nu_{\delta_{a}}.$ 
These measures on distinct level sets can be related to each other by 
means of the divergence theorem applied to the trace equation (2.3).

 The mass functions in (3.9)-(3.10) scale as distance. Thus, if 
$g_{\varepsilon}'  = \rho (x_{\varepsilon})^{-2}\cdot  
g_{\varepsilon},$ then $m' (E) = \rho (x_{\varepsilon})^{-1}m(E)$ and 
$m_{q}' (s) = \rho (x_{\varepsilon})^{-1}m_{q}(\rho 
(x_{\varepsilon})\cdot  s).$

\medskip

 When the blow-ups $g_{\varepsilon}' $ (sub)-converge, to the limit 
$(F, g_{o}', x)$, the rescaled measure $d\mu'$ and mass $m'$ converge 
to that of the limit potential $\nu_{a}$ satisfying (2.65). The same 
holds on each level set $\bar L_{\varepsilon} \rightarrow  \bar L 
\subset F$. Note that in case $\Delta\nu_{a}  =$ 0, the limit measure 
$d\mu'$ is the Riesz measure of the subharmonic function $\nu_{\bar L} 
= \sup(\nu_{a}, \nu_{a}(\bar L))$ on $F$. 

 However, when the blow-ups $g_{\varepsilon}'$ collapse at 
$x_{\varepsilon},$ the mass $m'$ as in (3.10) converges to 0. Thus, if 
the volume radius of $x_{\varepsilon}$ is small compared with the 
curvature radius, i.e. $\omega (x_{\varepsilon}) \leq  \mu_{o}\cdot 
\rho (x_{\varepsilon}),$ for some $\mu_{o} > $ 0 small but fixed, 
unwrap the collapse by passing to large finite covers as described in 
\S 1, so that $\omega (x_{\varepsilon}) \sim \rho (x_{\varepsilon})$, 
and define the measure $d\widetilde \mu$ as in (3.8) in the cover. This 
is essentially equivalent to renormalizing the mass (3.9) as follows:
\begin{equation} \label{e3.11}
\widetilde m_{\bar L}(E) = 
\frac{\rho^{2}(x_{\varepsilon})}{area_{g_{\varepsilon}}(L \cap 
B_{x_{\varepsilon}}(\tfrac{1}{2}\rho(x_{\varepsilon})))}\int_{E \cap 
\bar L}|\nabla\nu_{\delta_{a}}|dA_{g_{\varepsilon}}, 
\end{equation}
where $L$ is the level through $x_{\varepsilon}$. When the sequence 
based at $x_{\varepsilon}$ is not $\mu_{o}$-collapsed, then $\widetilde 
m(E) \sim  m(E)$, for sets $E$ whose diameter is on the order of $\rho 
(x_{\varepsilon})$, c.f. (2.79). Observe again that the mass (3.11) 
scales as a distance. In particular, in the blow-up scale 
$g_{\varepsilon}'  = \rho (x_{\varepsilon})^{-2}g_{\varepsilon}$ 
attached to $x_{\varepsilon},$ one has, on $\bar L$,
\begin{equation} \label{e3.12}
\widetilde m' (E) = \frac{1}{area_{g_{\varepsilon}'}(L \cap 
B_{x_{\varepsilon}}'(\tfrac{1}{2}))}\int_{E \cap \bar 
L}|\nabla\nu_{\delta_{a}}|dA_{g_{\varepsilon}'}. 
\end{equation}
Note from the definitions (3.1), (3.10) and (3.12) that $\widetilde 
m_{x_{\varepsilon}}' (1/2) = 1$, and correspondingly,
\begin{equation} \label{e3.13}
\widetilde m_{x_{\varepsilon}}(\tfrac{1}{2}\rho (x_{\varepsilon})) = 
\rho (x_{\varepsilon}). 
\end{equation}

 Each of these notions may be defined w.r.t. the 'curvature' $\delta $ 
of (2.10) in place of the 'affine' $\delta $ of (3.1). In situations 
where (3.5) holds, the corresponding measures and mass functions are 
uniformly bounded, in terms of $k_{o},$ with respect to each other. 
However, when (3.3) holds, the measures and mass functions w.r.t. 
$\delta $ become unbounded; this is why $\delta_{a}$ is used in place 
of $\delta ,$ for the time being.

\smallskip

 Next, we derive some estimates for the behavior of the limit potential 
$\nu_{a}$ on $F$. We begin with the following upper bound on $\nu_{a}$.

\begin{lemma} \label{l 3.2.}
  Let $x_{\varepsilon}$ be an allowable base point sequence satisfying 
(2.46). Then the potential function $\nu_{a}$ on the maximal limit $(F, 
g_{o}', x)$, normalized so that $\nu_{a}(x) = 0$, is bounded above on 
compact subsets of $\bar F,$ i.e. for any compact set $K \subset  \bar 
F,$
\begin{equation} \label{e3.14}
\nu_{a}(y) \leq  C, 
\end{equation}
for all $y\in K,$ where $C = C(K)$. The same estimate holds for the 
potentials $\nu_{\delta_{a}(x_{\varepsilon})}$, renormalized by 
additive constants, converging to $\nu_{a}$.
\end{lemma}

\noindent
{\bf Proof:}
 If $\delta_{a}(x_{\varepsilon}) >> \delta(x_{\varepsilon})$, i.e. 
(3.3) holds, then the limit function $\nu_{a}$ is affine and so (3.14) 
is immediate. Thus, one may assume that (3.5) holds and so replace 
$\delta_{a}$ by $\delta$. The functions $\nu_{\delta}$ thus 
sub-converge, modulo constants, to the limit function $\nu$.

 Recall from Theorem 2.11 that the potential $\nu$ satisfies the trace 
equation (2.65). If $\lambda > 0$, then as noted following Theorem 
2.11, $\nu \leq 0$ and so again (3.14) follows. Hence, one may assume 
that (2.65) is given by
\begin{equation} \label{e3.15}
\Delta\nu = \bar d_{o} \leq  0, 
\end{equation}
where $\bar d_{o} = \lim (\rho^{2}(x_{\varepsilon}) / 
\delta(x_{\varepsilon}))\bar d_{\varepsilon} >  -\infty$. Note that 
$\bar d_{\varepsilon} \rightarrow $ 0, while the ratio 
$\rho^{2}(x_{\varepsilon})/\delta(x_{\varepsilon})$ may remain bounded, 
or go to $\infty$; the product remains bounded as $\varepsilon  
\rightarrow 0$, since $x_{\varepsilon}$ is allowable.

 Since $\nu$ is smooth on $F$, $\nu$ can diverge to $+\infty$ on a 
compact set $K$ in $\bar F$ only on approach to $\partial F$. Let $p$ 
be any fixed point in $F \cap K$ and let $p_{\varepsilon} \in 
U^{\upsilon_{o}}$ be a sequence converging to $p$. Let $\hat 
\nu_{\delta(p_{\varepsilon})} = \nu_{\delta(p_{\varepsilon})} - 
\nu_{\delta(p_{\varepsilon})}(p_{\varepsilon})$ be the normalized 
potentials converging to the limit (normalized) potential $\nu_{p}$. 
Then it suffices to prove that $\hat \nu_{\delta(p_{\varepsilon})}$ is 
uniformly bounded above on the component $D' = D_{p_{\varepsilon}}'$ of 
$B_{p_{\varepsilon}}'(\frac{10}{9}) \cap U^{\upsilon_{o}}$ containing 
$p_{\varepsilon}$; here $\frac{10}{9}$ may be replaced by $1+\mu$, for 
any fixed $\mu > 0$, and the ball $B'$ is the ball in the metric 
$g_{\varepsilon}' = \rho(p_{\varepsilon})^{-2} \cdot g_{\varepsilon}$ 
based at $p_{\varepsilon}$.

  Suppose this is not the case, so that there exist $q_{\varepsilon} 
\in D_{p_{\varepsilon}}'$, for some sequence $p_{\varepsilon} 
\rightarrow p$, such that
\begin{equation} \label{e3.16}
\hat \nu_{\delta(p_{\varepsilon})}(q_{\varepsilon}) \rightarrow  \infty 
. 
\end{equation}
Then we claim that
\begin{equation} \label{e3.17}
\frac{\rho^{2}(q_{\varepsilon})}{\delta(q_{\varepsilon})} << 
\frac{\rho^{2}(p_{\varepsilon})}{\delta(p_{\varepsilon})}.
\end{equation}
To see this, recall that $\delta(p_{\varepsilon})$ is the $L^{2}$ 
average of $z$ on the rescaled ball $B_{p_{\varepsilon}}'(\frac{1}{2})$ 
based at $p_{\varepsilon}$. The relation (3.16) implies that 
$|D_{o}^{2}\nu_{\delta(p_{\varepsilon})}| >> 1$, near $q_{\varepsilon}$ 
and hence by the relation (2.44), $z/\delta(p_{\varepsilon}) >> 1$ near 
$q_{\varepsilon}$, in the scale based at $p_{\varepsilon}$. 
Transferring this estimate to the scale based at $q_{\varepsilon}$ 
gives (3.17).

  It follows from (3.17) and the definition of $\bar d_{o}$ in (3.15) 
that on the limit $(F, g_{o}', q)$, based at $q = \lim 
q_{\varepsilon}$, $g_{\varepsilon}' = \rho(q_{\varepsilon})^{-2} \cdot 
g_{\varepsilon}$, that the limit potential $\nu_{q} = \lim \hat 
\nu_{\delta(q_{\varepsilon})}$ is harmonic, (since $\bar 
d_{\varepsilon}$ is independent of base point). In particular, 
$\nu_{q}$ has no local maxima. Hence the approximating potentials $\hat 
\nu_{\delta(q_{\varepsilon})}$, restricted to the rescaled balls 
$B_{q_{\varepsilon}}'(1)$, have a  maximal value a definite amount 
larger than the value at $q_{\varepsilon}$. Note that this gives a 
contradiction if for instance the component of $\partial F$ containing 
$q$ is compact.

\smallskip

 Now observe this argument holds with respect to {\it any} allowable 
base point sequence $y_{\varepsilon}.$ Thus, given $x_{\varepsilon},$ 
let $L_{\varepsilon}$ be the $u$-level through $x_{\varepsilon}$ and 
let $R_{\varepsilon}$ be the component of the region bounded by 
$L^{\upsilon_{o}}$ and $L_{\varepsilon}$ containing $x_{\varepsilon}$. 
For each $y_{\varepsilon} \in L_{\varepsilon} \cap \partial 
R_{\varepsilon}$, consider the function $\hat 
\nu_{\delta(y_{\varepsilon})}$ on $D_{y_{\varepsilon}}'$ as above. The 
potential $\hat \nu_{\delta(y_{\varepsilon})}$ has a maximal value on 
the closure $\bar D_{y_{\varepsilon}}'$ of $D_{y_{\varepsilon}}'$, and 
so one thus obtains a function $\phi_{\varepsilon}$ on 
$L_{\varepsilon}$. The function $\phi_{\varepsilon}$ itself has a 
maximal value, achieved, (or arbitrarily close to being achieved, since 
$\Omega_{\varepsilon}$ is complete), at say $y_{\varepsilon}$. Let 
$q_{\varepsilon} \in \bar D_{y_{\varepsilon}}'$ be a point realizing 
the maximum of $\hat \nu_{\delta(y_{\varepsilon})}$. It follows that 
(3.16) and (3.17) hold, with $y_{\varepsilon}$ in place of 
$p_{\varepsilon}$. As before, this gives a contradiction if 
$q_{\varepsilon}$ is in the interior $D_{y_{\varepsilon}}'$. Suppose 
instead that $q_{\varepsilon} \in \partial D_{y_{\varepsilon}}'$. 
Observe that (3.17) gives 
\begin{equation} \label{e3.18}
\delta(q_{\varepsilon}) >> \rho'(q_{\varepsilon})^{2} 
\delta(y_{\varepsilon}).
\end{equation} 
Hence, if $y_{\varepsilon}$ is moved slightly to $y_{\varepsilon}'$ so 
that $q_{\varepsilon}' \in D_{y_{\varepsilon}}'$, then (3.18) is still 
preserved. It follows that the limit harmonic potential $\nu_{q} = \lim 
\hat \nu_{\delta(q_{\varepsilon})}$ as above has a local maximum at 
$q_{\varepsilon}$, giving again a contradiction. This completes the 
proof of Lemma 3.2.
{\endproof}

 By construction, on the maximal limit $(F, g_{o}', x)$, one has 
$B_{x}'(1) \subset  F$, c.f. also (2.19). Assume for the moment that 
$F$ is flat, so that $\partial F\cap\partial B_{x}'(1) \neq  
\emptyset$. If $B_{x}'(1)$ is simply connected, then the developing map 
$\mathcal{D} $ of the flat structure on $F$, c.f. [17], gives an 
isometric embedding of $B_{x}'(1)$ onto the unit ball $B_{0}(1) \subset 
 {\mathbb R}^{3},$ where $\mathcal{D} $ is normalized so that 
$\mathcal{D} (x) =$ 0. For any fixed $R_{o} > $ 1, let 
\begin{equation} \label{e3.19}
D(R_{o}) \subset  B_{x}'(R_{o}) \subset  \bar F 
\end{equation}
be a maximal domain in $B_{x}'(R_{o}),$ containing $B_{x}'(1),$ such 
that $\mathcal{D} $ is an embedding of $D(R_{o})$ onto its image in 
${\mathbb R}^{3}.$ Such a maximal domain is not unique, (as for 
instance cuts in a Riemann surface are not unique), but the results 
below are independent of the particular choice of maximal domain.

 A similar definition holds if $B_{x}'(1)$ is not simply connected. 
Thus, if $\Gamma  = \pi_{1}(B_{x}'(1)),$ then $\Gamma $ injects in 
$\pi_{1}(F)$ since the unit ball $\widetilde B_{\widetilde x}'(1)$ is 
convex in the universal cover $\widetilde F$ of $F$. The group $\Gamma$ 
is either ${\mathbb Z}$ or ${\mathbb Z} \oplus {\mathbb Z}$, acting by 
translation or twist on $\widetilde F$ or ${\mathbb R}^{3}$. The 
developing map descends to a map $\mathcal{D} : \widetilde F/\Gamma  
\rightarrow  {\mathbb R}^{3}/\Gamma ,$ and $\mathcal{D} $ induces an 
embedding of $B_{x}'(1)$ onto the unit ball in ${\mathbb R}^{3}/\Gamma 
.$ Thus, as in (3.19), define $D(R_{o})$ to be a maximal domain on 
which $\mathcal{D} $ is an embedding onto its image in ${\mathbb 
R}^{3}/\Gamma .$ 

 The reason for constructing the domains $D(R_{o})$ is that they are 
natural {\it compact} analogues of the geodesic balls $B_{x}'(R_{o})$; 
as noted preceding Remark 2.2, the balls $B_{x}'(R_{o})$ may not be 
compact, (or have compact closure), in general. The same definition 
holds if the limit $F$ is hyperbolic, using the developing map 
$\mathcal{D} $ of the hyperbolic structure. Domains $D_{y}(R_{o}),$ 
centered at a point $y\in F,$ may be defined in the same way. 

 The next two results give (local) upper bounds on the mass of 
$\nu_{a}$; these will be important in \S 4, (in Lemma 4.6).

\begin{proposition} \label{p 3.3.}
  For $(F, g_{o}', x, \nu_{a})$ as in Lemma 3.2, let $\bar L$ be any 
level set of the potential $\nu_{a}.$ Then there is a constant $M_{o},$ 
depending only on $R_{o},$ such that
\begin{equation} \label{e3.20}
m(D(R_{o})\cap\bar L) \leq  M_{o}, 
\end{equation}
where $m$ is the mass of the limit potential $\nu_{a}$ on the level set 
$\bar L.$
\end{proposition}

\noindent
{\bf Proof:}
 Assume the limit $F$ is flat; the proof in the hyperbolic case is the 
same. Thus, the domain $D(R_{o})$ may be viewed as a domain in 
${\mathbb R}^{3},$ (or ${\mathbb R}^{3}/\Gamma ).$

 As noted following Theorem 2.11, the function $\nu_{a}$ is uniformly 
controlled in $B_{x}'(\frac{1}{2})$. The function $\nu_{a}$ satisfies 
the elliptic equation (2.65) and is bounded above near $\partial F$ by 
Lemma 3.2. Standard properties of positive (or negative) solutions of 
elliptic equations in domains in ${\mathbb R}^{3}$ imply that $\nu_{a}$ 
is locally uniformly controlled, in that it satisfies a uniform local 
Harnack inequality and gradient estimates off $\partial F$, c.f. [9, 
Ch.8]. This implies that (3.20) holds for regions of $\bar L \cap 
D(R_{o})$ not too close to $\partial F$.

 In general, let $\bar H$ be a compact region within $D(R_{o}) \cap 
\bar L$ and assume w.l.o.g. that $\nu_{a}(\bar L) < \nu_{a}(L)$, where 
$L$ is the level through $x$. The flow lines of $\nabla \nu_{a}$ 
starting in $\bar H$ then tend to a local maximum of $\nu_{a}$. (The 
small set of flow lines ending in critical points of $\nu_{a}$ may be 
ignored in the following). For those flow lines $\gamma$ for which 
$dist(\gamma(s), \partial F) \geq t_{o}$, for some small $t_{o} > 0$, 
the arguments above give a uniform upper bound on $|\nabla 
\nu_{a}|(\gamma(s))$. Suppose instead a flow line $\gamma$ stays near 
$\partial F$. Then since $\nu_{a}$ is bounded above, $|\nabla 
\nu_{a}|(\gamma(s))$ is small, sufficiently far out along the flow. 
Hence, the region $R$ consisting of the union of these flow lines which 
terminate where $|\nabla \nu_{a}| \leq C$, for some fixed constant $C$, 
is a compact set in $F$. Then (3.20) follows by applying the divergence 
theorem to the equation (2.65) over the domain $K$.
{\endproof}

 Proposition 3.3 also holds on the sequence $(\Omega_{\varepsilon}, 
g_{\varepsilon}' , x_{\varepsilon})$ converging to $(F, g_{o}', x)$. 
Thus, let $D_{\varepsilon}(R_{o})$ be a sequence of domains in 
$B_{x_{\varepsilon}}'(R_{o}) \cap U^{\upsilon_{o}},$ for 
$U^{\upsilon_{o}}$ as in (1.18), converging in the Gromov-Hausdorff 
topology to the domain $D(R_{o})$ as in (3.19). In addition, 
$D_{\varepsilon}(R_{o})$ may be chosen to that if 
$\partial_{s}D_{\varepsilon}(R_{o})$ is any subset of $\partial 
D_{\varepsilon}(R_{o})$ converging to a proper subset of $\partial F 
\cap \bar D(R_{o})$, then $\partial_{s}D_{\varepsilon}(R_{o}) \subset 
L^{\upsilon_{o}}$. These conditions do not uniquely define 
$D_{\varepsilon}(R_{o})$, (even topologically), but this plays no role. 
If the sequence $(\Omega_{\varepsilon}, g_{\varepsilon}' , 
x_{\varepsilon})$ is collapsing, or is $\mu_{o}$-collapsed for 
$\varepsilon$ sufficiently small, define $D_{\varepsilon}$ so that in 
addition $\pi_{1}(B_{x_{\varepsilon}}'(1)) \subset 
\pi_{1}(D_{\varepsilon}(R_{o}))$, so that coverings unwrapping the 
collapse of $B_{x_{\varepsilon}}'(1)$ also unwrap 
$D_{\varepsilon}(R_{o})$. 

 Let $\bar L_{\varepsilon}$ be any level set of $u$ with $u(\bar 
L_{\varepsilon}) \rightarrow 1$ and let $m_{\delta_{a}}'$ be the mass 
in the scale $g_{\varepsilon}'$ of $\nu_{\delta_{a}(x_{\varepsilon})}$, 
lifted to covers if $g_{\varepsilon}'$ is sufficiently collapsed at 
$x_{\varepsilon}$.
\begin{corollary} \label{c 3.4.}
  There is a constant $M_{o} <  \infty ,$ such that on the levels $\bar 
L_{\varepsilon}$ as above,
\begin{equation} \label{e3.21}
m_{\delta_{a}}'(D_{\varepsilon}(R_{o})\cap\bar L_{\varepsilon}) \leq  
M_{o}. 
\end{equation}
\end{corollary}

\noindent
{\bf Proof:}
 The proof is identical to the proof of Proposition 3.3; alternately, 
(3.21) follows from (3.20) and the continuity of the mass to the limit.
{\endproof}

 The results above give local upper bounds for the mass of the 
potential, within small $g_{\varepsilon}$-neighborhoods of base points 
$x_{\varepsilon}.$ The following elementary result gives a crucial 
upper bound on the {\it  total}  mass, on the base scale 
$(\Omega_{\varepsilon}, g_{\varepsilon}).$

\begin{theorem} \label{t 3.5.}
  Let $L$ be any level set of u, $\delta  > $ 0 and set $\nu_{\delta} = 
(u-1)/\delta .$ Then on $(\Omega_{\varepsilon}, g_{\varepsilon}),$
\begin{equation} \label{e3.22}
m_{\delta}(L) \equiv  \int_{L}|\nabla\nu_{\delta}|dA_{g_{\varepsilon}} 
\leq  \frac{\bar \varepsilon\mathcal{Z}^{2}}{2\delta}. 
\end{equation}
\end{theorem}

\noindent
{\bf Proof:}
 Let $U^{+} = \{x\in\Omega_{\varepsilon}: u(x) \geq  u(L)\}$ and apply 
the divergence theorem to the trace equation (2.3), renormalized by 
$\delta$, over $U^{+}.$ A little algebra then gives
\begin{equation} \label{e3.23}
\int_{L}|\nabla\nu_{\delta}| = \int_{U^{+}} \bigl( 
\frac{\varepsilon\mathcal{Z}^{2}}{4\delta T} -  
\frac{1}{8}\frac{\sigma}{T} ( \frac{T^{2}u^{2} -  1}{\delta} ) -  
\frac{\varepsilon |z|^{2}}{4\delta T} \bigr)  \leq   \int_{U^{+}} 
\bigl( \frac{\varepsilon \mathcal{Z}^{2}}{4\delta T} -  
\frac{1}{8}\frac{\sigma}{T} ( \frac{T^{2}u^{2} -  1}{\delta}) \bigr) . 
\end{equation}
The left side of (3.23) is the total mass of the potential 
$\nu_{\delta}$ at this base scale $g_{\varepsilon}.$ Recall from \S 1 
that $Tu = w$ and $w = - s/\sigma ,$ where $\sigma $ is the $L^{2}$ 
norm of $s$; thus, the $L^{2}$ norm of $w$ equals 1. Since 
$vol_{g_{\varepsilon}}\Omega_{\varepsilon} =$ 1, this gives
$$\int_{\Omega_{\varepsilon}}(\frac{T^{2}u^{2} -  1}{\delta}) = 0. $$
Also, $T^{2}u^{2} - 1 \leq 0$ precisely on the set $W_{1} = \{w \leq 
1\}$. Hence
$$-\int_{U^{+}}(\frac{T^{2}u^{2} -  1}{\delta}) \leq  
-\int_{W_{1}}(\frac{T^{2}u^{2} -  1}{\delta}) = 
\int_{W^{1}}(\frac{T^{2}u^{2} -  1}{\delta}), $$
where $W^{1} = \{w \geq $ 1\}. Since $vol W^{1} \leq $ 1 and $u \leq $ 
1, it follows that
\begin{equation} \label{e3.24}
-\frac{1}{8}\frac{\sigma}{T}\int_{U^{+}}(\frac{T^{2}u^{2} -  
1}{\delta}) \leq  \frac{1}{8}\frac{\sigma}{T}(\frac{T^{2} -  
1}{\delta}) \leq  \frac{1}{4}\frac{\varepsilon\mathcal{Z}^{2}}{\delta 
T}, 
\end{equation}
where the last estimate follows from (1.8). This gives the bound (3.22).
{\endproof}

\begin{remark} \label{r 3.6.}
  Here a remark, not actually relevant to this paper, but relevant to 
the Sphere conjecture and the asymptotic geometry of 
$\mathcal{Z}_{c}^{2}$ solutions discussed in [1,2]. All of the analysis 
in \S 2-\S 3 carries over to the linearization at infinity of 
$\mathcal{Z}_{c}^{2}$ solutions $(N, g')$, provided $\sup_{N}u <  
\infty .$ In fact, this situation is considerably easier, since $su 
\equiv 0$ and there is no constant term in the Euler-Lagrange equations 
(1.14)-(1.15). In particular, one can define mass functions similar to 
(3.9), c.f. [2, (7.54)ff]. However, in this case, there is no upper 
bound estimate available for the mass - Theorem 3.5 is lacking. It is 
precisely the constant term $c = c_{\varepsilon}$ of (1.6) in the 
Euler-Lagrange equations which leads to this bound. This information, 
which will be seen to be crucial, is thus lost when passing to 
$\mathcal{Z}_{c}^{2}$ limits.
\end{remark}

\section{Uniformity Estimates.}
\setcounter{equation}{0}

 All of the results in \S 2 and \S 3, with the exception of Theorem 
3.5, are local; they concern the geometry of blow-ups $g_{\varepsilon}' 
$ and potentials $\nu_{\delta}$ in small neighborhoods of base points 
$x_{\varepsilon}.$ All of these results depend strongly on the size of 
$\delta ,$ and hence on the choice of base point $x_{\varepsilon}.$

 In order to relate these local behaviors at different base points, as 
well as relate them to the global estimate (3.22), it is important to 
obtain uniform control on $\delta$ on a given level set of $u$. 

 The purpose of this section is to prove the existence of base points 
$y_{\varepsilon}$ where $\delta ,$ (or more precisely $\delta_{a}),$ is 
not too small, both at $y_{\varepsilon}$ and at points on the same 
level set as $y_{\varepsilon}$ and within, say, unit 
$g_{\varepsilon}$-distance to $y_{\varepsilon}.$ The main result in 
this section is Theorem 4.4. Some applications to the structure of 
$\partial F$ are given in Proposition 4.8-Lemma 4.9.

\medskip

 To begin, the analysis is divided into two cases, according to the 
following:
\begin{definition} \label{d 4.1.}
 The sequence $(\Omega_{\varepsilon}, g_{\varepsilon}), \varepsilon  = 
\varepsilon_{i},$ satisfies the {\sf degeneration hypothesis}  if there 
is a sequence of base points $x_{\varepsilon}\in\Omega_{\varepsilon}$ 
such that, as $\varepsilon  \rightarrow $ 0,
\begin{equation} \label{e4.1}
|u- 1|(x_{\varepsilon}) \leq  \frac{1}{(\ln \bar \varepsilon)^{2}}, \ \ 
{\rm and} \ \ \rho (x_{\varepsilon}) \rightarrow  0. 
\end{equation}
The sequence $(\Omega_{\varepsilon}, g_{\varepsilon})$ satisfies the 
{\sf non-degeneration hypothesis} if there exists a constant $\rho_{o} 
> $ 0 such that whenever $|u- 1|(x_{\varepsilon}) \leq  1/(\ln \bar 
\varepsilon)^{2}$, then
\begin{equation} \label{e4.2}
\rho (x_{\varepsilon}) \geq  \rho_{o}. 
\end{equation}
\end{definition}

 Thus, the given sequence $(\Omega_{\varepsilon}, g_{\varepsilon}), 
\varepsilon  = \varepsilon_{i}$, satisfies either the degeneration or 
the non-degeneration hypothesis; if it satisfies the degeneration (or 
non-degeneration) hypothesis, then so do all subsequences. Of course if 
$(\Omega_{\varepsilon}, g_{\varepsilon})$ satisfies the degeneration 
hypothesis, then there still may well exist other base points 
$z_{\varepsilon}\in\Omega_{\varepsilon},$ distinct from 
$x_{\varepsilon},$ such that $|u- 1|(z_{\varepsilon}) \leq  1/(\ln \bar 
\varepsilon)^{2}$ but $\rho (z_{\varepsilon}) \geq  \rho_{o},$ for some 
$\rho_{o} > $ 0. The specific decay requirement (4.1) on $|u- 1|$ is 
used only in the proof of Proposition 4.8 below and for most purposes 
can be replaced by the weaker condition that $u(x_{\varepsilon}) 
\rightarrow $ 1 at any rate as $\varepsilon  \rightarrow $ 0.

 The analysis in these two cases is quite different, and in particular 
is much simpler when $(\Omega_{\varepsilon}, g_{\varepsilon})$ 
satisfies the non-degeneration hypothesis. The non-degeneration case 
(4.2) will be analysed in \S 7, and throughout \S 4-\S 5, it is assumed 
that $(\Omega_{\varepsilon}, g_{\varepsilon})$ satisfies the 
degeneration hypothesis. (It will be used however only in Theorem 4.4). 

 Further, in this section we use the assumption that $M$ is 
$\sigma$-tame. Thus, it will always be assumed that the sequence 
$(\Omega_{\varepsilon}, g_{\varepsilon}), \varepsilon  = 
\varepsilon_{i},$ is chosen to satisfy (1.36), i.e.
\begin{equation} \label{e4.3}
\bar \varepsilon\mathcal{Z}^{2} \equiv \bar 
\varepsilon\mathcal{Z}^{2}(g_{\varepsilon}) \leq  \bar 
\varepsilon^{4\mu}, 
\end{equation}
for some $\mu  <  \frac{1}{4}$ and $\varepsilon $ sufficiently small; 
(4.3) follows from (1.36) since $T \geq $ 1 and hence $T^{\mu} \leq $ 
T. Of course, (4.3) then remains valid for any subsequences.

\medskip

 Now for any given $\varepsilon  > $ 0 fixed, consider certain special 
base points as follows:

\begin{definition} \label{d 4.2.}
 A sequence of base points $x_{\varepsilon}$ is called {\sf admissible} 
 if it satisfies $\rho (x_{\varepsilon}) \rightarrow $ 0 together with 
the following conditions. Let $L = L_{\varepsilon}$ be the level set of 
$u$ containing $x_{\varepsilon}.$ Then
\begin{equation} \label{e4.4}
 \bar \varepsilon^{2\mu} \leq  |u(L)- 1| \leq  \frac{2}{(\ln \bar 
\varepsilon)^{2}}, 
\end{equation}
and, for the affine $\delta_{a}$ defined as in (3.1) and any point 
$p_{\varepsilon}\in B_{x_{\varepsilon}}(1)\cap L$,
\begin{equation} \label{e4.5}
\delta_{a}(p_{\varepsilon}) \geq  \delta_{o} \equiv  \bar 
\varepsilon^{2\mu}. 
\end{equation}
\end{definition}

 Observe that by Definition 2.7, admissible base points are allowable. 
Of course, $\delta_{o} \rightarrow $ 0 as $\varepsilon  \rightarrow $ 0 
and in fact, from (4.3),
\begin{equation} \label{e4.6}
\frac{\bar \varepsilon\mathcal{Z}^{2}}{\delta_{o}} \leq  \bar 
\varepsilon^{2\mu}\rightarrow  0, \ \ {\rm as} \ \ \varepsilon  
\rightarrow  0. 
\end{equation}
Hence $\bar \varepsilon\mathcal{Z}^{2}/ \delta (p_{\varepsilon}) 
\rightarrow $ 0, for all $p_{\varepsilon}\in L\cap 
B_{x_{\varepsilon}}(1);$ compare with (3.22). The estimates (4.5) and 
(4.6) play the central role in the arguments of \S 5. Although (4.5) is 
not required to hold globally on $L$, it is important to recognize that 
(4.5) holds on a region in $L$ of determined size about the base point 
$x_{\varepsilon}.$

 The definition of $\delta_{o}$ in (4.5), natural in view of (4.6), is 
just a convenient choice, and other choices, such as $\delta_{o} = \bar 
\varepsilon^{\mu},$ are equally possible. Nevertheless, it is necessary 
to make one choice, and so we choose (4.5).

\smallskip

 An important related class of base points are the following:

\begin{definition} \label {d 4.3.}
 A sequence of base points $x_{\varepsilon}$ is called {\sf 
well-separated}  (w.r.t. $U^{\upsilon_{o}}),$ if it satisfies (4.4) 
together with the following conditions: there exist scales 
$l_{\varepsilon} > $ 0 and constants $r_{o}, d_{o} > $ 0, (independent 
of $\varepsilon ),$ such that $l_{\varepsilon} \geq  r_{o}\rho 
(x_{\varepsilon}),$ (4.5) holds for all $p_{\varepsilon}\in 
B_{x_{\varepsilon}}^{\upsilon_{o}}(l_{\varepsilon})\cap L,$ and
\begin{equation} \label{e4.7}
A_{x_{\varepsilon}}^{\upsilon_{o}}(l_{\varepsilon}, 
(1+d_{o})l_{\varepsilon})\cap L = \emptyset  . 
\end{equation}
\end{definition}

 Here $A_{x}^{\upsilon_{o}}(r, s)$ is the geodesic annulus about $x$ of 
inner and outer radii $r, s$, in the metric space $(U^{\upsilon_{o}}, 
g_{\varepsilon}),$ where $U^{\upsilon_{o}} = \{x_{\varepsilon} \in 
\Omega: u(x_{\varepsilon}) \geq 1 - \upsilon_{o} \}$, as in Proposition 
1.1; similarly for the ball $B_{x}^{\upsilon_{o}}$. While it is 
possible to work with the (more natural) full annulus 
$A_{x_{\varepsilon}}$ in $(\Omega_{\varepsilon}, g_{\varepsilon}),$ it 
is somewhat simpler at this stage to work with the restricted annulus 
$A_{x_{\varepsilon}}^{\upsilon_{o}};$ it will be seen later (in Lemma 
5.7 below) that the difference between these two choices becomes 
irrelevant.

 The region $\mathcal{C}  = \mathcal{C} (\varepsilon ) = 
B_{x_{\varepsilon}}^{\upsilon_{o}}(l_{\varepsilon})\cap L$ gives rise 
to a {\it  compact}  and {\it  isolated}  region of the level set $L$ 
of $\nu_{a}$ in the limit $(F, g_{o}', x)$, where $g_{\varepsilon}'  = 
l_{\varepsilon}^{-2}g_{\varepsilon};$ compare with Remark 2.2. A 
well-separated base point satisfies all the properties of admissible 
base points, except that (4.5) may not hold at all points in the 1-ball 
about $x_{\varepsilon}$, (and $\rho(x_{\varepsilon})$ is not required 
to go to 0).

 Finally, a sequence of base points $x_{\varepsilon}$ is called {\sf  
preferred}  if it is well-separated with $l_{\varepsilon} \sim  \rho 
(x_{\varepsilon}),$ i.e. $l_{\varepsilon}/\rho (x_{\varepsilon})$ is 
bounded away from 0 and infinity as $\varepsilon  \rightarrow $ 0.

\medskip

 The purpose of this section is to prove the existence of at least one 
of these two kinds of base points, either admissible or preferred, near 
any base points $x_{\varepsilon}$ satisfying (4.1). In \S 5, it will 
then be proved that if $y_{\varepsilon}$ is an admissible base point, 
then there exist well-separated base points $z_{\varepsilon}$ near 
$y_{\varepsilon}.$ As in Remark 2.2, it is the well-separated base 
points which are naturally associated with 2-spheres.

\begin{theorem} \label{t 4.4.}
  Suppose $(\Omega_{\varepsilon}, g_{\varepsilon})$ satisfies the 
degeneration hypothesis (4.1), with $\varepsilon  = \varepsilon_{i}$ 
satisfying (4.3). Then for any base points $x_{\varepsilon}$ satisfying 
(4.1), there are base points $y_{\varepsilon}\in 
B_{x_{\varepsilon}}(1)$ which are either admissible or preferred.
\end{theorem}
 The proof will proceed in severals steps, via two Lemmas. Let $L = 
L(\varepsilon )$ be the the level set of $u$ through $x_{\varepsilon}.$ 
By Proposition 1.1, $t_{\upsilon_{o}}(x_{\varepsilon}) \rightarrow $ 0, 
where $t_{\upsilon_{o}} = dist(\cdot  , L^{\upsilon_{o}}).$ Let $L_{1} 
= L_{1}(\varepsilon )$ be the level set of $u$ defined by
\begin{equation} \label{e4.8}
u(L_{1}) = u(L) -  \bar \varepsilon^{2\mu}, 
\end{equation}
so that $u(L_{1}) \rightarrow $ 1, as $\varepsilon  \rightarrow $ 0 
also. By the converse of Proposition 1.1, c.f. (2.18), it follows that 
there are base points $y_{\varepsilon}\in 
B_{x_{\varepsilon}}(\frac{1}{10})\cap L_{1}$ such that $\rho 
(y_{\varepsilon}) \rightarrow $ 0 as well. Thus, $y_{\varepsilon}$ also 
satisfies (4.4).

 It will be proved that such base points $y_{\varepsilon}$ are either 
admissible, or there exist preferred base points $q_{\varepsilon}\in 
B_{y_{\varepsilon}}(1),$ i.e. either all $z_{\varepsilon}\in 
B_{y_{\varepsilon}}(1)\cap L_{1}$ satisfy (4.5) or there exist 
preferred base points $q_{\varepsilon}\in B_{y_{\varepsilon}}(1).$ Note 
that the Lipschitz property of $\rho $ implies that $\rho 
(z_{\varepsilon}) \leq $ 1, for all $z_{\varepsilon}\in 
B_{y_{\varepsilon}}(1).$ 

 Let $L_{2}$ be the level set of $u$ defined by
\begin{equation} \label{e4.9}
u(L_{2}) = u(L_{1}) -  \bar \varepsilon^{\mu}. 
\end{equation}
Consider arbitrary (piecewise) $C^{1}$ paths $\gamma $ in 
$U^{\upsilon_{o}}$, which start at some point in $L^{\upsilon_{o}},$ 
end at a point $z_{\varepsilon}$ in $B_{y_{\varepsilon}}(1)\cap L_{1},$ 
and which satisfy the following cone condition:
\begin{equation} \label{e4.10}
s \leq  C\cdot  dist(\gamma (s), L^{\upsilon_{o}}) = C \cdot 
t_{\upsilon_{o}}(\gamma(s)), 
\end{equation}
where $s$ is the arclength parameter of $\gamma $ and $C$ is any fixed 
constant. Observe that (4.10) is scale-invariant. It is clear that such 
paths exist; for instance one may choose $\gamma $ to be a shortest 
path in $U^{\upsilon_{o}}$ realizing the distance of $z_{\varepsilon}$ 
to $L^{\upsilon_{o}}.$ The condition (4.10) roughly means that $\gamma 
$ does not zigzag too much across the levels of $u$. For any such 
$\gamma ,$ let $s_{2}$ be the largest value of $s$ such that $\gamma 
(s_{2})\in L_{2}$ and let $a_{\varepsilon} = \gamma (s_{2}).$ Finally 
set $\gamma^{2} = \gamma|_{[s_{2}, b]},$ where $\gamma (b) = 
z_{\varepsilon}.$ 

\smallskip

 A simple but important idea in the proof is the following:

\begin{lemma} \label{l 4.5.}
  For any curve $\gamma $ as above, there exist points 
$p_{\varepsilon}\in\gamma^{2}$ such that
\begin{equation} \label{e4.11}
\delta_{a}(p_{\varepsilon}) \geq  \delta_{o}. 
\end{equation}
\end{lemma}

\noindent
{\bf Proof:}
 If not, then there is a curve $\gamma $ satisfying (4.10) such that, 
$\forall s \geq  s_{2}$, $\delta_{a}(\gamma^{2}(s)) \leq  \delta_{o},$ 
and hence by Lemma 3.1,
\begin{equation} \label{e4.12}
|du|(\gamma^{2}(s)) \leq  C_{o}\frac{\delta_{o}}{\rho (\gamma^{2}(s))}. 
\end{equation}
Observe again that the estimate (4.12) is scale-invariant. For the 
argument to follow, we work in the (base) scale $g_{\varepsilon}.$ For 
$t_{\upsilon_{o}}$ as above, Proposition 1.1 implies that
$$\frac{1}{\rho (\gamma^{2}(s))} \leq  \frac{c_{1}}{t_{\upsilon_{o}}} 
\leq  \frac{c_{2}}{s}, $$
where the $2^{\rm nd}$ inequality follows from (4.10). From (4.12), it 
then follows that, along $\gamma^{2}(s),$
$$|\frac{du}{ds}|  \leq  c_{3}\frac{\delta_{o}}{s}. $$
Integrating this along $\gamma^{2}$ gives
\begin{equation} \label{e4.13}
|u(L^{2}) -  u(L^{1})| \leq  c_{4}\delta_{o} \ln(\frac{1}{s_{2}}).
\end{equation}

 Now $s_{2} \geq  dist_{g_{\varepsilon}}(a_{\varepsilon}, 
L^{\upsilon_{o}}) \equiv  \lambda_{o}$ and again by (2.18), 
$\lambda_{o} \geq  R_{o}^{-1}\cdot \rho (a_{\varepsilon}).$ By 
Proposition 1.4, $\rho^{2} \geq  \kappa\cdot \bar \varepsilon$ at any 
base point and so in particular $\lambda_{o} \geq  \bar \varepsilon.$ 
Thus, using the definition of $\delta_{o},$ one has
\begin{equation} \label{e4.14}
|u(L^{2}) -  u(L^{1})| \leq  c_{5}\bar \varepsilon^{2\mu} \ln (\bar 
\varepsilon^{-1}) <<  \bar \varepsilon^{\mu}, 
\end{equation}
for $\varepsilon $ sufficiently small. This contradicts (4.9), and so 
proves the result.

{\endproof}

 Observe that the proof shows there are points $q_{\varepsilon}$ on 
$\gamma^{2}$ such that $\delta_{a}(q_{\varepsilon}) >>  \delta_{o}.$ 
More importantly, the proof also holds if $\gamma $ is a curve 
satisfying (4.10), with, for instance, 
\begin{equation} \label{e4.15}
C_{o} = C_{o}(\varepsilon) \leq  \varepsilon^{-\mu /2}, 
\end{equation}
in place of the $\varepsilon$-independent $C_{o}$ of (4.12).

 If all the points $p_{\varepsilon}$ from Lemma 4.5 may be chosen to be 
the end points $z_{\varepsilon}\in L_{1}\cap B_{y_{\varepsilon}}(1)$ of 
curves $\gamma ,$ then the result is proved. Thus, suppose not, so that 
there exists $z_{\varepsilon}\in L_{1}\cap B_{y_{\varepsilon}}(1)$ such 
that
\begin{equation} \label{e4.16}
\delta_{a}(z_{\varepsilon}) <  \delta_{o}. 
\end{equation}
We then construct a path $\gamma ,$ starting at $z_{\varepsilon},$ 
(i.e. in the reverse direction from previously), satisfying (4.10) and 
an analogue of (4.12). The point is to construct $\gamma $ with a 
controlled increase of $\delta_{a},$ i.e. the local mass, along 
$\gamma$, analogous to (4.15). Given the initial estimate (4.16), it is 
then shown that this path either terminates on $L_{2},$ giving a 
contradiction as in Lemma 4.5, or terminates on a preferred base point 
$q_{\varepsilon}.$ This will then complete the proof.

\smallskip

 The construction of $\gamma $ and the growth estimate for $\delta_{a}$ 
is by an inductive procedure. Thus, let $z_{\varepsilon} = 
z_{\varepsilon}^{1},$ and let $g_{\varepsilon}^{1} = \rho 
(z_{\varepsilon}^{1})^{-2}\cdot  g_{\varepsilon}$ be the rescaled 
metric based at $z_{\varepsilon}^{1}.$ Let $\delta_{1} = 
\delta_{a}(z_{\varepsilon}^{1})$ and let $\nu_{\delta_{1}}$ be the 
corresponding potential function. By adding a suitable constant, one 
may assume that $\nu_{\delta_{1}}(L_{1}) = - 1;$ recall that adding 
constants to $\nu_{\delta_{a}}$ does not affect the measure or mass of 
the potential.  Thus 
\begin{equation} \label{e4.17}
\nu_{\delta_{1}}(L_{2}) \sim  -\bar \varepsilon^{\mu}/\delta_{1} <<  - 
1, 
\end{equation}
for $\varepsilon $ small. In particular, it follows that $\rho^{1} 
\rightarrow $ 0 on $L_{2},$ where $\rho^{1}$ is the 
$g_{\varepsilon}^{1}$ curvature radius, and so $L_{2}$ converges, in 
the Hausdorff topology, to $\partial F,$ where $(F, g_{o}', z^{1})$ is 
the maximal limit based at $z_{\varepsilon}^{1}.$

 Define intermediate level sets $L^{j}$ by
\begin{equation} \label{e4.18}
\nu_{\delta_{1}}(L^{j+1}) = e^{k_{o}j}\nu_{\delta_{1}}(L^{j}), 
\end{equation}
with $L^{1} = L_{1}.$ Here $k_{o}$ is a fixed number which will be 
specified below, (see (4.23) and (4.26)). It follows that for 
$\varepsilon  > $ 0 fixed, the level $L_{2}$ is reached after at most 
$N$ iterations, i.e. $L_{2} \sim  L^{N},$ where
\begin{equation} \label{e4.19}
N = \frac{1}{k_{o}} \ln \frac{\bar \varepsilon^{\mu}}{\delta_{1}}. 
\end{equation}
The induction will be on $j$, i.e. down the level sets $L^{j},$ 
beginning with the levels $L^{1}$ and $L^{2}.$

\smallskip

 To start, recall the compact domains $D_{\varepsilon}(R_{o})$ defined 
following Proposition 3.3, with base point $z_{\varepsilon}^{1}$. If 
$B_{z_{\varepsilon}^{1}}^{1}(1)$ is $\mu_{o}$-collapsed, i.e. 
$\omega(z_{\varepsilon}^{1}) \leq \mu_{o}\rho(z_{\varepsilon}^{1})$, 
then the collapse is unwrapped in covering spaces, as described 
following Proposition 3.3. Here $\mu_{o}$ is a fixed parameter, e.g. 
$\mu_{o} = 1/100$. For convenience, set $R_{o} = 2$. Then Corollary 3.4 
gives
\begin{equation} \label{e4.20}
m(D_{\varepsilon}(2)\cap\bar L) \leq  M_{o}, 
\end{equation}
for any level set $\bar L$ of $u$, e.g. $\bar L = L^{1}$ or $L^{2}$, 
where $m = m_{\delta_{1}}$ is the mass w.r.t. $\delta_{1}$ in the scale 
$g_{\varepsilon}^{1}$.

\smallskip

 The basic idea in the inductive step, i.e. in the construction of the 
next base point $z_{\varepsilon}^{2}$ on the level $L^{2}$ is to 
examine the distribution of the measure $d\mu_{\delta_{1}}$ within 
$L^{2}\cap D_{\varepsilon}(2).$ By way of illustration, suppose $d\mu $ 
is a measure on $I = [0,1]$ of mass at most $M_{o}.$ If the support of 
$d\mu $ satisfies supp $d\mu  = I$, then it is clear that for all 0 $<  
t < $ 1, there exists $x\in$ supp$d\mu$ such that 
\begin{equation} \label{e4.21}
m_{d\mu}([x, x+t]) \leq  M_{o}\cdot  t, 
\end{equation}
(i.e. the analogue of (4.20) holds on the scale $t$). If on the other 
hand the mass of $d\mu $ is concentrated, (for example the Dirac 
measure at a point in the most extreme case), then (4.21) does not 
hold. The only situation in which (4.21) can fail is when there are 
gaps in supp $d\mu ,$ i.e. intervals $J\subset I$ on which 
$J\cap$supp$d\mu  = \emptyset$. This corresponds to the existence of 
well-separated base points.

\smallskip

 Define a level set $\bar L,$ for instance $\bar L = L^{2}$, to be {\it 
$d_{o}$-separated} (w.r.t. $U^{\upsilon_{o}}),$ if there is a point 
$p_{\varepsilon}\in\bar L$ and radius $s_{\varepsilon}$ such that
\begin{equation} \label{e4.22}
A_{p_{\varepsilon}}^{\upsilon_{o}}(s_{\varepsilon}, 
(1+d_{o})s_{\varepsilon}) \cap  \bar L = \emptyset  , 
\end{equation}
compare with (4.7). The condition (4.22) is scale-invariant, although 
we continue to work in the scale $g_{\varepsilon}^{1}$. The relation 
(4.22) implies that $s_{\varepsilon} \geq  
diam_{g_{\varepsilon}^{1}}C_{o},$ where $C_{o}$ is the component of 
$\bar L$ containing $p_{\varepsilon}$, and $diam$ is the extrinsic 
diameter, i.e. the diameter within $(U^{\upsilon_{o}}, 
g_{\varepsilon}^{1})$. Let $t(x) = t^{1}(x) = 
dist_{g_{\varepsilon}^{1}}(x, L^{\upsilon_{o}})$. Henceforth, assume 
that $k_{o}$ is chosen sufficiently large so that
\begin{equation} \label{e4.23}
t^{1}(x) \leq  \tfrac{1}{4}\mu_{o}, \forall x\in L^{2}\cap 
D_{\varepsilon}(2), 
\end{equation}
where $\mu_{o}$ is the collapse parameter, defined preceding (4.20).

\begin{lemma} \label{l 4.6.}
  Suppose the level $L^{2}$ has no $d_{o}$-separated components 
intersecting $D_{\varepsilon}(2)$, for $d_{o} \leq 1$. Then there exist 
base points $z_{\varepsilon}^{2}\in  L^{2}\cap D_{\varepsilon}(2)$ such 
that
\begin{equation} \label{e4.24}
m(L^{2}\cap D_{z_{\varepsilon}^{2}}(2t(z_{\varepsilon}^{2}))) \leq  
\mu_{1}M_{o}(1+d_{o})\cdot  t(z_{\varepsilon}^{2}), 
\end{equation}
where $\mu_{1} = \mu_{1}(\mu_{o}).$
\end{lemma}

\noindent
{\bf Proof:}
 Suppose first that $\hat L^{2} = L^{2}\cap D_{\varepsilon}(2)$ is path 
connected. Since $L^{2}$ has no $d_{o}$-separated components 
intersecting $D_{\varepsilon}(2)$, the diameter, (within 
$U^{\upsilon_{o}}$), satisfies $diam_{g_{\varepsilon}^{1}} \hat L^{2} 
\geq \mu_{o} / (1+d_{o}) \geq \mu_{o}/2$. A curve $\tau $ realizing the 
intrinsic diameter of $\hat L^{2}$ thus has length $l(\tau) \geq 
\mu_{o}/2$. 

 Now cover the path $\tau$ by a collection of smaller domains with 
disjoint interiors as follows. Let $s$ be the arclength parameter for 
$\tau ,$ with $p_{0} = \tau (0).$ Let $D_{1} = D_{p_{1}}(2t(p_{1})),$ 
where $p_{1}$ is the point on $\tau $ whose arclength distance to $\tau 
(0)$ equals $2t(p_{1});$ here $D_{1}$ is the domain as in (4.20), 
centered on $p_{1}.$ Let $s_{2}$ be the largest value of $s$ such that 
$\tau (s)\in D_{1}$ and set $p_{2} = \tau (s_{2});$ clearly $s_{2} >  
2t(p_{1}).$ Next let $D_{2} = D_{p_{3}}(2t(p_{3})),$ where $p_{3}$ is 
the point on $\tau $ whose arclength distance to $p_{2}$ equals 
$2t(p_{3}).$ Observe that $D_{1}\cap D_{2} = p_{2}.$ One then continues 
inductively in this way until the maximal $p_{k}$ is defined, for $k$ 
odd. It follows from (4.23) that the domain $V = \cup D_{j}$ covers at 
least $\frac{1}{2}$ the length $l(\tau)$ of $\tau$. 

 By construction, $m(V) = \sum m(D_{j})$. If for each $j$, $m(D_{j}) >  
\mu_{1}M_{o} t(p_{j}),$ then 
$$m(V) >  \mu_{1}M_{o}\sum t(p_{j}). $$
But by constuction, $\sum t(p_{j})$ is least $\frac{1}{2}l(\tau) \geq  
\frac{1}{4}\mu_{o}$. Thus, if $\mu_{1} >  4\mu_{o}^{-1},$ this 
contradicts (4.20), and hence there is some $j$ such that $m(D_{j}) 
\leq  \mu_{1}M_{o}t(p_{j})$. This gives (4.24) in this situation, with 
$z_{\varepsilon}^{2} = p_{j}$.

 If $\hat L^{2}$ is not path connected, let $\tau_{k}$ be a path 
realizing the intrinsic diameter of each path component $C_{k}$ of 
$\hat L^{2},$ and let $\tau $ be the union of the paths $\tau_{k}.$ 
Then perform the construction above on each component $\tau_{k}.$ Since 
there are no $d_{o}$-separated subsets, for any $k$, there must be a 
$C_{j}$ and point $q_{j}\in C_{j}$ such that 
$dist_{D_{\varepsilon}(2)}(q_{j}, C_{k}) \leq  
\frac{1+d_{o}}{2}(diam_{D_{\varepsilon}(2)}C_{k}).$ It follows that 
$l(\tau ) \geq  (1+d_{o})^{-1}\mu_{o},$ and the same argument as above 
then gives (4.24).

{\endproof}

 The domain $D_{\varepsilon}^{1} \equiv D_{\varepsilon}(2)$ converges, 
(in a subsequence), to a limit domain $D^{1} \equiv D_{o}(2)$ contained 
in the maximal limit $(F, g_{o}', z^{1})$, $z^{1} = \lim 
z_{\varepsilon}^{1}$; the limit $F$ is either flat or hyperbolic. Since 
the limit level $L^{2}\cap D_{o}(2) \subset  F$ has a definite 
diameter, bounded away from 0 and $\infty ,$ it follows that $L^{2}\cap 
D_{\varepsilon}(2)$ is $d_{o}$-separated only when (4.22) holds with
\begin{equation} \label{e4.25}
s_{\varepsilon} \sim  1 
\end{equation}
in the scale $g_{\varepsilon}^{1},$ i.e. $s_{\varepsilon}$ is bounded 
away from 0 and $\infty $ in this scale.

\medskip

 Lemma 4.6 provides the inductive step. Thus, fix any choice of $d_{o} 
\in (0, 1]$ and set
\begin{equation} \label{e4.26}
k_{o} = \ln[A_{o}\mu_{1}(1+d_{o})M_{o}], 
\end{equation}
where $A_{o} = a_{o}^{-1}$ is the area constant from (2.79). (Note that 
$k_{o}$ may be chosen to be large, by choosing $M_{o}$ large). Suppose 
$L^{2}$ has no $d_{o}$-separated components in $D_{\varepsilon}(2)$. 
Then Lemma 4.6 gives the existence of the next base point 
$z_{\varepsilon}^{2}\in L^{2}$ and domain $D_{\varepsilon}^{2} \equiv  
D_{z_{\varepsilon}^{2}}(2t(z_{\varepsilon}^{2}))$ satisfying (4.24). 
Observe that, by construction, 
$dist_{g_{\varepsilon}^{1}}(z_{\varepsilon}^{2}, z_{\varepsilon}^{1}) 
\leq $ 2 and $dist_{g_{\varepsilon}^{1}}(z_{\varepsilon}^{2}, L^{1}) 
\geq  \frac{1}{2},$ where the $2^{\rm nd}$ inequality follows from 
(4.23). Hence, one may choose a curve $\gamma_{1}(s)\subset  
U^{\upsilon_{o}}$ from $z_{\varepsilon}^{2}$ to $z_{\varepsilon}^{1},$ 
with arclength parameter $s$, such that
\begin{equation} \label{e4.27}
s \leq  c\cdot  dist_{g_{\varepsilon}^{1}}(\gamma_{1}(s), L^{2}), 
\end{equation}
i.e. the analogue of (4.10) holds. Further, by (3.7) one has, along 
$\gamma_{1},$
\begin{equation} \label{e4.28}
|du|(\gamma_{1}(s)) \leq  C\frac{\delta_{1}}{\rho (\gamma_{1}(s))}, 
\end{equation}
where $C$ is a fixed constant, possibly depending on $k_{o}$ but not on 
$\varepsilon .$

 Now at the endpoint $z_{\varepsilon}^{2}$ of $\gamma_{1},$ we claim 
that for $M_{1} = A_{o}\mu_{1}M_{o}(1+d_{o})$,
\begin{equation} \label{e4.29}
\delta_{2} \equiv  \delta_{a}(z_{\varepsilon}^{2}) \leq 
M_{1}\delta_{1}. 
\end{equation}
To see this, recall first that the limit based at $z_{\varepsilon}^{1}$ 
is flat or hyperbolic. As noted in (2.19ff), without loss of 
generality, one may assume that
\begin{equation} \label{e4.30}
\rho^{1}(z_{\varepsilon}^{2}) = (1+o(1))t(z_{\varepsilon}^{2}), 
\end{equation}
where $\rho^{1}$ is the curvature radius in the $g_{\varepsilon}^{1}$ 
scale. In the scale $g_{\varepsilon}^{2} = \rho 
(z_{\varepsilon}^{2})^{-2}\cdot  g_{\varepsilon}$ where 
$\rho^{2}(z_{\varepsilon}^{2}) =$ 1, by definition 
$\delta_{a}(z_{\varepsilon}^{2}) = \oint_{H(\frac{1}{2})}|\nabla u|,$ 
where $H(\frac{1}{2}) = B_{z_{\varepsilon}^{2}}(\frac{1}{2})\cap 
L^{2}.$ Further, in this scale $D_{\varepsilon}^{2} = 
D_{z_{\varepsilon}^{2}}(2)$ and hence (4.24) and (4.30) give 
$$m_{\delta_{1}}(D_{\varepsilon}^{2}) = \int_{D_{\varepsilon}^{2}\cap 
L^{2}}|\nabla u|/\delta_{1} \leq  \mu_{1}M_{o}(1+d_{o}). $$
Since $H(\frac{1}{2}) \subset  D_{\varepsilon}^{2},$ (again in the 
$g_{\varepsilon}^{2}$ scale), it follows that
$$\delta_{a}(z_{\varepsilon}^{2}) \leq  \bigl( 
\mu_{1}M_{o}(1+d_{o})/area_{g_{\varepsilon}^{2}}H(\tfrac{1}{2}) \bigr)/ 
\delta_{1}.$$
Via (2.79), this gives (4.29).

 We now repeat this construction at the level 2 in place of level 1. 
Thus, as above, at the base point $z_{\varepsilon}^{2}\in L^{2},$ 
rescale the metric based at $z_{\varepsilon}^{2},$ so that 
$g_{\varepsilon}^{2} = \rho (z_{\varepsilon}^{2})^{-2}\cdot  
g_{\varepsilon}$ and work with the potential $\nu_{\delta_{2}}$ in 
place of $\nu_{\delta_{1}}.$ From (4.24), (4.29) and (4.30), one has 
the bound, in the $g_{\varepsilon}^{2}$ scale,
\begin{equation} \label{e4.31}
m_{\delta_{2}}(D^{2}) \leq  M_{1}^{2}. 
\end{equation}
Assuming that there are no $d_{o}$-separated subsets within $\hat L^{3} 
= L^{3}\cap D_{\varepsilon}^{2},$ one may then construct again by Lemma 
4.6 the next level base point $z_{\varepsilon}^{3}$ and domain 
$D_{\varepsilon}^{3},$ as well as the next part of $\gamma ,$ i.e. 
$\gamma_{2},$ from $z_{\varepsilon}^{2}$ to $z_{\varepsilon}^{3}$ 
satisfying the analogue of (4.27). The estimate (4.28) translates to 
the estimate
$$|du|(\gamma_{2}(s)) \leq  CM_{1}^{2}\frac{\delta_{1}}{\rho 
(\gamma_{2}(s))}. $$
Thus, this process may be continued inductively until the level $L_{2}$ 
is reached, provided at each stage there are no $d_{o}$-separated 
subsets.

\smallskip

 Suppose then that this process continues to the level $L_{2},$ 
producing a curve $\gamma $ starting at $L_{2}$ and ending at 
$z_{\varepsilon}^{1}.$ Since each part $\gamma_{1}, \gamma_{2},$ etc. 
of $\gamma $ satisfies the analogue of (4.27), $\gamma $ satisfies the 
cone condition (4.10). Further, along $\gamma $ one has
\begin{equation} \label{e4.32}
\delta_{a}(\gamma (s)) \leq  C\cdot  M_{1}^{N}\delta_{1}, 
\end{equation}
where $N$ is given by (4.19). However, together with (4.26) this gives
$$\delta_{a}(\gamma (s)) \leq  C\cdot \bar \varepsilon^{\mu 
/2}\delta_{1}^{1/2} \leq  C\cdot \bar \varepsilon^{3\mu /2},$$ 
everywhere along $\gamma (s),$ where the $2^{\rm nd}$ inequality uses 
$\delta_{1} \leq  \bar \varepsilon^{2\mu}$ by (4.16); compare also with 
(4.15). The proof of Lemma 4.5 carries through as before, and again 
gives a contradiction. 

\smallskip

 Hence if (4.16) holds, (for some $z_{\varepsilon}),$ there must be a 
(first) level $L^{j-1},$ with base points $z_{\varepsilon}^{j-1}\subset 
L^{j-1},$ and domains $D_{\varepsilon}^{j-1}$ such that the next level 
$\hat L^{j} = L^{j}\cap D_{\varepsilon}^{j-1}$ is $d_{o}$-separated, 
i.e. there exist $z_{\varepsilon}^{j}\in \hat L^{j}$ such that (4.22) 
holds, with $s_{\varepsilon} \sim 1$, in the scale 
$g_{\varepsilon}^{j-1} = \rho (z_{\varepsilon}^{j-1})^{-2}\cdot  
g_{\varepsilon}$. Suppose that, for all $p_{\varepsilon}^{j}\in \hat 
L^{j}$,
\begin{equation} \label{e4.33}
\delta_{a}(p_{\varepsilon}^{j}) \geq  \min (\delta_{o}, 
\delta_{a}(z_{\varepsilon}^{j-1})).
\end{equation}
Then (2.78) implies that the potential function $\nu_{\delta_{a}^{j-1}} 
\equiv \nu_{\delta_{a}(z_{\varepsilon}^{j-1})}$ satisfies $|\nabla 
\nu_{\delta_{a}^{j-1}}|(p_{\varepsilon}^{j}) \geq \kappa_{o} / 
t(p_{\varepsilon}^{j})$ in the $g_{\varepsilon}^{j-1}$ scale. Since 
$\nu_{\delta_{a}^{j-1}}(L^{j}) /\nu_{\delta_{a}^{j-1}}(L^{j-1}) = 
e^{k_{o}}$ and $k_{o}$ is fixed, it follows that 
$\rho^{j-1}(z_{\varepsilon}^{j}) \geq \rho_{o} = \rho_{o}(\kappa_{o}) > 
0$. This means that $\rho^{j}(z_{\varepsilon}^{j}) \sim 1$, i.e. 
$s_{\varepsilon} \sim \rho(z_{\varepsilon}^{j})$ in the 
$g_{\varepsilon}$ scale, as $\varepsilon \rightarrow 0$. It follows 
that $z_{\varepsilon}^{j}$ is preferred, which completes the proof in 
this case.

 Suppose instead that for some $p_{\varepsilon}^{j}\in \hat L^{j}$,
$$\delta_{a}(p_{\varepsilon}^{j}) <  \max(\delta_{o}, 
\delta_{a}(z_{\varepsilon}^{j-1})).$$
In this case, the analogue of (4.16) holds, with $p_{\varepsilon}^{j}$ 
in place of $z_{\varepsilon}^{1}$ and one just continues the 
construction of the curve $\gamma $ joining $z_{\varepsilon}^{j-1}$ to 
$z_{\varepsilon}^{j} \equiv p_{\varepsilon}^{j}$. The point is that one 
must eventually come to points where (4.33) holds and the only way to 
do this is to come to a preferred base point, (unless $x_{\varepsilon}$ 
is admissible).

{\endproof}

 Theorem 4.4 implies that there may be many admissible or preferred 
base points $x_{\varepsilon}.$ In fact they are 1-dense in the region 
$S_{\varepsilon} \subset  \Omega_{\varepsilon}$ defined by (4.4), i.e. 
$\bar \varepsilon^{2\mu} \leq  |u-1| \leq  2/(\ln \bar 
\varepsilon)^{2},$ and $\rho  \leq  \chi_{\varepsilon},$ where 
$\chi_{\varepsilon}$ is any sequence with $\chi_{\varepsilon} 
\rightarrow $ 0 as $\varepsilon  = \varepsilon_{i} \rightarrow $ 0. It 
is obvious from the proof that such points are $\eta$-dense in 
$S_{\varepsilon},$ for any given $\eta  > $ 0. This lack of uniqueness 
plays no role however. Any level set $L$ of $u$ satisfying (4.4) and 
containing some admissible or preferred base point will be called an 
{\sf  admissible}  level set and denoted by
\begin{equation} \label{e4.34}
L = L_{o}. 
\end{equation}

 Let $\nu_{\delta_{o}} = (u-1)/\delta_{o},$ for $\delta_{o}$ as in 
(4.5), and let $m_{o}, \widetilde m_{o}$ be the mass of the potential 
$\nu_{\delta_{o}},$ as in (3.9) and (3.11) respectively. The following 
uniform lower bound on the mass is essentially immediate.

\begin{proposition} \label{p 4.7.}
  Let $x_{\varepsilon}$ be any admissible or preferred base point. Then 
there is a constant $\kappa_{o} > $ 0, independent of $\varepsilon  
\leq $ 1, such that for any $p_{\varepsilon}\in L_{o}\cap 
B_{x_{\varepsilon}}(1),$ or $p_{\varepsilon} = x_{\varepsilon}$ 
respectively,
\begin{equation} \label{e4.35}
\widetilde m_{o}(\rho (p_{\varepsilon})) \geq  \kappa_{o}\cdot \rho 
(p_{\varepsilon}). 
\end{equation}
In particular, if $(\Omega_{\varepsilon}, g_{\varepsilon})$ is not 
$\mu_{o}$ collapsed at $x_{\varepsilon}$, i.e. $\omega(x_{\varepsilon}) 
\geq \mu_{o} \cdot \rho(x_{\varepsilon})$, then
\begin{equation} \label{e4.36}
m_{o}(\rho(x_{\varepsilon})) \geq \kappa_{o} \cdot 
\rho(x_{\varepsilon}),
\end{equation}
where $\kappa_{o}$ depends only on $\mu_{o}$.
\end{proposition}

\noindent
{\bf Proof:}
  By construction, i.e. (4.5), $\delta_{a}(p_{\varepsilon}) \geq  
\delta_{o},$ and hence
$$|d\nu_{\delta_{a}(p_{\varepsilon})}| \leq  |d\nu_{\delta_{o}}|, $$
for all $p_{\varepsilon}$ as above. Hence, the mass defined w.r.t. 
$\delta_{o}$ dominates the mass defined w.r.t. 
$\delta_{a}(p_{\varepsilon}).$ The estimates (4.35) and (4.36) are then 
an immediate consequence of (3.13).

{\endproof}

 We close this section with some applications of the proof of Theorem 
4.4 to the structure of $\partial F;$ these will be useful in \S 5.3 
and \S 6. Let $(F, g_{o}', x)$ be a maximal limit, (flat or 
hyperbolic), of $(\Omega_{\varepsilon}, g_{\varepsilon}' , 
x_{\varepsilon}),$ with $x_{\varepsilon}$ allowable, and let $\nu $ be 
the limit potential, $\nu  = \lim \nu_{\delta_{a}}, \delta_{a} = 
\delta_{a}(x_{\varepsilon}).$ If $(\Omega_{\varepsilon}, 
g_{\varepsilon})$ is $\mu_{o}$-collapsed at $x_{\varepsilon},$ i.e. 
$\omega (x_{\varepsilon}) \leq  \mu_{o}\cdot \rho (x_{\varepsilon}),$ 
it is assumed that the collapse is unwrapped in covering spaces, so 
that $\omega (x_{\varepsilon}) \sim \rho (x_{\varepsilon})$. Define
\begin{equation} \label{e4.37}
I_{\infty} = \overline{\cap_{N}\{U^{-N}\}}, 
\end{equation}
where $U^{-N}$ is the $-N$ superlevel set of $\nu$, $N\in{\mathbb 
Z}^{+}$, and the closure is taken in $\bar F$. Let $d\mu_{-N}$ be the 
measure of $\nu$ on the level $L^{-N}$, as in (3.8), but restricted to 
$D_{x}(\infty)$, c.f. (3.19), i.e. the maximal domain on which the 
developing map is an embedding. Thus $d\mu_{-N}$ is supported on the 
level set $L_{-N}\cap D_{x}(\infty)$.

 Let $d\mu_{\infty}$ be the weak limit of the measures $d\mu_{-N}$ in 
$\bar F$ as $N \rightarrow  \infty .$ The limit exists since the 
measures $d\mu_{-N}$ are related via the divergence theorem applied to 
the trace equation (2.65). The measure $d\mu_{\infty}$ is supported on 
$\partial F.$

\begin{proposition} \label{p 4.8.}
   Suppose the base points $x_{\varepsilon}$ are allowable and satisfy 
$|u-1|(x_{\varepsilon}) \leq 2/ (\ln \bar \varepsilon)^{2}$, and let 
$(F, g_{o}' , x, \nu)$ be the maximal limit as above. Then
\begin{equation} \label{e4.38}
\partial F = I_{\infty}. 
\end{equation}
There is a constant $m_{o}$ such that if B(r) is any r-ball in $\bar F$ 
with center in $\partial F,$ then $\forall r \leq $ 1,
\begin{equation} \label{e4.39}
m_{o}\cdot  r \leq  m_{\infty}(B(r)\cap\partial F), 
\end{equation}
where $m_{\infty}$ is the mass of $d\mu_{\infty}.$ Further for any unit 
ball B(1) in $\bar F$ centered at a point in $\partial F,$
\begin{equation} \label{e4.40}
m_{\infty}(B(1)\cap\partial F) \leq  M_{o}, 
\end{equation}
for some constant $M_{o} < \infty$. The constants $m_{o}$ and $M_{o}$ 
are independent of the base point on $\partial F.$ In particular, 
$dim_{\mathcal{H}}\partial F \leq $ 1 and $\mathcal{H}_{1}(\partial 
F\cap B(1)) \leq  M_{o},$ where $\mathcal{H}_{1}$ is the 1-dimensional 
Hausdorff measure of the metric space $(\bar F, g_{o}' ).$ 
\end{proposition}

\noindent
{\bf Proof:}
 It is clear that $I_{\infty} \subset  \partial F$ and so only the 
reverse inclusion in (4.38) needs to be proved. Via Proposition 2.1, if 
$y_{\varepsilon}\in U^{\upsilon_{o}}$ converges to $y\in\bar F$ (in the 
Hausdorff topology) then $\nu_{\delta_{a}}(y_{\varepsilon}) \rightarrow 
 -\infty$ if and only if $y\in\partial F.$ This does not immediately 
imply (4.38), (since $\nu_{\delta_{a}}$ does not apriori converge to 
$\nu $ at $\partial F).$ However, if for any $y\in\partial F$ there 
exist $y_{\varepsilon}\rightarrow y$ as above such that, for some 
$m_{1}>0$,
\begin{equation} \label{e4.41}
\widetilde m_{\delta_{a}}(\tfrac{1}{2}\rho (y_{\varepsilon})) \geq  
m_{1}\cdot \rho (y_{\varepsilon}), 
\end{equation}
then (4.38) does follow. Here $\widetilde m_{\delta_{a}}$ is the mass 
of $\nu_{\delta_{a}}$, $\delta_{a} = \delta_{a}(x_{\varepsilon}),$ on 
the level through $y_{\varepsilon}$, as in (3.11). Moreover (4.41) 
implies (4.39) provided the constant $m_{1}$ is independent of 
$y\in\partial F.$ Now if 
\begin{equation} \label{e4.42}
\delta_{a}(y_{\varepsilon}) \geq  c_{o}\cdot 
\delta_{a}(x_{\varepsilon}), 
\end{equation}
for some $c_{o} > $ 0, then $\widetilde m_{\delta_{a}(y_{\varepsilon})} 
\leq  c_{o}^{-1} \widetilde m_{\delta_{a}(x_{\varepsilon})}.$ But by 
(3.13), $\widetilde m_{\delta_{a}(y_{\varepsilon})}(\frac{1}{2}\rho 
(y_{\varepsilon})) \geq  \rho (y_{\varepsilon}).$ Hence, (4.41) follows 
from (4.42).

 The proof of (4.42) essentially follows from the proof of Lemma 4.5. 
Thus, consider paths $\gamma  = \gamma_{y,\varepsilon}$ in 
$(\Omega_{\varepsilon}, g_{\varepsilon})$, parametrized by arclength 
$s$, starting at a point $q_{\varepsilon}\in L^{\upsilon_{o}}$ with 
$q_{\varepsilon} \rightarrow y$, ending at $x_{\varepsilon},$ and 
satisfying the cone condition (4.10). If (4.42) were false along any 
such curve $\gamma ,$ then there is a constant $C_{o} <  \infty $ such 
that, for all $s$ sufficiently small and hence for all $s \leq S_{o}$, 
$S_{o} < \infty$,
$$\delta_{a}(\gamma (s)) \leq  C_{o}\cdot \delta_{a}(x_{\varepsilon}). 
$$
Let $L_{\varepsilon}$ be the $u$-level set through $x_{\varepsilon}$. 
As in (4.12)-(4.14), it then follows that
\begin{equation} \label{e4.43}
|u(L^{\upsilon_{o}}) -  u(L_{\varepsilon})| \leq  
C_{1}\delta_{a}(x_{\varepsilon}) \cdot |\ln \bar \varepsilon|, 
\end{equation}
where $C_{1}$ depends only on $C_{o}$ and $S_{o}$. However, in general 
one has the estimate
\begin{equation} \label{e4.44}
\delta_{a}(x_{\varepsilon}) \leq  c_{2}\cdot |u- 1|(x_{\varepsilon}), 
\end{equation}
for generic $x_{\varepsilon}$ in the sense of (2.78), for some uniform 
$c_{2} <  \infty .$ To see this, (4.44) is equivalent to 
$|\nu_{\delta_{a}(x_{\varepsilon})}|(x_{\varepsilon}) \geq  
c_{2}^{-1}.$ If this estimate failed, i.e. 
$|\nu_{\delta_{a}(x_{\varepsilon})}|(x_{\varepsilon}) << $ 1, then the 
limit function $\nu_{a}$ would be 0 at generic base points, 
contradicting the fact that $\nu_{a}$ is non-trivial, c.f. Theorem 
2.11. 

 Since, by hypothesis, $|u- 1|(x_{\varepsilon}) \leq  2/|\ln \bar 
\varepsilon|^{2},$ (4.43) and (4.44) give $|u(L^{\upsilon_{o}})- u(L)| 
\rightarrow $ 0 as $\varepsilon  \rightarrow $ 0, a contradiction. This 
establishes (4.42) and hence (4.38)-(4.39), provided the point $y \in 
\partial F$ is accessible, in the sense that $y = \lim 
y_{\varepsilon},$ with $y_{\varepsilon}$ joined to $x_{\varepsilon}$ by 
a curve satisfying the cone condition (4.10). Let $\partial_{a}F 
\subset \partial F$ denote the subset of accessible points.

\smallskip

 By the definition of the measure $d\mu_{\infty},$ the estimate (4.40) 
follows immediately from Proposition 3.3 for balls $B_{y}(1)$ with 
center $y \in \partial F$ and $dist_{g_{o}'}(y, x) \leq 2$, (for 
example). Together with the results above, this implies that 
$\partial_{a}F$ has Hausdorff dimension at most 1, with locally finite 
$\mathcal{H}_{1}$, within $\bar B_{x}(2)$, c.f. [13, Ch.4] for 
discussion of Hausdorff measures. In turn, this now implies that all 
points of $\partial F$ within $\bar B_{x}(2)$ are accessible. Further, 
the arguments above may be applied at any other base point $x'\in F,$ 
with $x_{\varepsilon}' \rightarrow x'$ and with 
$\delta_{a}(x_{\varepsilon})$ replaced by 
$\delta_{a}(x_{\varepsilon}')$. The limit of 
$\delta_{a}(x_{\varepsilon}')/\delta_{a}(x_{\varepsilon})$ as 
$\varepsilon \rightarrow 0$ varies continuously as $x'$ varies over 
$F$, c.f. (the end of) Appendix B, and hence it follows that $\partial 
F$ everywhere has Hausdorff dimension at most 1, with locally finite 
$\mathcal{H}_{1}$.

\smallskip

 To see that (4.39) holds globally on $F$, i.e. $m_{o}$ is independent 
of center point of $B(r)$, the local estimates above prove that 
(4.39)-(4.40) hold in any ball $B_{x'}(2),$ when $\delta_{a}' $ is used 
in place of $\delta_{a}$ to define the local mass $m_{\infty},$ and 
with $m_{o}, M_{o}$ then independent of $x'$. This means that $\partial 
F$ is not too dense in $F$, and it follows that there exists $r_{o} = 
r_{o}(m_{o}, M_{o}) > 0$ such that any point $x'\in F$ with $dist(x', 
\partial F) \geq  r_{o}$ may be joined to $x$ via a path $\gamma$ with 
$dist(\gamma (s), \partial F) \geq  r_{o}.$ Hence, since any $y \in F$ 
satisfies $y \in B_{x'}(2)$ for some $x'$ as above, there also exist 
paths $\gamma_{\varepsilon}$ from $y_{\varepsilon}$ to 
$x_{\varepsilon}$ satisfying the cone condition (4.10); the constant 
$C$ in (4.10) of course depends on $y$, but $c_{o}$ in (4.42) is 
independent of $y$, in the limit $\varepsilon  \rightarrow 0$.

 Finally, to prove that (4.40) holds globally, i.e $M_{o}$ is 
independent of the base point, suppose $x' $ is a base point in $F$ 
such that $m_{\infty}(\partial F\cap B_{x'}(1)) >> 1$, where of course 
$m_{\infty}$ is the mass w.r.t. the potential $\nu_{\delta_{a}}$ based 
at $x$. This implies that $\delta_{a}(x_{\varepsilon}' 
)/\delta_{a}(x_{\varepsilon}) >> $ 1 in the limit on $F$. Interchanging 
the roles of $x$ and $x' ,$ it follows that $m_{\infty}' (\partial 
F\cap B_{x}(1)) << $ 1, where $m_{\infty}' $ is the mass w.r.t. the 
potential $\nu_{\delta_{a}'}$ based at $x'$. This contradicts the fact 
that (4.39) holds globally, (at the base point $x' $ in place of $x$), 
which gives (4.40).
{\endproof}

 We point out that the assumption $|u-1|(x_{\varepsilon}) \leq 2 /(\ln 
\bar \varepsilon)^{2}$ was used only in (4.43)-(4.44). Thus, 
Proposition 4.8 remains valid under the assumption 
$\delta_{a}(x_{\varepsilon}) \leq 2 /(\ln \bar \varepsilon)^{2}$.

\smallskip

  Recall that a subset of ${\mathbb R}^{3}$ is {\it  polar}  if it is 
contained in a countable collection of compact sets, each of which has 
capacity 0. These are precisely the sets where non-positive subharmonic 
functions in domains in ${\mathbb R}^{3}$ can become infinite, (i.e. 
$-\infty ),$ c.f. [12]. A polar set has Hausdorff dimension at most 1. 
For the discussion below, suppose that the maximal limit $F$ is flat; 
the discussion holds equally well for hyperbolic limits, with obvious 
changes. (The following remarks are not actually needed elsewhere in 
the paper).

\smallskip

 The boundary $\partial F$ may be divided into two parts, a {\it  
regular}  part $\partial_{r}F$ and a {\it  singular}  part, where 
$p\in\partial_{r}F$ if and only if $p\in\partial F$ and there is a 
neighborhood $U$ of $p$ in $\bar F$ such that $U$ embeds in ${\mathbb 
R}^{3}.$ Note that such an extension of $F$ past $\partial_{r}F$ is 
unique. By Proposition 4.8, locally the regular boundary is contained 
in a closed set of finite $\mathcal{H}_{1}$ measure in ${\mathbb 
R}^{3}$ and is of course a polar set. In particular, if $\partial_{s}F 
= \emptyset  $ and $F$ is flat, then $\bar F$ is a complete, smooth 
flat manifold and hence
\begin{equation} \label{e4.45}
\bar F = {\mathbb R}^{3}/\Gamma , 
\end{equation}
where $\Gamma  \subset Isom({\mathbb R}^{3})$ is a discrete group, 
maybe trivial, acting freely on ${\mathbb R}^{3}.$ 

 The local structure of $\partial_{s}F$ is determined by the holonomy 
of the developing map $\mathcal{D} : \widetilde F \rightarrow  {\mathbb 
R}^{3},$ where $\widetilde F$ is the universal cover of $F$. For 
simplicity, assume here that $F$ is maximally extended as a smooth flat 
manifold so that $\partial F = \partial_{s}F.$ For $q\in\partial F$ and 
$r$ sufficiently small, let $T_{r}$ denote the component of the 
$r$-tubular neighborhood of $\partial F$ containing $q$, and let 
$\partial T_{r}\subset F$ denote its boundary. Also, let $S_{q}(r')$ 
denote the geodesic sphere about $q$, $r' \geq r$. For generic $r$, 
$r'$, $S_{q}(r')\cap\partial T_{r}$ is a collection of embedded curves. 
Any such curve $\gamma$ is either a loop, (linking $\partial F),$ or a 
complete embedding of ${\mathbb R} ,$ i.e. $\gamma $ has no endpoints. 
The former case occurs when $B_{q}(r')$ is compact, while the latter 
occurs when $B_{q}(r')$ is non-compact. Let $\widetilde \gamma$ be a 
lift of $\gamma $ to $\widetilde F,$ so that $\widetilde \gamma$ is an 
embedding of ${\mathbb R} $ in $\widetilde F.$ Then $\mathcal{D} 
(\widetilde \gamma)$ is a closed loop in ${\mathbb R}^{3},$ and 
$\mathcal{D}|_{\widetilde \gamma}$ acts as a covering projection 
${\mathbb R}  \rightarrow  S^{1}.$ The loop $\mathcal{D} (\widetilde 
\gamma)$ determines an isometric ${\mathbb Z}$-action on ${\mathbb 
R}^{3}.$ 

 Any isometric ${\mathbb Z}$-action on ${\mathbb R}^{3}$ is either 
rotational, i.e. rotation about an axis $A$, or a translation or a 
twist, i.e. a rotation about and translation along an axis $A$. In the 
latter two cases, the corresponding quotient spaces are the product 
$S^{1}\times {\mathbb R}^{2}$ and twisted product 
$S^{1}\times_{\alpha}{\mathbb R}^{2}$ respectively, where $\alpha$ is 
the rotation angle. These quotients are everywhere smooth and hence do 
not contribute $\partial_{s} F$. The only source of singularity to 
$\partial F$ is that given by rotational holonomies. It follows that 
$\partial_{s}F$ is given locally by a collection of lines, 
corresponding to rotation axes. A neighborhood of each line in the 
closure $\bar F$ is a neighborhood of (0,0) in ${\mathbb R}\times C,$ 
where $C$ is a flat cone of cone angle $\alpha ;$ the angle $\alpha $ 
may assume any value in (0, $\infty ],$ the value $\alpha  = \infty $ 
corresponding to the universal cover of $({\mathbb R}^{2}\setminus 
\{0\}).$

\smallskip

 We close this section with the following useful result.
\begin{lemma} \label{l 4.9.}
  For $F$ as in Proposition 4.8, let $t(x) = dist_{g_{o}'}(x, \partial 
F).$ Then $t$ is unbounded on $F$.
\end{lemma}

\noindent
{\bf Proof:}
 The proof is by contradiction, so suppose that $T_{o} = \sup_{F}t < 
\infty$. By Lemma 3.2, the limit potential $\nu_{a}$ is then bounded 
above on $F$, $\nu_{a} \leq C$, for some $C < \infty$. The potential 
$\nu_{a}$ satisfies the elliptic equation (2.65). Suppose first that 
$\Delta \nu_{a} = 0$, and consider a maximizing sequence $\{y_{j}\}$ 
for $\nu_{a}$, i.e. $\nu_{a}(y_{j}) \rightarrow \sup_{F}\nu_{a}$. The 
maximum principle implies that $\nu_{a}$ approches a constant function 
near $y_{j}$, for $j$ sufficiently large, which in turn implies that 
the mass of $\nu_{a}$ in $B_{y_{j}}(1)$ approaches 0. However, this 
contradicts the uniform lower bound (4.39) for the mass, again since 
$t(y_{j}) \leq T_{o} < \infty$.

  If $\nu_{a}$ is not harmonic, then Theorem 2.11 implies that 
$\nu_{a}$ is superharmonic. Let $h$ be the harmonic function on $F$ 
determined by the Riesz measure $d\mu_{\infty}$ on $\partial F$. Since 
$\nu_{a}$ is superharmonic, $h \leq \nu_{a} \leq C$. Then exactly the 
same argument as above with $h$ in place of $\nu_{a}$, gives the same 
contradiction.
{\endproof}

\section{Mass Gaps and Essential Spheres.}
\setcounter{equation}{0}

 This is the most important part of the paper, and culminates in 
Theorems 5.9 and 5.12, which prove the existence of essential 2-spheres 
in $M$ under two mild hypotheses. The first of these is the 
degeneration hypothesis (4.1) and the second, (which has two parts), is 
a non-collapse hypothesis.

 In \S 5.1, we describe the non-collapse hypothesis (first part), and 
derive some elementary consequences for the structure of the level set 
$L_{o}$ of admissible base points. Next, the work in \S 5.2 proves that 
if $x_{\varepsilon}$ is any admissible base point, then there exists a 
well-separated base point $y_{\varepsilon}$ on the same level $L_{o},$ 
(naturally constructed from the initial $x_{\varepsilon}),$ c.f. 
Theorem 5.4. This crucial result is obtained by examining the mass 
distribution on the level $L_{o}$ near $x_{\varepsilon},$ using the 
global mass estimate from Theorem 3.5 in combination with the local 
mass bound from Proposition 4.7.

 In \S 5.3, it is then shown that well-separated base points are 
naturally associated to 2-spheres in $M$, under the $2^{\rm nd}$ 
non-collapse hypothesis, c.f. Theorem 5.9. This leads to Theorem 5.12, 
proving that such 2-spheres are essential in $M$. 

 It is assumed throughout \S 5 that $M$ is $\sigma$-tame. Thus, the 
sequence $(\Omega_{\varepsilon}, g_{\varepsilon}), \varepsilon  = 
\varepsilon_{i},$ is chosen to satisfy (4.3).

\medskip

{\bf 5.1.}
 Theorem 4.4 gives the existence of special base points on 
$(\Omega_{\varepsilon}, g_{\varepsilon}),$ namely admissible or 
preferred. In this subsection, we derive some elementary consequences 
for the geometry near such base points, in the presence of the 
following assumption.

{\bf Non-Collapse Assumption I.}
 {\it There exists an admissible or preferred base point 
$x_{\varepsilon}$ and a constant $\mu_{o} > 0$, independent of 
$\varepsilon ,$ such that for all $p_{\varepsilon}\in L_{o}\cap 
B_{x_{\varepsilon}}(1),$ or $p_{\varepsilon} = x_{\varepsilon}$ 
respectively,
\begin{equation} \label{e5.1}
\frac{vol B_{p_{\varepsilon}}(\rho (p_{\varepsilon}))}{\rho 
(p_{\varepsilon})^{3}} \geq  \mu_{o}. 
\end{equation}
where $L_{o}$ is the level set through $x_{\varepsilon},$ as in (4.34)}.

\smallskip

 The bound (5.1) means that $(\Omega_{\varepsilon}, g_{\varepsilon})$ 
is not volume collapsing within $\rho$-balls centered at 
$p_{\varepsilon}$, i.e. the rescaled metrics $(\Omega_{\varepsilon}, 
g_{\varepsilon}' , p_{\varepsilon})$, $g_{\varepsilon}'  = \rho 
(p_{\varepsilon})^{-2}\cdot  g_{\varepsilon}$, have a uniform lower 
bound on the volume radius at $p_{\varepsilon},$ and so one has 
convergence to the limit $(F, g_{o}', p)$, (in a subsequence). It is 
clear that (5.1) implies that there exists $a_{o} = a_{o}(\mu_{o}) > 0$ 
such that, for $H_{p_{\varepsilon}}(r) = B_{p_{\varepsilon}}(r)\cap 
L_{o},$
\begin{equation} \label{e5.2}
\frac{area H_{p_{\varepsilon}}(\rho (p_{\varepsilon}))}{\rho 
(p_{\varepsilon})^{2}} \geq  a_{o}. 
\end{equation}
As always, (5.2) is understood to hold for generic points 
$p_{\varepsilon},$ as in (2.79), (although it follows from (5.5) below 
that it holds in fact for all $p_{\varepsilon}$). The ball 
$B_{x_{\varepsilon}}(1)$ may be replaced by 
$B_{x_{\varepsilon}}(\eta)$, for any fixed $\eta > 0$, throughout \S 5. 
In \S 6, we will deal separately with the possibility of collapse, 
where (5.1) does not hold, (for {\it all} admissible or preferred base 
points). The bound (5.1) is part of the assumptions of {\it all}  of 
the results to follow in \S 5, but will not be mentioned explicitly 
each time.

\smallskip

 Referring to Proposition 4.7, the following Corollary is now obvious.

\begin{corollary} \label{c 5.1.}
  Let $x_{\varepsilon}$ be an admissible or preferred base point 
satisfying (5.1). There is a constant $\kappa_{1} > $ 0, depending only 
on $\kappa_{o}$ and $\mu_{o},$ and in particular independent of 
$\varepsilon ,$ such that for any $q_{\varepsilon}\in L_{o}\cap 
B_{x_{\varepsilon}}(1),$ respectively $q_{\varepsilon} = 
x_{\varepsilon},$
\begin{equation} \label{e5.3}
m_{o}(\rho (q_{\varepsilon})) \geq  \kappa_{1}\cdot \rho 
(q_{\varepsilon}). 
\end{equation}
\end{corollary}

\noindent
{\bf Proof:}
 This is an immediate consequence of (4.35) and (5.1), c.f. (4.36).
{\endproof}

\begin{remark} \label{r 5.2.}
  Observe that under the assumption (5.1), one has
\begin{equation} \label{e5.4}
\rho (q_{\varepsilon}) \leq  \tfrac{1}{2}\kappa_{1}^{-1}\delta_{o} = 
\tfrac{1}{2}\kappa_{1}^{-1}\bar \varepsilon^{2\mu} \rightarrow  0, \ \ 
{\rm as} \ \ \varepsilon  \rightarrow  0, 
\end{equation}
for all $q_{\varepsilon}\in L_{o}\cap B_{x_{\varepsilon}}(1)$ or 
$q_{\varepsilon} = x_{\varepsilon}$ respectively. To see this, by 
(5.3), the mass of $\nu_{\delta_{o}}$ in $B_{q_{\varepsilon}}(\rho 
(q_{\varepsilon}))\cap L_{o}$ is at least $\kappa_{1}\rho 
(q_{\varepsilon}).$ On the other hand, Theorem 3.5 implies that the 
total mass $m_{o}$ of $L_{o}$ is bounded by $\frac{1}{2}\delta_{o}$, 
which gives (5.4). In particular, by (4.5), it follows that
\begin{equation} \label{e5.5}
\rho^{2}(q_{\varepsilon}) <<  \delta_{a}(q_{\varepsilon}), 
\end{equation}
corresponding to Case (i) of (2.12). In the opposite direction, 
Proposition 1.4 gives a uniform lower bound on $\rho ,$ i.e. $\rho^{2} 
\geq  \bar \varepsilon$ everywhere and so in particular on $L_{o}.$
\end{remark}

\smallskip

 Theorem 3.5 and (5.3) also combine to give the following result.

\begin{corollary} \label{c 5.3.}
  Let $x_{\varepsilon}$ be an admissible or preferred base point and 
let $C_{o}$ be any component of $L_{o}\cap B_{x_{\varepsilon}}(1)$, 
respectively $L_{o}\cap B_{x_{\varepsilon}}(\rho (x_{\varepsilon})).$ 
Then, $\forall p_{\varepsilon} \in C_{o}$, there exists $r_{o} > 0$ s.t.
\begin{equation} \label{e5.6}
diam_{g_{\varepsilon}}C_{o} \geq  r_{o}\cdot \rho (p_{\varepsilon}), 
\end{equation}
where $diam$ is the extrinsic diameter in $(\Omega_{\varepsilon}, 
g_{\varepsilon})$, and $r_{o}$ is independent of $C_{o}$ and 
$\varepsilon$. Hence
\begin{equation} \label{e5.7}
diam_{g_{\varepsilon}}C_{o} \geq  r_{1}\bar \varepsilon^{1/2}. 
\end{equation}
Further,
\begin{equation} \label{e5.8}
diam_{g_{\varepsilon}}C_{o} \leq  5\kappa_{1}^{-1}\bar 
\varepsilon^{2\mu} \rightarrow  0, \ \ {\rm as} \ \ \varepsilon  
\rightarrow  0. 
\end{equation}
\end{corollary}

\noindent
{\bf Proof:}
 The proof is the same in the admissible and well-separated cases, so 
assume $x_{\varepsilon}$ is admissible. The estimate (5.6) follows 
exactly as in (2.79) from the uniform non-triviality of the limit 
potential $\nu_{a}$. Note that (5.5) implies that the limit potential 
is harmonic, and so has no local maxima, (or minima). Thus (2.79) and 
(5.6) hold for all $p_{\varepsilon} \in C_{o}$, and not just generic 
$p_{\varepsilon}$ in the sense of (2.78). The estimate (5.7) is an 
immediate consequence of (5.6) and Proposition 1.4.

  To prove (5.8), consider the cover $\mathcal{U}$ of $C_{o}$ by 
$\rho$-balls $B_{z}(\rho (z))$, $z \in C_{o}$. By a standard covering 
argument, c.f. [10, B33.3] or [13, \S 2.1], one may select a finite 
subcollection of disjoint balls $B_{z_{j}}(\rho(z_{j}))$ such that the 
balls $B_{z_{j}}(5\rho(z_{j}))$ cover $C_{o}$. Then
\begin{equation} \label{e5.9}
diam_{g_{\varepsilon}}(C_{o}) \leq  5\sum\rho (z_{j}) \leq  
5\kappa_{1}^{-1}\sum m_{o}(B_{z_{j}}(\rho (z_{j}))) \leq 
5\kappa_{1}^{-1}\bar \varepsilon^{2\mu},
\end{equation}
where the second estimate uses (5.3) and the last uses (3.22) and (4.6).
{\endproof}

 Of course, combining (5.7) and (5.8) gives an upper bound, {\it  
depending on}  $\varepsilon ,$ on the number of components of $L_{o}$ 
in the ball $B_{x_{\varepsilon}}(1) \subset  (\Omega_{\varepsilon}, 
g_{\varepsilon}).$

\smallskip

 It is worthwhile to keep in mind that Corollary 5.3 does not assert 
that $diam_{g_{\varepsilon}}C_{o} \sim  \rho (q_{\varepsilon}),$ for 
some $q_{\varepsilon}\in C_{o}.$ In fact, it may well be the case that 
for all $q_{\varepsilon}\in C_{o}$, $diam_{g_{\varepsilon}}C_{o} >>  
\rho (q_{\varepsilon})$. Nevertheless, Corollary 5.3 does show that all 
the components of $L_{o}$ in $B_{x_{\varepsilon}}(1)$ are being crushed 
(metrically) to points as $\varepsilon  \rightarrow $ 0 and not to more 
complicated spaces. However, the {\it  distribution}  of these points 
might be very complicated; for instance they could apriori become 
arbitrarily dense in $B_{x_{\varepsilon}}(1)$ as $\varepsilon  
\rightarrow $ 0, compare with Remark 1.3.

\medskip

{\bf 5.2.}
 We now come to perhaps the central technical result of the paper, 
proving the existence of gaps of a {\it  definite relative size}  
(independent of $\varepsilon )$ in the level set $L_{o}$ near an 
admissible base point. This result implies that if $x_{\varepsilon}$ is 
admissible, then there is a well-separated base point 
$y_{\varepsilon}\in L_{o}$ near to $x_{\varepsilon},$ c.f (4.7); in 
fact $x_{\varepsilon}$ is well-separated w.r.t. all of 
$\Omega_{\varepsilon},$ and not only w.r.t. $U^{\upsilon_{o}}.$ This 
allows one to pass from admissible base points to well-separated base 
points.

\begin{theorem} \label{t 5.4.}
  For any admissible base point $x_{\varepsilon}\in L_{o} = 
L_{o}(\varepsilon )$ there exist points $y_{\varepsilon}\in 
B_{x_{\varepsilon}}(1)\cap L_{o}$ and scales $\lambda_{\varepsilon} = 
\lambda_{\varepsilon}(y_{\varepsilon})\in $(0, 1), 
$\lambda_{\varepsilon} \geq  r_{o}\cdot \rho (y_{\varepsilon})$, for a 
fixed $r_{o} > 0$, such that 
\begin{equation} \label{e5.10}
L_{o} \cap  A_{y_{\varepsilon}}(\lambda_{\varepsilon}, 
(1+d_{o})\lambda_{\varepsilon}) = \emptyset  , 
\end{equation}
where $d_{o} = d_{o}(\mu ) > $ 0 is independent of $\varepsilon $ and 
$y_{\varepsilon}.$ 
\end{theorem}

\noindent
{\bf Proof:}
 Given an admissible base point $x_{\varepsilon}\in L_{o},$ consider 
the unit ball $B = B_{x_{\varepsilon}}(1)$ w.r.t. $g_{\varepsilon}.$ To 
each component $C_{j}$ of $L_{o}$ contained in $B$ associate a 
subinterval $I_{j} = [a_{j}, b_{j}]$ of $I = [0,1]$ by setting 
$$a_{j} = min\{r: C_{j}\cap B_{x_{\varepsilon}}(r) \neq  \emptyset \}, 
\ \  b_{j} = max\{r: C_{j}\cap B_{x_{\varepsilon}}(r) \neq  \emptyset 
\}. $$
By (5.6) and (5.7), the length $l(I_{j})$ satisfies 
\begin{equation} \label{e5.11}
l(I_{j}) = b_{j} -  a_{j} \geq  r_{o}\rho (q_{j}) \geq  r_{1}\bar 
\varepsilon^{1/2}, \ \ {\rm for \ some} \ \  q_{j}\in C_{j}. 
\end{equation}
Further, to each $I_{j}$ associate a mass $m$ by setting $m(I_{j}) 
\equiv m_{o}(C_{j})$. Then the same estimates as in (5.9) give
\begin{equation} \label{e5.12}
m(I_{j}) \geq  \kappa_{2}l(I_{j}), 
\end{equation}
where $\kappa_{2} = \kappa_{2}(\kappa_{1})$. The intervals $I_{j}$ may 
of course intersect, but this will only add further to the mass. Since 
the total mass of $L_{o}(1) \equiv  L_{o}\cap B_{x_{\varepsilon}}(1)$ 
is small, i.e. 
\begin{equation} \label{e5.13}
m(L_{o}(1)) \leq  \bar \varepsilon^{2\mu}, 
\end{equation}
by (3.22), the total Lebesque measure of $I_{o} \equiv  \cup I_{j}$ is 
very small in $I$. Hence there exist gaps between the intervals 
$I_{j}.$ However, there may be very many, very small gaps. To estimate 
the gap size, we essentially consider $I_{o}$ as obtained from $I$ by 
the usual Cantor set construction of deleting open intervals from $I$.

 Thus, consider $I_{o}$ as obtained from $I$ by inductively dividing 
intervals, starting at $I$, into two parts obtained by removing an open 
interval of {\it  relative}  length $\leq  d_{o};$ the parameter 
$d_{o}$ will be estimated below. Initially then, remove an open 
interval of length $d_{1}$ from $I$, so there are left two closed 
intervals, of total length $(1- d_{1}).$ Remove from either or both of 
these again an interval of total relative length $d_{2},$ so that total 
length of removed intervals is then $d_{1} + d_{2}(1- d_{1}) = 1 - 
(1-d_{1})(1-d_{2})$. Repeat this process inductively $N$ times, until 
the collection of intervals $I_{o}$ is reached. 

 A straightforward computation then shows that the total length of the 
removed intervals equals $1 -  \prod_{1}^{N}(1 - d_{i})$. It follows 
that the total length of the remaining closed intervals, and hence the 
total mass of $I_{o}$ by (5.12), satisfies
\begin{equation} \label{e5.14}
m(I_{o}) \geq  \kappa_{2}\prod_{1}^{N}(1 - d_{i}). 
\end{equation}
Let $\bar d = \sup_{i}d_{i}$, so that (5.14) gives
$$m(I_{o}) \geq  \kappa_{2}(1 -  \bar d)^{N}. $$
Then from (5.13) one obtains the bound
\begin{equation} \label{e5.15}
(1 -  \bar d)^{N} \leq  \kappa_{2}^{-1}\bar \varepsilon^{2\mu}. 
\end{equation}
On the other hand, since each repetition of the process above produces 
a new closed interval of length $<  \frac{1}{2}$ that of its 
predecessor, some of the final closed intervals $I_{j}$ in $I_{o}$ must 
have length $\leq  2^{-N}.$ But by (5.11), the length of each interval 
$I_{j}$ is at least $\bar \varepsilon,$ (in fact $\bar 
\varepsilon^{1/2}).$ Hence, one also has
\begin{equation} \label{e5.16}
\bar \varepsilon \leq  2^{-N}. 
\end{equation}
An elementary calculation combining (5.15) and (5.16) then gives
$$\bar d >  1 -  \kappa_{2}^{-1/N}2^{- 2\mu}.$$ 

 This implies that there are gaps in $I_{o}$ at {\it  some}  scale, of 
relative size at least 
\begin{equation} \label{e5.17}
d_{o} = 1 -  \kappa_{2}^{-1}2^{- 2\mu}. 
\end{equation}
By construction, this means that there exist points $y_{\varepsilon}\in 
L_{o}\cap B_{x_{\varepsilon}}(1)$ and numbers $\lambda_{\varepsilon}\in 
(0,1]$ such that (5.10) holds. The fact that $\lambda_{\varepsilon} 
\geq  r_{o}\cdot \rho (y_{\varepsilon})$ follows from (5.6).

{\endproof}

\begin{remark} \label{r 5.5.}
  It is clear from the proof of Theorem 5.4 that there could be many 
possible choices of base points $y_{\varepsilon}$ and factors 
$\lambda_{\varepsilon}$ satisfying (5.10). For the purposes of this 
paper, this lack of uniqueness, already occuring in \S 4, plays no 
fundamental role.

 However, it is not asserted that there exist gaps forming as in (5.10) 
at some definite or determinable scale, (although see Lemma 5.6). For 
instance, it is not claimed that one can choose $\lambda_{\varepsilon}$ 
to be bounded away from 0 as $\varepsilon  \rightarrow $ 0, i.e. that 
there are gaps forming in $L_{o}$ of a definite size w.r.t. 
$g_{\varepsilon},$ compare with Remark 1.3. At the other, small-scale 
extreme, it is not known that gaps form on the scale of $\rho 
(y_{\varepsilon}),$ i.e. whether $\lambda_{\varepsilon} \sim  \rho 
(y_{\varepsilon}),$ as is the case for preferred base points.

 It is the existence of such gaps provided by Theorem 5.4 which is 
important, and not so much the scale at which they are visible. Without 
the $\sigma$-tame assumption in Theorem 5.4 the size of all relative 
gaps could approach 0 as $\varepsilon  \rightarrow $ 0, as in the 
construction of Cantor sets of measure 0 in $[0,1]$ with Hausdorff 
dimension 1, c.f. [13, \S 4.10]. This is in fact the only place in the 
paper where the $\sigma$-tame assumption is used in an essential way; 
(the arguments in \S 4 and in \S 5 preceding Theorem 5.4 can be 
modified to hold without the $\sigma$-tame assumption).
\end{remark}

\smallskip

 The following (more technical) observation based on the proof of 
Theorem 5.4 will be used later in the final choice of base points.

\begin{lemma} \label{l 5.6.}
  Let $x_{\varepsilon}$ be any admissible base point. Then there exist 
well-separated base points $y_{\varepsilon}\in 
B_{x_{\varepsilon}}(1)\cap L_{o}$ satisfying (5.10), together with
\begin{equation} \label{e5.18}
\lambda_{\varepsilon} \leq  \bar \varepsilon^{3\mu /2}, 
\end{equation}
and such that
\begin{equation} \label{e5.19}
m_{o}(B_{y_{\varepsilon}}(\lambda_{\varepsilon})) \geq  \bar 
\varepsilon^{\mu /2}\lambda_{\varepsilon}. 
\end{equation}
\end{lemma}

\noindent
{\bf Proof:}
 To prove (5.18), observe that the proof of Theorem 5.4 began with the 
ball $B_{x_{\varepsilon}}(1).$ However, it can be applied equally well 
to balls $B_{x_{\varepsilon}}(s_{\varepsilon}),$ with $s_{\varepsilon} 
\rightarrow $ 0 sufficiently slowly as $\varepsilon  \rightarrow $ 0. 
For example, let $s_{\varepsilon} = \bar \varepsilon^{3\mu /2}$ and 
rescale the balls $B_{x_{\varepsilon}}(s_{\varepsilon})$ to size 1, 
i.e. set $\bar g_{\varepsilon} = s_{\varepsilon}^{-2}\cdot  
g_{\varepsilon}.$ Since the mass scales as a distance, $\bar 
m_{o}((B_{x_{\varepsilon}}(1)) \leq  \bar \varepsilon^{-3\mu /2}\bar 
\varepsilon^{2\mu} \leq  \bar \varepsilon^{\mu /2},$ where $\bar m_{o}$ 
is the mass $m_{o}$ w.r.t the metric $\bar g_{\varepsilon}.$ Now 
starting with (5.13), with $\bar \varepsilon^{\mu / 2}$ in place of 
$\bar \varepsilon^{2\mu}$, the same arguments as before in Theorem 5.4 
apply to give (5.10), with $\lambda_{\varepsilon}$ satisfying (5.18). 
(The estimate for $1 - d_{o}$ is smaller than that in (5.17) by a 
factor of 2). 

 The proof of (5.19) is similar. Thus, let 
$A_{y_{\varepsilon}}(\lambda_{\varepsilon})$ be an annulus satisfying 
(5.10) and (5.18) and consider the rescaled metrics $\bar 
g_{\varepsilon} = \lambda_{\varepsilon}^{-2}\cdot  g_{\varepsilon}.$ If 
there is any fixed constant $\kappa  > $ 0, say $\kappa  = \mu /2,$ s.t.
\begin{equation} \label{e5.20}
\bar m_{o}(\bar B_{y_{\varepsilon}}(1)) \leq  \bar 
\varepsilon^{\kappa}, 
\end{equation}
then exactly the same arguments as in the proof of (5.18) above may be 
repeated at this scale to produce a new smaller annulus 
$A_{y_{\varepsilon}^{1}}(\lambda_{\varepsilon}^{1}) \subset  
B_{y_{\varepsilon}}(\lambda_{\varepsilon}), \lambda_{\varepsilon}^{1} 
<<  \lambda_{\varepsilon},$ satisfying the gap condition (5.10), with 
$\lambda_{\varepsilon}^{1}$ in place of $\lambda_{\varepsilon};$ as 
above, here $d_{o} = d_{o}(\kappa ) = d_{o}(\mu ).$ In this new scale 
$g_{\varepsilon}^{1} = (\lambda_{\varepsilon}^{1})^{2}g_{\varepsilon}$, 
the mass is of course much larger. Hence, by iterating this procedure a 
finite number of times, (depending on $\varepsilon ),$ one obtains base 
points $y_{\varepsilon}$ and scale-factors $\lambda_{\varepsilon}$ 
satisfying (5.10) and (5.18)-(5.19).

{\endproof}

 Theorem 5.4 and Lemma 5.6 prove the existence of well-separated base 
points, (of course under the assumptions (4.1), (5.1) and the 
$\sigma$-tame assumption), satisfying the properties (5.10) and 
(5.18)-(5.19). Preferred base points automatically satisfy 
(5.18)-(5.19); for such base points, since $\lambda  \sim  \rho ,$ 
(5.18) follows from (5.4) while (5.19) follows from (5.3). From now on, 
it is assumed that the well-separated base point $y_{\varepsilon}$ has 
been chosen to satisfy (5.18)-(5.19). One further observation is needed 
for the final determination of the base point. 

\smallskip

 The level $L_{o}$ disconnects the manifold $\Omega_{\varepsilon}$ into 
two parts,
\begin{equation} \label{e5.21}
\Omega_{\varepsilon} = U_{o} \cup  U^{o}, 
\end{equation}
where $U_{o} = \{x\in\Omega_{\varepsilon}: u(x) <  u(L_{o})\}$, $U^{o} 
= \{x\in\Omega_{\varepsilon}: u(x) >  u(L_{o})\}$; of course neither of 
these sets is necessarily connected. Hence each component 
$C_{y_{\varepsilon}}(\lambda_{\varepsilon})$ of any annulus 
$A_{y_{\varepsilon}}(\lambda_{\varepsilon},(1+d_{o})\lambda_{\varepsilon
})$ as in Theorem 5.4 is contained either in $U_{o}$ or $U^{o}.$

\smallskip

 The next result asserts the existence of such components contained in 
$U^{o}.$

\begin{lemma} \label{l 5.7.}
  Let $A_{y_{\varepsilon}}(\lambda_{\varepsilon}, 
(1+d_{o})\lambda_{\varepsilon}))$ be any annulus satisfying (5.10), 
with $\lambda_{\varepsilon}$ satisfying (5.18). Then there is a 
component $\mathcal{C}_{\varepsilon}  = 
\mathcal{C}_{y_{\varepsilon}}(\lambda_{\varepsilon})$ of 
$A_{y_{\varepsilon}}(\lambda_{\varepsilon}, 
(1+d_{o})\lambda_{\varepsilon})$ such that
\begin{equation} \label{e5.22}
\mathcal{C}_{\varepsilon}  \subset  U^{o}. 
\end{equation}
\end{lemma}

\noindent
{\bf Proof:}
 Let $s_{\varepsilon} = (1+\frac{d_{o}}{2})\lambda_{\varepsilon}.$ By 
(5.10), the ball $B(s_{\varepsilon}) = 
B_{y_{\varepsilon}}(s_{\varepsilon})$ satisfies $\partial 
B(s_{\varepsilon})\cap L_{o} = \emptyset  ,$ so that any component of 
$\partial B(s_{\varepsilon})$ is then contained either in $U_{o}$ or 
$U^{o}.$ It suffices then to prove that some component of $\partial 
B(s_{\varepsilon})$ is contained in $U^{o}.$

 To see this, suppose instead the full boundary $\partial 
B(s_{\varepsilon}) \subset  U_{o}.$ Then some part of the level set 
$L_{o}\cap B(s_{\varepsilon})$ bounds a connected domain $D = 
D_{\varepsilon} \subset  U^{o},$ with $D \subset\subset  
B(s_{\varepsilon}), \partial D \subset  L_{o}.$  Recall that by (5.4), 
$\rho  \leq  \bar \varepsilon^{2\mu}$ on $L_{o}\cap 
B_{y_{\varepsilon}}(\frac{1}{2}),$ (ignoring constant factors here). 
Hence, by the Lipschitz property of $\rho , \rho  \leq  \bar 
\varepsilon^{3\mu /2}$ everywhere in $B(s_{\varepsilon}),$ and so, at 
all base points in $D$, 
\begin{equation} \label{e5.23}
\rho^{2}/\delta_{o} \rightarrow  0. 
\end{equation}

 Choose a point $z_{\varepsilon}\in D$ which realizes the maximal value 
of $u$ on $D$, and work in the scale $g_{\varepsilon}' $ where $\rho' 
(z_{\varepsilon}) =$ 1. Since $\rho  \sim  dist (L^{\upsilon_{o}}, 
\cdot )$ and $L^{\upsilon_{o}}$ lies outside $D$, it follows that 
$dist_{g_{\varepsilon}'}(z_{\varepsilon}, L_{o}) \leq $ 2 for 
$\varepsilon $ sufficiently small. Consider the potential function 
$\nu_{\delta_{o}} = (u-1)/\delta_{o}.$ If the mass $m_{o}$ of 
$\nu_{\delta_{o}}$ w.r.t. $g_{\varepsilon}'$ is uniformly bounded over 
$B_{z_{\varepsilon}}'(1/2)$, then $\nu_{\delta_{o}}$ subconverges, 
modulo constants, to a limit $\nu_{o}.$ If the mass $m_{o}$ is very 
large, then renormalize it by dividing by 
$m_{o}(B_{z_{\varepsilon}}'(1/2))$. (This has effect of increasing 
$\delta ).$ In both cases, the estimate (5.23) shows that Case (i) of 
(2.12) holds and so Proposition 2.3 implies that the limit function 
$\nu_{o}$ satisfies $\Delta\nu_{o} =$ 0. Since, $\nu_{o}$ has a maximal 
value at $z = \lim z_{\varepsilon},$ it follows that $\nu_{o}$ is 
constant.

 However, the estimate (4.35) on $L_{o}$ passes to the limit, and shows 
that the limit potential $\nu_{o}$ has non-zero mass on $L_{o}.$ This 
gives a contradiction.
{\endproof}

 Lemma 5.7 applies equally well to preferred base points; such base 
points are only well-separated w.r.t. $U^{\upsilon_{o}},$ but Lemma 5.7 
implies that there is a full component $\mathcal{C}_{\varepsilon}$ of 
the annulus $A^{\upsilon_{o}}$ contained in $U^{o}$ and hence 
$\mathcal{C}_{\varepsilon}\cap L^{\upsilon_{o}} = \emptyset  .$

\medskip

{\bf 5.3.}
 We are now ready to to fix the choice of the base point and prove, 
under a $2^{\rm nd}$ natural non-collapse assumption, that such base 
points give rise to geometrically natural 2-spheres in $M$. 

\begin{definition} \label{d 5.8.}
  A {\sf distinguished} base point is a base point 
$q_{\varepsilon}\in\mathcal{C_{\varepsilon}}$ as in (5.22), with 
$q_{\varepsilon} \in 
A_{y_{\varepsilon}}((1+\frac{1}{4}d_{o})\lambda_{\varepsilon}, 
(1+\frac{3}{4}d_{o})\lambda_{\varepsilon})$, where $y_{\varepsilon}$, 
$\lambda_{\varepsilon}$ satisfy the conclusions of Lemma 5.6.
\end{definition}

 For the remainder of \S 5, we work with a choice of distinguished base 
point. Again there are apriori many possible choices here, and the 
results to follow in \S 5.3 apply to any such choice. To any 
distinguished base point $q_{\varepsilon}$ is associated the (center) 
base point $y_{\varepsilon}$ and component $\mathcal{C}_{\varepsilon} 
\subset U^{o}$ of $A_{y_{\varepsilon}}(\lambda_{\varepsilon}, 
(1+d_{o})\lambda_{\varepsilon}),$ with $q_{\varepsilon} \in 
\mathcal{C}_{\varepsilon}$.

 Consider the rescaled metric
\begin{equation} \label{e5.24}
\bar g_{\varepsilon} = \lambda_{\varepsilon}^{-2}\cdot  
g_{\varepsilon}, 
\end{equation}
so that $A_{\varepsilon}(\lambda_{\varepsilon}) = 
A_{y_{\varepsilon}}(\lambda_{\varepsilon}, 
(1+d_{o})\lambda_{\varepsilon})$ now becomes $\bar A_{\varepsilon}(1) = 
\bar A_{y_{\varepsilon}}(1, 1+d_{o})$ and 
$B_{\varepsilon}(\lambda_{\varepsilon}) = 
B_{y_{\varepsilon}}(\lambda_{\varepsilon})$ becomes $\bar 
B_{\varepsilon}(1) = \bar B_{y_{\varepsilon}}(1).$ For 
$\mathcal{C}_{\varepsilon}  \subset  \bar A_{\varepsilon}(1)$ as above 
satisfying (5.22), Proposition 1.1 implies
\begin{equation} \label{e5.25}
\bar \rho(p_{\varepsilon}) \geq  d_{o}/4, 
\end{equation}
for all $p_{\varepsilon}\in\bar A_{\varepsilon}(1+\frac{1}{4}d_{o}, 
1+\frac{3}{4}d_{o})\cap\mathcal{C}_{\varepsilon} .$ On the other hand, 
since $\lambda  \geq  r_{o}\cdot \rho (y_{\varepsilon}),$
$$\bar \rho(y_{\varepsilon}) \leq  r_{o}^{-1}, $$
and hence by the Lipschitz property of $\rho ,$
$$\bar \rho(q_{\varepsilon}) \leq  2r_{o}^{-1}.$$
Thus, the rescaling (5.24) is uniformly homothetic to the natural 
blow-up metric based at $q_{\varepsilon},$ i.e.
\begin{equation} \label{e5.26}
g_{\varepsilon}'  = \rho (q_{\varepsilon})^{-2}\cdot  g_{\varepsilon}. 
\end{equation}
 We now make the following $2^{\rm nd}$ non-collapse assumption.

\smallskip

{\bf Non-Collapse Assumption II.}
 {\it For some distinguished base point $q_{\varepsilon}$, there exists 
$\mu_{o} > 0$ and some sequence $K_{\varepsilon} \rightarrow \infty$, 
(slowly) as $\varepsilon \rightarrow 0$, such that, for all $ 1 \leq K 
\leq K_{\varepsilon}$}
\begin{equation} \label{e5.27}
\frac{B_{q_{\varepsilon}}(K\rho(q_{\varepsilon}))}{(K\rho(q_{\varepsilon
}))^{3}} \geq \mu_{o}.
\end{equation}

\smallskip

  The condition (5.27) means that the rescalings (5.26) based at 
$q_{\varepsilon}$ do not collapse on large scales in the metric (5.26), 
i.e. $\omega'(q_{\varepsilon}) >> 1$, as $\varepsilon \rightarrow 0$. 
Although similar, the assumption NCA II does not quite follow from NCA 
I, since $q_{\varepsilon} \notin L_{o}$ and the collapse is on scales 
large compared with $\rho(q_{\varepsilon})$. However, the assumption 
NCA II is much weaker than NCA I, since it applies to the behavior at a 
single base point $q_{\varepsilon}$, while NCA I applies to all 
$p_{\varepsilon} \in B_{x_{\varepsilon}}(1) \cap L_{o}$. 

  Under the NCA II, the pointed sequence $(\Omega_{\varepsilon}, 
g_{\varepsilon}' , q_{\varepsilon})$ (sub)-converges to a maximal flat 
limit $(F, g_{o}', q)$. Similarly, the rescalings 
$(\Omega_{\varepsilon}, \bar g_{\varepsilon}, q_{\varepsilon})$ 
(sub)-converge to a maximal flat limit $(F, \bar g_{o}, q)$, and $(F, 
\bar g_{o})$ is homothetic to $(F, g_{o}')$. Recall that the 
convergence of the blow-ups $(\Omega_{\varepsilon}, g_{\varepsilon}' , 
q_{\varepsilon})$ to the limit $(F, g_{o}', q)$ is smooth away from 
$\partial F$ and that $\partial F$ is formed by the limits of points 
$z_{\varepsilon}$ where the curvature is concentrating at a higher 
rate, i.e. $\rho' (z_{\varepsilon}) \rightarrow $ 0, c.f. the 
discussion preceding (1.23). Of course $\partial F \neq  \emptyset  $ 
and by Proposition 2.1, $\partial F$ is the $q_{\varepsilon}$-based 
Hausdorff limit of the level $L^{\upsilon_{o}}$; the same holds w.r.t. 
$(F, \bar g_{o})$. The center points $y_{\varepsilon}$ converge to a 
point $y\in\bar F$ with $q\in \bar B_{y} (1+d_{o})$, while the domains 
$\mathcal{C_{\varepsilon}}$ converge to the domain $\mathcal{C}_{o} 
\subset \bar A_{y}(1)$.

\smallskip

 The results above now combine easily to give the existence of natural 
2-spheres.

\begin{theorem} \label{t 5.9.}
  Let $q_{\varepsilon}$ be a distinguished base point of 
$(\Omega_{\varepsilon}, g_{\varepsilon})$ satisfying the non-collapse 
assumption NCA II and let $(F, g_{o}', q)$ be a maximal flat blow-up 
limit constructed above. Then there is a non-empty compact subset 
$\Sigma_{o}$ of $\partial F$, such that either $\partial F\setminus 
\Sigma_{o} = \emptyset  $ or
\begin{equation} \label{e5.28}
dist_{g_{o}'}(\Sigma_{o}, \partial F\setminus \Sigma_{o}) \geq  d_{o} > 
0. 
\end{equation}
Further, $\mathcal{C}_{o}$ is embedded in ${\mathbb R}^{3},$
\begin{equation} \label{e5.29}
\mathcal{C}_{o} \subset  {\mathbb R}^{3}, 
\end{equation}
and $\mathcal{C}_{o}$ is isometric to a spherical annulus $S^{2}\times 
I$ about a point in ${\mathbb R}^{3}.$
\end{theorem}
\noindent
{\bf Proof:}
For convenience we work with the metric $(F, \bar g_{o})$ in place of 
$(F, g_{o}')$; as noted above, these metrics are homothetic. Returning 
to the sequence $(\Omega_{\varepsilon}, \bar g_{\varepsilon})$ as in 
(5.24), suppose first that the annuli $\bar A_{\varepsilon}(1,1+d_{o})$ 
are connected, (in a subsequence); more generally, one may suppose that 
$\bar A_{\varepsilon}(1,1+d_{o}) \subset  U^{o}$. Then we claim that 
$L^{\upsilon_{o}}\cap\bar B_{\varepsilon}(1) \neq  \emptyset  ,$ and 
so, (by Proposition 1.1 as usual),
\begin{equation} \label{e5.30}
\bar \rho \rightarrow  0, 
\end{equation}
somewhere within $\bar B_{\varepsilon}(1).$ For if (5.30) does not 
hold, then $L^{\upsilon_{o}}$ lies outside $\bar B_{\varepsilon}(1)$, 
and $\bar g_{\varepsilon}$ converges smoothly on $\bar 
B_{\varepsilon}(1)$ to the limit $\bar B_{y}(1)$. Since $\bar 
A_{\varepsilon}(1,1+d_{o}) \subset  U^{o}, L_{o} \subset  \bar 
B_{\varepsilon}(1)$ bounds a domain $D_{\varepsilon}$ in $\bar 
B_{\varepsilon}(1)$ and hence $\nu $ has a  minimal value in 
$D_{\varepsilon}$. This gives the same contradiction as in Lemma 5.7 
above. 

 Thus, a singularity forms at a point $z_{\varepsilon}\in  \bar 
B_{\varepsilon}(1),$ i.e. there are points $z_{\varepsilon}$ within 
$\bar B_{\varepsilon}(1)$ of much higher curvature concentration than 
that within $\mathcal{C}_{\varepsilon} .$ In particular, there is a 
non-empty component of $\partial F$ within the limit $\bar B_{y}(1)$ of 
$\bar B_{\varepsilon}(1)$. As noted above, $\bar \rho $ is bounded away 
from 0 in $\mathcal{C}_{\varepsilon}$, so that 
\begin{equation} \label{e5.31}
\mathcal{C}_{o} \cap  \partial F = \emptyset .
\end{equation}

 Suppose instead $L^{\upsilon_{o}}\cap\bar B_{\varepsilon}(1) = 
\emptyset$, so that the balls $\bar B_{\varepsilon}(1)$ converge 
smoothly to the limit $\bar B_{y}(1)$ and surrounding annulus $\bar 
A_{y}(1,1+d_{o})$ in $F$. This limit annulus is not connected, since 
the annuli $\bar A_{\varepsilon}(1,1+d_{o})$ are not connected. Hence, 
the limit $\bar B_{y}(1)$ is a geodesic ball in a smooth flat manifold, 
with disconnected boundary. This can only occur if $\bar B_{y}(1)$ is 
isometric to a ball in a flat product $F_{prod} = {\mathbb R}\times 
S^{1}\times S^{1}$. Now by Lemma 4.9, the distance function $t$ to 
$\partial F$, (w.r.t. $\bar g_{o}$), is unbounded on $F$, while by 
(2.19), $t = \rho $ on $F$. Hence, one may choose a sequence of points 
$y_{i}$ with $t(y_{i}) \rightarrow  \infty ,$ approximate them by 
points $y_{\varepsilon_{i}}\in  (\Omega_{\varepsilon_{i}}, 
g_{\varepsilon_{i}}),$ and consider the rescalings based at 
$y_{\varepsilon_{i}}.$ This has the effect of blowing down $F$. 
However, the flat product structure of $(F, \bar g_{o})$ on $\bar 
B_{y}(1)$ extends to all of $(F, \bar g_{o})$ and thus it follows that 
the rescalings based at $y_{\varepsilon_{i}}$ are collapsing. This 
contradicts the NCA II, (5.27). Actually, to apply Lemma 4.9, one needs 
to know that $|u-1|(q_{\varepsilon}) \geq \bar \varepsilon^{4\mu}$. 
By (4.44), it suffices to show that $\delta_{a}(q_{\varepsilon}) \geq 
\bar \varepsilon^{4\mu}$, which follows from (5.41) below; (the proof 
of (5.41) is independent of Theorem 5.9).

  Thus $\Sigma_{o} = \partial F\cap \bar B_{y}(1)$ is non-empty and the 
estimate (5.28), (w.r.t. $\bar g_{o}$), follows as in (5.25). The 
domain $\mathcal{C}_{o}$ from (5.31) is a connected geodesic annulus 
about the point $y\in\bar F,$ separating $\Sigma_{o}$ from the rest of 
$\partial F.$

 Suppose now that $\bar B_{y}(1+d_{o})$ is simply connected. Then the 
developing map $\mathcal{D} $ is an isometric immersion of $\bar 
B_{y}(1+d_{o})$ into the ball $B_{0}(1+d_{o})$ about $0$ in ${\mathbb 
R}^{3}$. To prove (5.29) and the fact that $\mathcal{C}_{o}$ is a 
spherical annulus, it suffices to prove that $\mathcal{D} $ is an 
embedding, or equivalently $\mathcal{D} $ has no holonomy on $\bar 
B_{y}(1+d_{o})$. As in the discussion following Proposition 4.8, each 
holonomy transformation is either a rotation, twist or translation. Now 
$\mathcal{D} $ can have no rotational holonomy, generated by rotation 
about an axis within $\bar B_{y}(1+d_{o})$. For such an axis, 
necessarily a line, must be contained in $\partial F,$ which is 
impossible by the gap estimate (5.28). Hence the only holonomy 
transformations are those given by twists or translations. These both 
introduce non-trivial loops, and thus $\bar B_{y}(1+d_{o})$ could not 
be simply connected.

 On the other hand, if $\bar B_{y}(1+d_{o})$ is not simply connected, 
then the arguments above imply that the developing map $\mathcal{D}$ 
maps $F$ into a non-trivial quotient ${\mathbb R}^{3}/\Gamma ,$ for 
some $\Gamma $ acting freely and isometrically on ${\mathbb R}^{3};$ 
hence $\Gamma  = {\mathbb Z} $ or ${\mathbb Z}\oplus{\mathbb Z}$ acts 
by twists and/or translations. Via Lemma 4.9 as above, this contradicts 
the NCA II. This completes the proof.
{\endproof}

 Since the convergence of $(\Omega_{\varepsilon}, g_{\varepsilon}' , 
q_{\varepsilon})$ to the limit $(F, g_{o}', q)$ is smooth away from 
$\partial F,$ there are 2-spheres $S^{2} = S_{\varepsilon}^{2}$ 
embedded in $(\Omega_{\varepsilon}, g_{\varepsilon}' )$ converging 
smoothly to the limit $S^{2} \subset  \mathcal{C}_{o} \subset $ $(F, 
g_{o}')$. The metric $g_{\varepsilon}' $ is just a rescaling of 
$g_{\varepsilon},$ and so the metrics $g_{\varepsilon}$ are crushing 
the 2-spheres $S_{\varepsilon}^{2}$ to {\it  points}. Of course, since 
$\Omega_{\varepsilon}$ is weakly embedded in $M$, the spheres 
$S_{\varepsilon}^{2}$ are also embedded in $M$, i.e.
\begin{equation} \label{e5.32}
S_{\varepsilon}^{2} \subset  M. 
\end{equation}

 It follows from Theorem 5.9 that the 2-sphere $S_{\varepsilon}^{2} 
\subset  \Omega_{\varepsilon},$ although it converges to a point as 
$\varepsilon  \rightarrow $ 0, surrounds the formation of a 
singularity, in that "within" $S_{\varepsilon}^{2}$ there are regions 
of much higher curvature concentration than the curvature concentration 
on $S_{\varepsilon}^{2}$ itself.

\medskip

{\bf 5.4.}
 In this subsection, we prove the central result of the paper, Theorem 
5.12, stating that the 2-spheres $S_{\varepsilon}^{2}$ constructed 
above in (5.32) are essential 2-spheres in $M$. 

 To begin, one needs to understand the behavior of the potential 
functions $\nu_{\delta}$ on the blow-ups $g_{\varepsilon}' $ and the 
limit $F$ in Theorem 5.9. To do this, recall that the spheres 
$S_{\varepsilon}^{2}$ are contained in the connected region 
$\mathcal{C}_{\varepsilon} = A_{\varepsilon}(\lambda_{\varepsilon}) = 
A_{y_{\varepsilon}}(\lambda_{\varepsilon}, 
(1+d_{o})\lambda_{\varepsilon})$, where $q_{\varepsilon}$ is a 
distinguished base point; the annulus is connected by the proof of 
Theorem 5.9. In particular, (5.18)-(5.19) hold. Let $m_{o}'$ be the 
mass of the potential $\nu_{\delta_{o}} = (u-1)/\delta_{o}, \delta_{o} 
= \bar \varepsilon^{2\mu}$, on the level set $L_{o}$ containing 
$y_{\varepsilon}$, w.r.t. $g_{\varepsilon}'$ in (5.26). Also let 
$B_{\varepsilon}' = B_{y_{\varepsilon}}'(\lambda_{\varepsilon}')$, 
where $\lambda_{\varepsilon}' = \lambda_{\varepsilon} \cdot 
\rho(q_{\varepsilon})^{-1}$. Since the metrics $\bar g_{\varepsilon}$ 
and $g_{\varepsilon}'$ in (5.24) and (5.26) are uniformly homothetic, 
$\lambda_{\varepsilon}' \rightarrow \lambda_{o}' \in (0, \infty)$, as 
$\varepsilon \rightarrow 0$.

 Now there are three possibilities for the behavior of 
$m_{o}'(B_{\varepsilon}')$, as $\varepsilon  \rightarrow $ 0, namely
\begin{eqnarray}
& (i)   & m_{o}'(B_{\varepsilon}') \rightarrow  \infty , \nonumber \\
& (ii)  & n_{o} \leq  m_{o}'(B_{\varepsilon}') \leq  n_{o}^{-1},\\
& (iii) & m_{o}'(B_{\varepsilon}') \rightarrow  0, \nonumber
\end{eqnarray}
where $n_{o}$ is some positive constant. Observe that in case (iii), 
(5.19) implies $m_{o}'(B_{\varepsilon}') \geq  \bar \varepsilon^{\mu 
/2}.$ On the other hand, one has the natural potential function w.r.t. 
the base point $q_{\varepsilon},$ i.e. as in (2.17),
\begin{equation} \label{e5.34}
\nu_{\delta (q_{\varepsilon})} = (u-1)/\delta (q_{\varepsilon}), 
\end{equation}
and its associated mass function $m_{q_{\varepsilon}}$ on $L_{o}$, 
again in the $g_{\varepsilon}$ scale. It is important to note that here 
we are using the 'curvature' $\delta $ as in (2.10), and not 
$\delta_{a}$ from (3.1), used throughout \S 4 and so far in \S 5. Of 
course, for any fixed $\varepsilon ,$ the potentials $\nu_{\delta_{o}}$ 
and $\nu_{\delta (q_{\varepsilon})}$ are proportional, as are the mass 
functions.

 We claim first that the two potentials and mass functions above have a 
uniformly bounded proportionality constant, when (5.33)(ii) holds.

\begin{lemma} \label{l 5.10.}
  If (5.33)(ii) holds, then there is a constant $n_{1} = n_{1}(n_{o}) > 
$ 0, independent of $\varepsilon$, s.t.
\begin{equation} \label{e5.35}
n_{1}\delta_{o} \leq  \delta (q_{\varepsilon}) \leq  
n_{1}^{-1}\delta_{o}. 
\end{equation}
\end{lemma}

\noindent
{\bf Proof:}
 By assumption, the mass $m_{o}'$ of $\nu_{\delta_{o}}$ in 
$B_{\varepsilon}'$ is on the order of 1. On the other hand, the mass 
$m_{\delta_{a}}$ of $\nu_{\delta_{a}}, \delta_{a} = 
\delta_{a}(q_{\varepsilon}),$ is also on the order of 1 in 
$B_{\varepsilon}'$; this follows from (3.13), the NCA II in (5.27), 
together with Proposition 3.3. Hence, from the definitions of the 
masses, it follows that
\begin{equation} \label{e5.36}
\delta_{a} \sim  \delta_{o}, 
\end{equation}
i.e. the ratio is uniformly bounded away from 0 and $\infty $ as 
$\varepsilon  \rightarrow $ 0. Now recall from (3.2)-(3.3) that 
$\delta_{a} \geq  c_{1}\cdot \delta $ and $\delta_{a} >>  \delta $ if 
and only if the limit $\nu_{a}$ of $\nu_{\delta_{a}}$ is an {\it 
affine}  function; here of course $\delta  = \delta (q_{\varepsilon}).$ 
 However, an affine function has no gap in the support of the (limit) 
measure $d\mu ,$ or equivalently no mass gap. Hence by (5.28) one must 
have $\delta_{a} \leq  c_{1}\cdot \delta ,$ so that
\begin{equation} \label{e5.37}
\delta_{a} \sim  \delta . 
\end{equation}
Combining the estimates (5.36)-(5.37) gives (5.35).
{\endproof}

 Observe that Lemma 5.10 implies that
\begin{equation} \label{e5.38}
m_{\delta (q_{\varepsilon})}(B_{\varepsilon}') \sim  1, \ \ {\rm as} \ 
\  \varepsilon  \rightarrow  0, 
\end{equation}
so that the part of the level $L_{o}$ in $B_{\varepsilon}'$ has a 
definite amount of mass for the potential $\nu_{\delta 
(q_{\varepsilon})}.$ Since one has convergence to the limit $F$, and 
not collapse, observe that (5.5) holds here also, for the same reasons. 
Thus,
\begin{equation} \label{e5.39}
\rho^{2}(q_{\varepsilon}) <<  \delta (q_{\varepsilon}), 
\end{equation}
so that the base points $q_{\varepsilon}$ satisfy Case (i) of (2.12).

 In the cases (5.33)(i) or (iii), one needs to renormalize the mass by 
setting
\begin{equation} \label{e5.40}
\delta_{1} = \delta_{o}\cdot  m_{o}' (B_{\varepsilon}'), 
\end{equation}
so that the mass of the potential $\nu_{\delta_{1}} = (u- 
1)/\delta_{1}$ in $B_{\varepsilon}'$ is on the order of 1. By the same 
reasoning as in Lemma 5.10, it then follows that $\delta_{1} \sim 
\delta_{a}(q_{\varepsilon}) \sim  \delta (q_{\varepsilon})$. We claim 
that in both cases (5.39) still holds. In the situation of (5.33)(i), 
this is obvious, since then $\delta (q_{\varepsilon}) \geq  
\delta_{o}.$ Suppose instead (5.33)(iii) holds. Then applying the 
estimate (5.19) to (5.40) gives
\begin{equation} \label{e5.41}
\delta_{1} \geq  \bar \varepsilon^{2\mu +\mu /2}.
\end{equation}
On the other hand, on $L_{o}\cap B_{y_{\varepsilon}}(\frac{1}{2})$, in 
the $g_{\varepsilon}$ scale, one has $\rho  \leq  \bar 
\varepsilon^{2\mu}$ by (5.4). Since $q_{\varepsilon} \in  
B_{y_{\varepsilon}}((1+d_{o})\lambda_{\varepsilon})$ and 
$\lambda_{\varepsilon} \leq  \bar \varepsilon^{3\mu /2}$ by (5.18), the 
Lipschitz property of $\rho ,$ c.f. \S 1, thus gives
$$\rho (q_{\varepsilon}) \leq  \bar \varepsilon^{3\mu /2}. $$
Hence, for $\rho  = \rho (q_{\varepsilon}),$
\begin{equation} \label{e5.42}
\frac{\rho^{2}}{\delta_{1}} \leq  \frac{\bar \varepsilon^{3\mu}}{\bar 
\varepsilon^{2\mu +\mu /2}}  \rightarrow  0. 
\end{equation}

 These arguments thus imply that the distinguished base points 
$q_{\varepsilon}$ necessarily satisfy Case (i) of (2.12), i.e. satisfy 
the conclusions of Proposition 2.3. Further, (5.37) holds, so that the 
potentials $\nu_{\delta}$ converge to the limit potential $\nu $ modulo 
constants. In sum, one has:

\begin{corollary} \label{c 5.11.}
  For $q_{\varepsilon}$ a distinguished base point as in Theorem 5.9, 
the linearized equations of the blow-up limit $(F, g_{o}' , \nu , q)$ 
are the linearized static vacuum Einstein equations
\begin{equation} \label{e5.43}
r'  = D^{2}\nu , \ \  \Delta\nu  = 0. 
\end{equation}
For $\varepsilon$ sufficiently small, and within 
$B_{y}((1+d_{o})\lambda_{o}') \cap F$, $\lambda_{o}' = \lim 
\lambda_{\varepsilon}' \rho(q_{\varepsilon})^{-1}$, the metric 
$g_{\varepsilon}' $ has the form,
\begin{equation} \label{e5.44}
g_{\varepsilon}'  = (1 -  2\nu\delta )g_{o}'  + o(\delta ), 
\end{equation}
modulo diffeomorphisms, where $g_{o}'$ is the flat metric on $F$. The 
potential function $\nu $ is a non-constant harmonic function and is 
the limit of $\nu_{\delta (q_{\varepsilon})}$, mod constants.
\end{corollary}
{\endproof}

 We are now in position to prove the main result of the paper; this 
result constitutes the most important part in the proof of the Main 
Theorem of \S 0.

\begin{theorem} \label{t 5.12.}
  Let $(\Omega_{\varepsilon}, g_{\varepsilon})$ be a sequence of 
minimizing pairs for $M$, with $\varepsilon  = \varepsilon_{i} 
\rightarrow $ 0 satisfying (1.36) and suppose the degeneration 
hypothesis (4.1) and non-collapse assumptions NCA I, NCA II hold on 
some subsequence.

 Then, for $\varepsilon $ sufficiently small, the 2-sphere 
$S_{\varepsilon}^{2}$ in $M$ given by (5.32) is essential in M. In 
particular, $M$ is a reducible 3-manifold.
\end{theorem}

\noindent
{\bf Proof:}
 Suppose that $S^{2} = S_{\varepsilon}^{2}$ is inessential in $M$, i.e. 
$S^{2}$ bounds a 3-ball in $M$. By constructing a suitable metric 
glueing of a 3-ball onto $S^{2}$ in $M$ as a comparison metric, we will 
obtain a contradiction to the fact that $g_{\varepsilon}$ is a 
minimizer of $I_{\varepsilon}^{~-}$, for $\varepsilon$ sufficiently 
small.

 Given the preceding results of \S 5, most all of the work to prove 
Theorem 5.12 has already been done, explicitly for this purpose, in 
[1]. First, as in [1, Thm.3.1], the proof needs to be divided into two 
cases, according to whether or not $S^{2}$ bounds a 3-ball $B^{3}$ on 
the inside or outside. The 2-sphere $S^{2}$ is said to bound $B^{3}$ on 
the {\it  inside}  if the component $B^{3}$ of $M\setminus S^{2}$ 
contains the center base point $y_{\varepsilon}$ associated to the 
distinguished base point $q_{\varepsilon};$ otherwise $S^{2}$ bounds on 
the {\it  outside}. 

 Since $\mathcal{C}_{\varepsilon}  \subset  U^{o}$ by (5.22), one has 
$\nu_{\delta} >  \nu_{\delta}(L_{o})$ on $\mathcal{C}_{\varepsilon}$, 
$\delta  = \delta (q_{\varepsilon}),$ while Theorem 5.9 implies that 
$L_{o}\cap B_{y_{\varepsilon}}'(\lambda_{\varepsilon}') \neq  \emptyset 
 $ and $L^{\upsilon_{o}}\cap 
B_{y_{\varepsilon}}'(\lambda_{\varepsilon}') \neq  \emptyset$. Since 
$\lambda_{\varepsilon}' \rightarrow \lambda_{o}' \in (0, \infty)$, the 
factor $\lambda_{\varepsilon}'$ will be ignored, (i.e. set to 1), in 
the following; (this is the same as working with $\bar g_{\varepsilon}$ 
in place of $g_{\varepsilon}'$). By Corollary 5.11, the potential 
functions $\nu_{\delta}$ converge, modulo constants, to the limit 
potential $\nu $ on $(F, g_{o}')$ and by Proposition 4.8, $\partial F = 
\{\nu  = -\infty\}.$ It follows that there is a level set $L'  = L' 
(\varepsilon )$ of $\nu_{\delta}$ such that $L'\cap  
B_{y_{\varepsilon}}' (1+d_{o}) \subset\subset  B_{y_{\varepsilon}}' 
(1+d_{o}),$ i.e.
\begin{equation} \label{e5.45}
L'\cap  \partial B_{y_{\varepsilon}}' (1+d_{o}) = \emptyset  , 
\end{equation}
with $dist_{g_{\varepsilon}'}(L_{o}, L' ) \geq  d_{1},$ for some fixed 
$d_{1} > $ 0. 

 Theorem 5.9 shows that the ball $B_{y}(1+d_{o})$ is embedded in 
${\mathbb R}^{3},$ and hence $L'  \subset  {\mathbb R}^{3}$ also. Thus 
$L' $ bounds a compact domain $W \subset  {\mathbb R}^{3}, \partial W = 
L' ,$ containing the limit $\partial F\cap B_{y}' (1).$ In particular, 
$\nu (x) \leq  \nu (L' ),$ for all $x\in W.$ Observe also that $W$ is 
contained in the 3-ball $B(1+d_{o})\subset  {\mathbb R}^{3}.$

\smallskip

{\bf Case I.}
 $S^{2}$ bounds on the inside. The proof of this case is given in [1, 
Thm.4.2], (which was developed just for this purpose). The idea here is 
to glue the {\it flat} metric $g_{o}'$ onto the metric 
$g_{\varepsilon}'$ along the compact level set $L' .$ Namely the 
2-sphere $\bar S^{2} = S^{2}(1+\frac{1}{2}d_{o}) \subset  
\mathcal{C}_{o} \subset  F \subset  {\mathbb R}^{3}$ bounds a flat 
3-ball $\bar B^{3} \subset  {\mathbb R}^{3}$ on the inside. Since 
$S_{\varepsilon}^{2}$ bounds a 3-ball $B^{3} \subset  M$ on the inside, 
one may construct a comparison metric $\bar g_{\varepsilon}$ to 
$g_{\varepsilon}' $ essentially by replacing 
$g_{\varepsilon}'|_{B^{3}}$ with the flat metric $g_{o}'|_{\bar 
B^{3}}.$ For technical reasons, it is preferable to carry out this 
replacement argument on $(W, L')$ in place of $(B^{3}, S^{2}).$ This 
construction is explained in complete detail in [1, (4.26)ff], and we 
refer there for details. (Note that the level $L'$ above corresponds to 
$L_{o}$ in [1] and the level set results related to (5.45) are the same 
as the assumptions [1,(4.5)-(4.7)]). The upshot is that, for 
$\varepsilon $ sufficiently small, the comparison metric $\bar 
g_{\varepsilon}$ satisfies $I_{\varepsilon}^{~-}(\bar g_{\varepsilon}) 
<  I_{\varepsilon}^{~-}(g_{\varepsilon}' ),$ contradicting the 
minimizing property of $g_{\varepsilon}' .$

\smallskip

{\bf Case II.}
 $S^{2}$ bounds on the outside. The proof in this situation consists of 
two (related) subcases, both based on the construction of a comparison 
metric on the outside. The 2-sphere $S^{2}$ in $M$ given by (5.32) 
bounds a 3-ball in $M$ on the outside, and so in particular $S^{2}$ 
separates $M$ into two components, $M = M'\cup_{S^{2}}B^{3},$ the 
outside one given by $B^{3}.$ The comparison metrics are essentially 
the same as those already constructed, in essentially the same 
circumstances, in [1].

{\bf (i).}
 Suppose first that, for some $\varepsilon  > $ 0 sufficiently small,
\begin{equation} \label{e5.46}
B^{3}\cap\partial F \neq  \emptyset  , 
\end{equation}
so that there are other components of $\partial F$ outside $S^{2}$ in 
$F$. Of course, $\partial F\cap\mathcal{C}_{o} = \emptyset  $ and 
$\partial F\cap\partial B^{3} = \emptyset  .$ 

 We explain the idea first formally on the limit $F$, and then pass to 
the metrics $(\Omega_{\varepsilon}, g_{\varepsilon}).$ The round 
2-sphere $S^{2} = S^{2}(1+\frac{d_{o}}{2}) \subset  {\mathbb R}^{3}$ 
isometrically embeds in a round 3-sphere $S^{3}(R)$ of large radius 
$R$, $R = 10(1+\frac{d_{o}}{2})$ for instance. The surface $S^{2}$ 
disconnects $S^{3}(R)$ into a smaller, inside 3-ball $B^{-}$ and a 
large, outside 3-ball $B^{+},$ so that, as above, $M = 
M'\cup_{S^{2}}B^{+}.$ Let $\bar g_{o}$ be the round metric on 
$B^{+}\subset  S^{3}(R).$ This metric joins continuously to the flat 
metric $g_{o}' $ along the seam $S^{2}.$ The resulting metric across 
$S^{2}$ is then $C^{o},$ but not $C^{1}.$ However, it is clear that 
this metric may be smoothed in a small neighborhood of the seam to give 
a smooth metric $\widetilde g_{o}$ with non-negative scalar curvature, 
$\widetilde s_{o} \geq $ 0.

 Observe that $vol_{\widetilde g_{o}}B^{+} < \infty ,$ while 
$vol_{g_{o}'}(F\setminus B(1+\frac{d_{o}}{2})) = \infty ,$ (since $F$ 
is unbounded). Thus, the volume of the comparison metric is much 
smaller than that of $(F, g_{o}')$. Similarly, $z \equiv $ 0 on 
$(B^{+}, \bar g_{o})$ and the smoothing $\widetilde g_{o}$ has 
curvature $\widetilde z_{o}$ bounded near the seam, (depending only on 
the fixed smoothing). On the other hand, by (5.46), the metric $g_{o}' 
$ formally has $|z| = \infty $ on $\partial F,$ so that 
$\mathcal{Z}^{2}$ is much larger for $(F, g_{o})$ than for the 
comparison $\widetilde g_{o}.$ Finally, since the smoothing $\widetilde 
g_{o}$ satisfies $\widetilde s_{o} \geq $ 0, one has $(\widetilde 
s_{o})^{-} \equiv $ 0, so that the comparison gives no added 
contribution to $\mathcal{S}_{-}^{2}.$ Thus, formally, the comparison 
metric has a smaller value of $I_{\varepsilon}^{~-}$ than $g_{o}' .$

 It is now a simple matter to make this formal reasoning rigorous. 
Thus, choose $\varepsilon $ sufficiently small so that 
$g_{\varepsilon}' $ is very close to the limit $(F, g_{o})$ on a large 
compact domain in $F$. In particular $g_{\varepsilon}' $ is almost flat 
away from $\partial F.$ Since $(S^{2}, g_{\varepsilon}' )$ has Gauss 
curvature strictly larger than that of $S^{3}(R),$ the Weyl embedding 
theorem, c.f. [16], implies that $(S^{2}, g_{\varepsilon}' )$ 
isometrically embeds in $S^{3}(R)$ and hence bounds the domain $B^{+} = 
(B^{+})_{\varepsilon}$ as above. It follows that the metric $\bar 
g_{\varepsilon} = g_{\varepsilon}'\cup\bar g_{o},$ with 
$g_{\varepsilon}' $ restricted to $\Omega_{\varepsilon}\setminus 
B^{3}$, $\partial B^{3} = S^{2},$ is a $C^{o}$ metric on 
$\Omega_{\varepsilon}.$ As above, this metric may be smoothed to a 
metric $\widetilde g_{\varepsilon},$ satisfying $(\widetilde 
s_{\varepsilon})^{-} \geq  (s_{\varepsilon}' )^{-}$ everywhere. By 
(5.46), the curvature $\mathcal{Z}^{2}$ of $g_{\varepsilon}'$ is 
arbitrarily large, for $\varepsilon$ sufficiently small. It then 
follows by the same arguments as above in the limit case that
\begin{equation} \label{e5.47}
I_{\varepsilon}^{~-}(\widetilde g_{\varepsilon}) <  
I_{\varepsilon}^{~-}(g_{\varepsilon}' ),
\end{equation}
contradicting the minimizing property of $(\Omega_{\varepsilon}, 
g_{\varepsilon}' ).$

{\bf (ii).}
 Thus, one may assume that
\begin{equation} \label{e5.48}
B^{3}\cap\partial F = \emptyset  . 
\end{equation}
Hence $\partial F$ is compact and, (by Theorem 5.9), 
\begin{equation} \label{e5.49}
F = {\mathbb R}^{3}\setminus \partial F, 
\end{equation}
with $\partial F$ of Hausdorff dimension at most 1. It follows that the 
harmonic potential $\nu $ extends to a global subharmonic function on 
${\mathbb R}^{3}$ and so has the representation, c.f. [12],
\begin{equation} \label{e5.50}
\nu (x) = -\int_{\partial F}\frac{1}{|x- y|}d\mu_{y} + h, 
\end{equation}
where $h$ is harmonic on ${\mathbb R}^{3}$, and $d\mu_{y}$ is the Riesz 
measure of $\nu $ on $\partial F$, as preceding Proposition 4.8. Note 
that all level sets of the Newtonian potential in (5.50) are compact, 
since $\partial F$ is compact. Similarly, the function $\nu $ has a 
compact level set $L' $ in $B_{y}(1+d_{o}),$ by (5.45). It is then 
elementary to see that $h$ also has a level set with a compact 
component. Since $h$ is harmonic, this forces $h$ to be constant. 
Hence, by adding a suitable constant to $\nu ,$ one has
\begin{equation} \label{e5.51}
\nu (x) = -\int_{\partial F}\frac{1}{|x- y|}d\mu_{y}. 
\end{equation}
It follows that $\nu $ has an expansion of the form
\begin{equation} \label{e5.52}
\nu  = - m/r + O(r^{-2}), 
\end{equation}
where $r(x) = |x|$ in ${\mathbb R}^{3}$ and $m > 0$. In particular, the 
metric $g_{\varepsilon}'$ has the form (5.44), with $\nu$ as in (5.52). 
The expression (5.51) also gives $|D^{k}\nu| = O(r^{-k-1})$, for all 
$k$. (Note that by rescaling down the flat limit $F$ by arbitrarily 
large factors, one may obtain a new flat blow-up limit $\hat F$ with 
$\partial\hat F = pt$ with $\hat \nu = - m/r;$ this new limit 
corresponds {\it  formally}  to blow-ups with an asymptotically flat 
end, c.f. Remark 1.6).

 The argument in this case is again by a glueing to a large 3-ball on 
the outside, similar to (i) above. This situation is a little more 
complicated however, since there are no components to $\partial F$ 
outside $S^{2},$ so that the comparison of the $\mathcal{Z}^{2}$ term 
is more delicate. However, the construction of a comparison metric 
$\widetilde g_{\varepsilon}$ satisfying (5.47) in this case is exactly 
that already given in Case II of [1, Thm.3.1], and so we refer there 
for further details.

 In both cases one thus has a contradiction, completing the proof.
{\endproof}

\smallskip

 Theorem 5.12 implies that if $M$ is $\sigma$-tame and irreducible then 
on some sequence of minimizing pairs $(\Omega_{\varepsilon_{i}}, 
g_{\varepsilon_{i}})$ satisfying (4.3), either the degeneration 
assumption (4.1) fails, or one of the non-collapse assumptions NCA I or 
NCA II fails. To complete the proof of the Main Theorem, it remains to 
understand the situation where one of these assumptions fails. This is 
done in the final two sections.

\section{The Collapse Situation.}
\setcounter{equation}{0}

 In this section, we analyse in general the situation where there is 
collapse at base points $x_{\varepsilon}\in (\Omega_{\varepsilon}, 
g_{\varepsilon})$, $u(x_{\varepsilon}) \rightarrow 1$, on the scale of 
the curvature radius $\rho (x_{\varepsilon}).$ This situation applies 
for instance when one of the two non-collapse assumptions NCA I or NCA 
II does not hold. The main point is that collapse at $x_{\varepsilon}$ 
implies non-degeneration within bounded $g_{\varepsilon}$-distance to 
$x_{\varepsilon},$ c.f. Theorem 6.5. The remaining analysis of the 
collapse situation is then completed in \S 7 in connection with the 
non-degeneration hypothesis.

\smallskip

 Let $x_{\varepsilon}$ be any sequence of base points with 
$u(x_{\varepsilon}) \rightarrow $ 1. As discussed in \S 2, in this 
generality one may have $\rho (x_{\varepsilon}) \rightarrow $ 0, $\rho 
(x_{\varepsilon}) \rightarrow  \rho_{o} > $ 0, or $\rho 
(x_{\varepsilon}) \rightarrow  \infty ,$ (in the case $\sigma (M) =$ 
0). Consider such sequences which collapse on the scale of $\rho 
(x_{\varepsilon}),$ i.e. $\omega (x_{\varepsilon}) <<  \rho 
(x_{\varepsilon})$ as $\varepsilon  \rightarrow $ 0, where $\omega $ is 
the volume radius. As explained in \S 1, the collapse of 
$(\Omega_{\varepsilon}, g_{\varepsilon}' )$, $g_{\varepsilon}'  = \rho 
(x_{\varepsilon})^{-2}\cdot  g_{\varepsilon},$ at $x_{\varepsilon}$ may 
then be unwrapped by passing to sufficiently large, (depending on 
$\varepsilon ),$ finite covering spaces and on passing to a convergent 
subsequence, one obtains a maximal limit $(F, g_{o}', x)$. As 
previously, let $\nu_{\delta_{a}} = (u-1)/\delta_{a}(x_{\varepsilon})$, 
for $\delta_{a}(x_{\varepsilon})$ as in (3.1) and recall from (3.2) 
that $\delta_{a} \geq  c\cdot \delta$. For the remainder of the paper 
we work with $\delta_{a}$ as in (3.1), in place of $\delta $ in (2.10). 
By Theorem 2.11, the potential functions $\nu_{\delta_{a}}$ 
sub-converge to a limit potential function $\nu_{a}$, modulo addition 
of constant functions. The limit $(F, g_{o})$ has a free isometric 
$S^{1}$ or $T^{2}$ action, leaving the potential $\nu_{a}$ invariant. 

 The limit $(F, g_{o}', x)$ may be either flat or hyperbolic, i.e. of 
constant negative curvature, c.f. \S 2.1. The former case occurs 
whenever $\rho  \rightarrow $ 0, or if $\rho $ is bounded and $\sigma 
(M) =$ 0. Conversely, if $\sigma (M) < $ 0, then the limit is 
hyperbolic precisely when $\rho  \geq  \rho_{o},$ for some $\rho_{o} > 
$ 0. If $\rho  \rightarrow  \infty ,$ forcing $\sigma (M) =$ 0, limits 
could apriori be either flat or hyperbolic, c.f. \S 7.3 and Appendix A.

\smallskip

 Isometric $S^{1}$ actions on flat manifolds $F$ may be classified as 
follows. First, by lifting to a cover if necessary, the $S^{1}$ action 
may be converted to a proper ${\mathbb R}$-action. Since $F$ is flat, 
an isometric ${\mathbb R}$-action is completely determined by its 
behavior in a neighborhood $U$ of any given point $x\in F.$ One may 
assume that $U$ is embedded in ${\mathbb R}^{3}$ and hence the action 
on $U$ extends uniquely to an ${\mathbb R}$-action on ${\mathbb 
R}^{3}.$ As discussed following Proposition 4.8, the ${\mathbb 
R}$-action on ${\mathbb R}^{3}$ is generated either by a rotation about 
an axis $A$, a twist about an axis $A$, or a translation. Thus, the 
action is free everywhere except in the case of rotation, where the 
axis $A$ consists of fixed points. The same classification holds for 
$S^{1}$ actions on hyperbolic manifolds.

 If the collapse at $x_{\varepsilon}$ is a rank 2 collapse, (or a 
sequence of rank 1 collapses converging to a rank 2 collapse), then the 
the limit $(F, g_{o}', x)$ has a free isometric $T^{2}$ action, 
obtained by unwrapping both collapsing circles. The results to follow 
hold for any $S^{1}$ action with $S^{1} \subset T^{2}.$

  We begin with the following result.
\begin{lemma} \label{l 6.1}
Any collapsing sequence of base points $x_{\varepsilon}$ with 
$u(x_{\varepsilon}) \rightarrow 1$ is allowable, w.r.t. $\delta_{a}$, 
(in place of $\delta$), c.f. Definition 2.7.
\end{lemma}

\noindent
{\bf Proof:}
Base points ${x_{\varepsilon}}$ are not allowable only if, on a limit 
$(F, g_{o}', x)$, there exists a solution $\bar \nu_{\infty}$ to the 
equation 
\begin{equation} \label{e6.1}
D^{2} \bar \nu_{\infty} = g_{o}',
\end{equation}
obtained by renormalizing $\nu_{\delta_{a}}$, c.f. (2.51)-(2.53). 
However, in the collapse case, $\bar \nu_{\infty}$ must, in addition, 
be invariant under the isometric $S^{1}$-action. The equation (6.1) 
only has $S^{1}$-invariant solutions in the case of a rotational 
$S^{1}$-action, and in this case, the only solution is of the form 
$\bar \nu_{\infty} = \frac{1}{2} \cdot r^{2}$, where $r$ is the 
distance from some point on the rotation axis, (modulo affine 
functions). However, then $|\nabla \bar \nu_{\infty}|(x) = 1$, which 
implies, via the renormalization, that $|\nabla \nu_{a}| \rightarrow 
\infty$ at and near $x_{\varepsilon}$; this contradicts the definition 
of $\nu_{\delta_{a}}$, c.f. (3.13).
{\endproof}

 Throughout this section, let 
\begin{equation} \label{e6.2}
L^{o}=\{x_{\varepsilon}\in (\Omega_{\varepsilon}, g_{\varepsilon}): |u- 
1|(x_{\varepsilon}) = 1/(\ln \bar \varepsilon)^{2}\},
\end{equation}
and $U^{o}$ the corresponding $u$-superlevel set where $|u- 1| \leq  
1/(\ln \bar \varepsilon)^{2},$ i.e. $u \geq  1 - 1/(\ln \bar 
\varepsilon)^{2}$, as in (4.1); these should not be confused with the 
levels $L_{o}$ from (4.34) and their superlevels as in (5.21).

 It is useful to separate the cases of flat and hyperbolic limits. We 
begin with the flat case.

\begin{proposition} \label{p 6.2.}
  Let $x_{\varepsilon}$ be a collapsing sequence of base points in 
$U^{o}$, and let $(F, g_{o}', x, \nu)$ be a maximal flat limit as 
above. Then $\partial F$ consists of a finite number of orbits of the 
free $S^{1}$ action, together possibly with the full axis $A$. In fact, 
there is a constant $M <  \infty$, independent of $x$, such that
\begin{equation} \label{e6.3}
\#(\partial F) \leq  M, 
\end{equation}
where \# denotes the number of components. In particular, away from A, 
$\partial F$ is compact.
\end{proposition}

\noindent
{\bf Proof:}
 Since $x_{\varepsilon}$ is allowable, Proposition 4.8 and Lemma 4.9 
hold on the limit $(F, g_{o}', x)$. The potential $\nu $ is invariant 
under the $S^{1}$ action and by Proposition 4.8, $\partial F = 
I_{\infty}$. Hence $\partial F$ is a union of orbits of the action, so 
that the $S^{1}$ action extends to $\bar F.$

 The measure $d\mu_{\infty}$ of the potential $\nu $ at $\partial F$ is 
also $S^{1}$ invariant, c.f. (4.37ff). Hence, for orbits which are not 
fixed points of the action, i.e. all cases except the axis in the 
rotational case, this measure is a multiple of Lebesque measure on a 
circle. Proposition 4.8 then implies that the number of orbits not in 
$A$ is locally finite, i.e. there is a uniform bound on the number of 
orbits in any compact set in $F$.

 To obtain the global bound (6.3), recall from Lemma 4.9 that $t = 
dist_{g_{o}'}(\partial F, \cdot  )$ is unbounded on $F$. Hence, one may 
choose a point $y\in F$ with $t(y) >> $ 1, as well as base points 
$y_{\varepsilon}\in\Omega_{\varepsilon}$ such that $y_{\varepsilon} 
\rightarrow  y$. The metrics $g_{\varepsilon}'  = \rho 
(y_{\varepsilon})^{-2}\cdot  g_{\varepsilon}$ then converge to a 
rescaling $(F, g_{1}', y)$ of $(F, g_{o}', x)$, (again unwrapping 
collapse). Applying the arguments above to $(F, g_{1}', y)$ gives a 
uniform bound on the number of orbits in balls of fixed size about $y$ 
in $(F, g_{1}', y)$. Since $y$ is arbitrary, this implies (6.3).

 Finally, in the case of a rotational action, note that the full axis 
$A$ is necessarily contained in $\partial F,$ since the $S^{1}$ action 
is free on $F$ itself.
{\endproof}

 The next result shows that collapse at $x_{\varepsilon}$ as in 
Proposition 6.2 propagates to collapse at larger scales in a natural 
way. Recall the definition of the levels $L^{\upsilon_{o}}$ and 
$U^{\upsilon_{o}}$ from (1.17).

\begin{lemma} \label{l 6.3.}
  For $x_{\varepsilon}$ as above, suppose $y_{\varepsilon}$ are base 
points in the component of $U^{\upsilon_{o}}$ containing 
$x_{\varepsilon}$ such that $u(y_{\varepsilon}) \rightarrow $ 1 and 
\begin{equation} \label{e6.4}
dist_{g_{\varepsilon}}(y_{\varepsilon}, x_{\varepsilon}) \leq  c\cdot  
dist_{g_{\varepsilon}}(y_{\varepsilon}, L^{\upsilon_{o}}), 
\end{equation}
for some constant $c > 0$. Then there are constants $C_{\varepsilon},$ 
with $C_{\varepsilon} \rightarrow  \infty $ as $\varepsilon  
\rightarrow $ 0 such that if
\begin{equation} \label{e6.5}
dist_{g_{\varepsilon}}(y_{\varepsilon}, x_{\varepsilon}) \leq  
C_{\varepsilon}, 
\end{equation}
then the metrics $(\Omega_{\varepsilon}, g_{\varepsilon})$ also 
collapse at $y_{\varepsilon},$ i.e. $\omega (y_{\varepsilon}) <<  \rho 
(y_{\varepsilon}).$
\end{lemma}

\noindent
{\bf Proof:}
 Suppose first that $\rho (x_{\varepsilon}) \geq  \rho_{o},$ for some 
$\rho_{o} > $ 0, so that by Proposition 1.1, 
$dist_{g_{\varepsilon}}(x_{\varepsilon}, L^{\upsilon_{o}}) \geq  
\rho_{1},$ for some $\rho_{1} > $ 0. There is a curve $\gamma $ joining 
$x_{\varepsilon}, y_{\varepsilon},$ satisfying (6.4) and hence $\rho  
\geq  \rho_{2} > $ 0 everywhere along $\gamma ,$ so that the curvature 
of $(\Omega_{\varepsilon}, g_{\varepsilon})$ is everywhere bounded in a 
neighorhood of $\gamma .$ The result then follows from elementary 
comparison geometry, c.f. [10, Ch.8C] or [9, Ch.9.1].

 If instead $\rho (x_{\varepsilon}) \rightarrow $ 0, we work with the 
metrics $\widetilde g_{\varepsilon} = u^{2}g_{\varepsilon},$ as in the 
proof of Proposition 2.3. Note that $\widetilde g_{\varepsilon}$ is 
uniformly quasi-isometric to $g_{\varepsilon}$ within 
$U^{\upsilon_{o}},$ and hence the volume radii w.r.t. $g_{\varepsilon}$ 
and $\widetilde g_{\varepsilon}$ are uniformly equivalent within 
$U^{\upsilon_{o}}.$ Similarly, Proposition 1.1 implies that the $L^{2}$ 
curvature radii $\rho$, $\widetilde \rho$, are uniformly equivalent and 
hence it suffices to establish collapse w.r.t. $\widetilde 
g_{\varepsilon}.$ The Ricci curvature of $\widetilde g = \widetilde 
g_{\varepsilon}$ is given by
$$\widetilde r = 2(d\ln u)^{2} + \tfrac{1}{2}sg + \frac{2}{u}cg + 
\varepsilon\nabla\mathcal{Z}^{2}; $$
this follows from (1.38), via the Euler-Lagrange equations 
(1.14)-(1.15). Hence
\begin{equation} \label{e6.6}
\widetilde r \geq  \tfrac{1}{2}sg + \varepsilon\nabla\mathcal{Z}^{2}. 
\end{equation}
From (1.18), one has $|\varepsilon\nabla\mathcal{Z}^{2}| \leq  
c\varepsilon t_{\upsilon_{o}}^{-4},$ where $t_{\upsilon_{o}} = 
dist_{g_{\varepsilon}}(L^{\upsilon_{o}}, \cdot  )$ so that (6.6) gives
\begin{equation} \label{e6.7}
\widetilde r \geq  \tfrac{1}{2}su^{-2}\widetilde g -  c\varepsilon 
t_{\upsilon_{o}}^{-4}\cdot \widetilde g. 
\end{equation}
The standard volume comparison theorem, c.f. [10, Ch.5A] or [14, 
9.1.3], implies that if the Ricci curvature of $\widetilde 
g_{\varepsilon}$ is uniformly bounded below, then collapse at 
$x_{\varepsilon}$ implies collapse at $y_{\varepsilon}$ satisfying 
(6.4)-(6.5). Note that the scalar curvature $s$ of $g_{\varepsilon}$ in 
(6.7) is uniformly bounded below, so only the second term in (6.7) 
needs to be considered. To understand this, if $t_{\upsilon_{o}}$ is 
sufficiently small, pass to the blow-up scale where $\widetilde \rho' 
(x_{\varepsilon}) =$ 1. Then (6.7) translates to
$$(\widetilde r)'  \geq  \tfrac{1}{2}s' u^{-2}(\widetilde g)'  -  
c\alpha (t_{\upsilon_{o}}' )^{-4}\cdot  (\widetilde g)'.$$ 
In this scale, $t_{\upsilon_{o}}' (x_{\varepsilon}) \rightarrow $ 1. 
Since $\alpha  \rightarrow $ 0, the Ricci curvature $(\widetilde r)'$ 
is almost non-negative; its negative part decays as $(t_{\upsilon_{o}}' 
)^{-4}$, for $(t_{\upsilon_{o}}')$ large. Again, the volume comparison 
theorem implies that collapse at the base point $x_{\varepsilon}$ 
implies collapse in cones as in (6.4) for $t_{\upsilon_{o}}' $ 
arbitrarily large. This procedure may then be iterated until the base 
scale $g_{\varepsilon}$ is reached. The fact that one may allow 
$C_{\varepsilon} \rightarrow  \infty $ as $\varepsilon  \rightarrow $ 0 
then follows again from volume comparison at the base scale 
$g_{\varepsilon}.$
{\endproof}

 Lemma 6.3 implies that collapse at some initial $x_{\varepsilon} \in 
U^{o}$ propagates to collapse at any $y_{\varepsilon}$, 
$u(y_{\varepsilon}) \rightarrow 1$, in any connected cone region 
defined by (6.4)-(6.5). The limits $(F, g_{o}', y)$ - or $(F, g_{o}, 
y)$ when no rescaling is involved - at such base points 
$y_{\varepsilon}$ may either be flat or hyperbolic. Hence, one may 
iterate this process to increase $t_{\upsilon_{o}}$, i.e. replace 
$x_{\varepsilon}$ by $y_{\varepsilon}$, etc, provided the limits are 
flat, and provided the conclusions of Proposition 6.2 hold.

 Suppose that all limits based at such $y_{\varepsilon}$ satisfying 
(6.4)-(6.5) are flat; this is the case for instance if $\sigma (M) = 
0$. Hence of course $\rho 
(y_{\varepsilon})/t_{\upsilon_{o}}(y_{\varepsilon}) \rightarrow  1$ as 
$\varepsilon  \rightarrow 0$. We then claim that
\begin{equation} \label{e6.8}
\sup_{U^{\upsilon_{o}}}t_{\upsilon_{o}} \rightarrow  \infty , \ \ {\rm 
as} \ \  \varepsilon  \rightarrow  0;
\end{equation}
in fact (6.8) holds in cone regions about $x_{\varepsilon}$ as in 
(6.4). This estimate may be viewed as an analogue of Proposition 6.2 at 
the {\it base} scale $(\Omega_{\varepsilon}, g_{\varepsilon}, 
x_{\varepsilon})$.

  To prove the claim, start with an initial collapsing sequence 
$x_{\varepsilon} \in U^{o}$. On the maximal flat limit $(F, g_{o}', 
x)$, Theorem 2.11 shows that either
\begin{equation} \label{e6.9}
\Delta \nu_{a} = 0, \ \ {\rm or} \ \ \Delta \nu_{a} = -d_{o} < 0,
\end{equation}
where $d_{o} = \lim \rho^{2}|\bar d_{\varepsilon}| / \delta_{a}$; here 
$\rho^{2} / \delta_{a}$ is evaluated at $x_{\varepsilon}$ and $\bar 
d_{\varepsilon}$ is from (2.5). The fact that $x_{\varepsilon}$ is 
allowable by Lemma 6.1 implies that there is a constant $D_{o} < 
\infty$, independent of $x_{\varepsilon}$, such that $d_{o} \leq D_{o}$.

  Now if $\Delta \nu_{a} = 0$ on $(F, g_{o}', x)$, so that $\nu_{a}$ is 
an $S^{1}$-invariant harmonic function, then it is a standard 
consequence of the maximum principle that there are cones $C \subset F$ 
as in (6.4), (with $\partial F$ in place of $L^{\upsilon_{o}}$), such 
that $\nu_{a}(y) > \nu_{a}(x)$ for some $y \in C$ with $t'(y) = 
\rho'(y)$ arbitrarily large. Such points $y$ may be approximated by 
points $y_{\varepsilon} \in U^{o} \subset \Omega_{\varepsilon}$ and 
hence the process above increasing $t_{\upsilon_{o}}$ may be iterated 
at the much larger scale based at $y_{\varepsilon}$ in place of 
$x_{\varepsilon}$. This proves the estimate (6.8) in case the first 
alternative in (6.9) always holds in $U^{o}$. In fact, the same 
argument holds if 
\begin{equation} \label{e6.10}
\Delta \nu_{a} = -\mu,
\end{equation}
where $\mu$ is sufficiently small, since all limits $\nu_{a}$ are 
non-trivial by Theorem 2.11.

  Suppose instead, possibly at some stage of the iteration above, one 
has collapsing base points $x_{\varepsilon} \in U^{o}$ such that 
$\Delta \nu_{a} = -d_{o} < 0$ on the limit $(F, g_{o}', x)$. Hence for 
$\varepsilon$ small, $(\rho^{2} / \delta_{a})(x_{\varepsilon}) \geq 
\frac{1}{2}d_{o} |\bar d_{\varepsilon}^{-1}|$. On the other hand, since 
$d_{o} \leq D_{o}$, one has $(\rho^{2} / \delta_{a})(y_{\varepsilon}) 
\leq 2D_{o} |\bar d_{\varepsilon}^{-1}|$, for $\varepsilon$ small, for 
any collapsing base points $y_{\varepsilon}$ with $u(y_{\varepsilon}) 
\rightarrow 1$. Combining these estimates gives 
\begin{equation} \label{e6.11}
\delta_{a}(y_{\varepsilon}) \geq \frac{d_{o}}{4D_{o}} 
(\rho')^{2}(y_{\varepsilon}) \delta_{a}(x_{\varepsilon}), 
\end{equation}
where $\rho'$ is the $L^{2}$ curvature radius in the scale 
$g_{\varepsilon}'$ based at $x_{\varepsilon}$. Hence, 
$\delta_{a}(y_{\varepsilon})$ is essentially increasing quadratically 
in the distance $\rho'$. However, Proposition 6.2 implies that $\rho'$ 
is unbounded on the limit $(F, g_{o}', x)$. If $y_{\varepsilon}$ are 
chosen so that $y_{\varepsilon} \rightarrow y \in F$ with $\rho'(y) >> 
1$, then (6.11) implies that the mass of the potential $\nu_{y} = \lim 
\nu_{\delta_{a}(y_{\varepsilon})}$ based at $y$ on $(F, g_{o}', y)$ 
becomes arbitrarily small. This contradicts Proposition 4.8. Hence 
(6.8) holds in all cases, as claimed.

\smallskip

  As usual, Proposition 1.1 implies that whenever $t_{\upsilon_{o}} 
\geq  t_{o},$ for any fixed $t_{o} > $ 0, one has
\begin{equation} \label{e6.12}
\rho  \geq  \rho_{o} >  0, 
\end{equation}
with $\rho_{o}$ depending only on $t_{o}.$ Thus, there is no 
degeneration in such regions.

\smallskip

 Next, we turn to the case of hyperbolic limits of collapsing sequences 
$\{x_{\varepsilon}\}$. Under the assumption $\rho (x_{\varepsilon}) 
\leq  K_{o},$ for some $K_{o} <  \infty ,$ this can occur only when 
$\sigma(M) < 0$ and $\rho (x_{\varepsilon}) \geq  \rho_{o},$ for some 
$\rho_{o} > $ 0. Thus, essentially {\it  no}  rescalings of 
$g_{\varepsilon}$ are involved and one may set $g_{\varepsilon}'  = 
g_{\varepsilon}.$ The following is the analogue of Proposition 6.2.

\begin{proposition} \label{p 6.4.}
  Let $(F, g_{o}, x, \nu )$ be a maximal hyperbolic limit arising from 
a collapse at $x_{\varepsilon} \in U^{o}$. Then $\partial F$ consists 
of a uniformly locally finite number of orbits of the $S^{1}$ action, 
together possibly with the full axis $A$. Thus, there exists $M <  
\infty$, independent of $x$, such that for any geodesic 1-ball $B(1)$ 
in $\bar F$,
\begin{equation} \label{e6.13}
\#(\partial F\cap B(1)) \leq  M. 
\end{equation}
\end{proposition}

\noindent
{\bf Proof:}
 The proof is identical to the first part of Proposition 6.2. 

{\endproof}

 Unlike the flat case, the existence of a global bound in (6.13) is not 
asserted in the hyperbolic case; the hyperbolic case is not 
scale-invariant, as is the flat case. However, Lemma 4.9 implies that 
(6.8) holds also in this case.

\smallskip

 The results above may be summarized in the following result on the 
structure of the base metrics $(\Omega_{\varepsilon}, g_{\varepsilon})$ 
near collapsing base points.

\begin{theorem} \label{t 6.5.}
  Let $x_{\varepsilon} \in U^{o}$ be a collapsing sequence of base 
points. Then for $\varepsilon $ small, the metrics 
$(\Omega_{\varepsilon}, g_{\varepsilon})$ do not degenerate in the 
component of $U^{\upsilon_{o}}$ containing $x_{\varepsilon}$ where 
\begin{equation} \label{e6.14}
t_{\upsilon_{o}} \geq  t, \ \ dist(x_{\varepsilon}, y_{\varepsilon}) 
\leq  C_{\varepsilon}, 
\end{equation}
for any $t > 0$ and $C_{\varepsilon}$ as in (6.5). For any $t > 0$, the 
region defined by (6.14) is non-empty and all limits of 
$(\Omega_{\varepsilon}, g_{\varepsilon}, y_{\varepsilon})$ based in 
this region are constant curvature manifolds $(F, g_{o}, y)$ with 
scalar curvature $s_{o} \equiv  \sigma (M),$ and with a free isometric 
$S^{1}$ or $T^{2}$ action. 

 If $\sigma(M) = 0$, all limits at base points in the region (6.14) 
with $t_{\upsilon_{o}} \in [t_{o}, T_{o}]$, for some $T_{o} < \infty$, 
are flat. If $\sigma(M) < 0$, then all limits based in the same region 
are hyperbolic, and such limits have boundary $\partial F$ which is 
either empty, or consists of a uniformly locally bounded number of 
orbits of the action.
\end{theorem}
{\endproof}

 Theorem 6.5 allows one to pass to the situation where one has 
non-degeneration. The analysis in this situation is carried out next.

\section{Non-Degeneration Situations.}
\setcounter{equation}{0}

 In this final section, we complete the proof of the Main Theorem. 
Throughout this section, it is assumed that $M$ is $\sigma$-tame and 
irreducible. These assumptions are not actually always necessary, but 
making these assumptions uniformly simplifies the overall logic of the 
arguments.

 Recall from \S 4 that on $(\Omega_{\varepsilon}, g_{\varepsilon}),$ 
either the non-degeneration hypothesis (4.2) or the degeneration 
hypothesis (4.1) holds. In addition of course, either $\sigma (M) < $ 0 
or $\sigma (M) =$ 0. The analysis is divided accordingly into three 
subsections. 

 In \S 7.1 we analyse the situation where non-degeneration (4.2) holds 
with $\sigma (M) < $ 0, while \S 7.2 analyses the situation where 
degeneration (4.1) holds with $\sigma (M) < $ 0. In \S 7.3 these two 
cases are handled together when $\sigma (M) = 0$. The Main Theorem is 
proved in Theorem 7.6, (the case $\sigma (M) < 0$), and in Theorem 7.7, 
(the case $\sigma (M) = 0$). 

\smallskip

$\mathbf{7.1. \sigma (M) < 0.}$

 In this subsection, we examine the structure of 
$(\Omega_{\varepsilon}, g_{\varepsilon})$ when $\sigma (M) < $ 0 and 
the non-degeneration hypothesis holds. For convenience, recall the 
non-degeneration hypothesis (4.2): namely there exists a constant 
$\rho_{o} > $ 0 such that if $u(x_{\varepsilon}) \geq  1 - 1/(\ln \bar 
\varepsilon)^{2}$, i.e. $x_{\varepsilon} \in U^{o}$ as following (6.2), 
then
\begin{equation} \label{e7.1}
\rho (x_{\varepsilon}) \geq  \rho_{o}, 
\end{equation}
on $(\Omega_{\varepsilon}, g_{\varepsilon});$ here $\varepsilon  = 
\varepsilon_{i}$ is some sequence satisfying (4.3), with 
$\varepsilon_{i} \rightarrow $ 0.

 If $\rho  \geq  \rho_{o},$ for some $\rho_{o} > $ 0, on {\it  all}  of 
$(\Omega_{\varepsilon}, g_{\varepsilon}),$ then $M$ is tame and the 
Main Theorem is proved in [1,Thm.0.2]. In fact in this situation, the 
sequence $(\Omega_{\varepsilon}, g_{\varepsilon})$ is constant in 
$\varepsilon $ and $g_{\varepsilon} = g_{o}$ is hyperbolic, c.f. the 
discussion in \S 1. The work below is thus a generalization of that in 
[1]. 

\smallskip

 Let $L^{o}$ be the level set as in (6.2), with $U^{o}$ the superlevel 
set as above. We claim that $vol_{g_{\varepsilon}}U^{o} \rightarrow $ 1 
as $\varepsilon  \rightarrow $ 0. This follows from the following more 
general result, which will also be needed later. Let $L^{1}$ be the 
level of $u$ given by 
\begin{equation} \label{e7.2}
|u(L^{1})- 1| = \bar \varepsilon^{2\mu} 
\end{equation}
with $U^{1} = \{u \geq  1 - \bar \varepsilon^{2\mu}\}$ the superlevel 
set of $L_{1}.$ Hence $U^{1}\subset U^{o}$, (for $\varepsilon$ small). 

\begin{lemma} \label{l 7.1.}
  For $U^{1}$ as above,
\begin{equation} \label{e7.3}
vol_{g_{\varepsilon}}U^{1} \rightarrow  1, \ \ {\rm as} \ \  
\varepsilon  \rightarrow  0. 
\end{equation}
\end{lemma}

\noindent
{\bf Proof:}
 To see this, (1.8) and (4.3) imply that 1 $\leq  T \leq  1+\bar 
\varepsilon^{4\mu}.$ Since $w = uT$ and $u \leq  1$, this gives $w \leq 
 1+\bar \varepsilon^{4\mu}$. Let $v_{1}$ be the volume of the set 
$U_{1} = \{u \leq  1-\bar \varepsilon^{2\mu}\}$ and $v_{2}$ be the 
volume of $U^{1} = \{u \geq  1-\bar \varepsilon^{2\mu}\},$ so that 
$v_{1}+v_{2} = 1$. Then $w \leq 1 - \bar \varepsilon^{2\mu} + c\bar 
\varepsilon^{4\mu}$ on $U_{1}$ while $w \leq 1 + c\bar 
\varepsilon^{4\mu}$ on $U^{1}$. Since, by definition, the $L^{2}$ norm 
of $w$ equals 1, one obtains
\begin{equation} \label{e7.4}
1 \leq  v_{1}(1-\bar \varepsilon^{2\mu} + c \bar 
\varepsilon^{4\mu})^{2} + v_{2}(1+c\bar \varepsilon^{4\mu})^{2}, 
\end{equation}
which gives $v_{1}\bar \varepsilon^{2\mu} \leq  v_{2}c\bar 
\varepsilon^{4\mu},$ for $\varepsilon $ small. It follows that $v_{2} 
\geq  1- c\bar \varepsilon^{2\mu},$ which gives (7.3).

{\endproof}

 Let $V^{o}$ be the $\rho_{o}/2$ tubular neighborhood of $U^{o}.$ We 
claim that the filling radius, (c.f. [10,3.35]), of $L^{o}$ in $U_{o} = 
\Omega_{\varepsilon} \setminus U^{o}$ is at least $\rho_{o}/2,$ i.e. no 
collection of components of $L^{o}$ bound within $V^{o} \cap U_{o}$. 
For if not, then there is a domain $D_{\varepsilon}\subset V^{o}\cap 
U_{o}$ with $\partial D_{\varepsilon} \subset L^{o}$ and hence $u$ 
achieves a minimal value at some point $q_{\varepsilon}\in 
D_{\varepsilon}$. Note that $\rho  \geq  \rho_{o}/4$ everywhere on 
$D_{\varepsilon}.$ Let $(F, g_{o}, q)$ be a maximal limit based at 
$q_{\varepsilon}.$ It follows that the limit potential function $\nu $ 
has a minimal value at $q$. This contradicts Theorem 2.11, via the 
trace equation (2.65).

 Since $U^{o} \subset  V^{o}, vol_{g_{\varepsilon}}V^{o} \rightarrow $ 
1 also, and so $vol_{g_{\varepsilon}}(V^{o}\setminus U^{o}) \rightarrow 
$ 0. It follows in particular that
\begin{equation} \label{e7.5}
vol_{g_{\varepsilon}}T^{o} \rightarrow  0, 
\end{equation}
where $T^{o}$ is the $\rho_{o}/4$ tubular neighborhood of $L^{o}.$

\smallskip

 Now let $x_{\varepsilon}$ be any base point in $U^{o},$ so that $\rho 
(x_{\varepsilon}) \geq  \rho_{o}.$ Any maximal limit $(H_{x}, g_{o}, 
x)$ of $(\Omega_{\varepsilon}, g_{\varepsilon}, x_{\varepsilon}), 
\varepsilon  = \varepsilon_{i},$ is of constant curvature, with scalar 
curvature $s_{o} \equiv  \sigma (M)$ and with $u \equiv $ 1, c.f. 
Proposition 1.2. In case the sequence collapses at $x_{\varepsilon},$ 
the limit has a free isometric $S^{1}$ or $T^{2}$ action, i.e. is a 
rank 1 or 2 hyperbolic cusp. If 
$dist_{g_{\varepsilon}}(x_{\varepsilon}, L^{o}) \rightarrow  \infty ,$ 
the limit $(H_{x}, g_{o})$ is complete, while if 
$dist_{g_{\varepsilon}}(x_{\varepsilon}, L^{o})$ remains bounded, the 
limit may either be complete or incomplete; however there is no 
boundary of $H_{x}$ within the Hausdorff limit of the region $V^{o}$ 
above. Observe also that (7.5) implies that all limits $(H_{x}, g_{o}, 
x)$ for $x_{\varepsilon}\in T^{o}$ are necessarily collapsing. Hence, 
if $(\Omega_{\varepsilon}, g_{\varepsilon})$ is not collapsing at 
$x_{\varepsilon},$ then the standard volume comparison theorem, c.f. 
[10, Ch.5A], or [14, 9.1.3], implies that the limit $H_{x}$ is complete.

 The main result of this subsection is the following:

\begin{theorem} \label{t 7.2.}
   Suppose $\sigma (M) < $ 0, $M$ is irreducible and $\sigma$-tame, and 
the sequence $(\Omega_{\varepsilon}, g_{\varepsilon})$, $\varepsilon = 
\varepsilon_{i}$, satisfies the non-degeneration hypothesis (7.1). Then 
$(\Omega_{\varepsilon}, g_{\varepsilon})$ is constant in $\varepsilon$, 
i.e. $(\Omega_{\varepsilon}, g_{\varepsilon}) = (H, g_{o})$, where $H = 
\ \stackrel{\cdot}{\cup} H_{k}$ is a finite collection of complete, 
connected, hyperbolic manifolds $H_{k}$, $1 \leq  k \leq  q$, $q = q(M) 
<  \infty$, and $H$ is embedded in $M$. The metric $g_{o}$ is 
hyperbolic, satisfying
\begin{equation} \label{e7.6}
s_{g_{o}} = \sigma(M), \ \ vol_{g_{o}}H = 1. 
\end{equation}
\end{theorem}

\noindent
{\bf Proof:}
 Suppose first that there is a (complete) component $C_{\varepsilon}$ 
of $\Omega_{\varepsilon}$ such that $C_{\varepsilon}\subset U^{o},$ so 
that the sequence of metrics $g_{\varepsilon}|_{C_{\varepsilon}}$ is 
tame. As mentioned at the beginning of \S 7, the work in 
[1,Prop.2.5,Rmk.2.10] then implies that 
$g_{\varepsilon}|_{C_{\varepsilon}}$ is a constant sequence, so that 
$(C_{\varepsilon}, g_{\varepsilon}) = (C_{o}, g_{o})$ is a complete 
hyperbolic manifold of finite volume.

 Thus, $g_{\varepsilon}$ can only vary on the components of $U^{o} 
\subset \Omega_{\varepsilon}$ on which $\partial U^{o} \neq  \emptyset  
.$ In this situation, we claim first that not all of $U^{o}$ is 
collapsing, i.e. there are base points $x_{\varepsilon}\in U^{o}$ such 
that $(\Omega_{\varepsilon}, g_{\varepsilon})$ does not collapse at 
$x_{\varepsilon}.$ For suppose instead all of $U^{o}$ is collapsing as 
$\varepsilon  \rightarrow $ 0. Then $U^{o}$ is a graph manifold, and 
all based limits, when the collapse is unwrapped, are rank 1 or rank 2 
hyperbolic cusps. (Of course, any domain of bounded diameter about any 
such $x_{\varepsilon}$ has volume converging to 0, due to the collapse; 
however, there could be a large number of such regions, corresponding 
to a large number of limit cusps, whose total volume approaches 1).

 Observe that the volume of geodesic balls in the expanding end of a 
hyperbolic cusp is proportional to the volume of the corresponding 
geodesic sphere. It then follows from (7.3) and (7.5) that most all (in 
terms of volume) of the components of $U^{o}$ intersecting $\partial 
U^{o}$ are cusps which are expanding into the interior of $U^{o}.$ 
However, this implies that, for any $\varepsilon  >$ 0 sufficiently 
small, there must exist points realizing the maximal diameter of the 
orbits (circles or tori) of the F-structure, (i.e. graph manifold 
structure), within $U^{o}$. Choosing such points as base points, it 
follows that the limiting hyperbolic cusp also has a maximal value for 
the diameter of the associated $S^{1}$ or $T^{2}$ orbits. However, this 
is impossible, and thus proves the claim.

\smallskip

 Let $x_{\varepsilon}$ be any non-collapsing sequence of base points in 
$U^{o}.$ It follows then from the observation preceding Theorem 7.2 
that $(\Omega_{\varepsilon}, g_{\varepsilon}, x_{\varepsilon})$ 
sub-converges to a {\it complete} hyperbolic manifold $(H_{x}, g_{o})$, 
of finite volume; the level set $L^{o}$ satisfies 
$dist_{g_{\varepsilon}}(L^{o}, x_{\varepsilon}) \rightarrow  \infty $ 
as $\varepsilon  \rightarrow $ 0 and so has no limit on $H_{x}.$ The 
convergence is smooth and uniform on compact subsets. Any such manifold 
$H_{x}$ is of finite topological type, with a finite number of cusp 
ends each of the form $T^{2}\times {\mathbb R}^{+}$. Hence $H_{x}$ 
embeds in $\Omega_{\varepsilon}$ and also in $M$. By the Margulis 
Lemma, c.f. [17], there is a fixed lower bound on the volume of any 
such $H_{x}.$

 Now repeat this argument on $\Omega_{\varepsilon} \setminus H_{x}$. It 
follows from the Margulis Lemma that after a finite number of 
iterations, one obtains a finite number of disjoint limit hyperbolic 
manifolds $H_{k},$ with $H = \ \stackrel{\cdot}{\cup} H_{k}$ satisfying 
$vol_{g_{o}}H = 1$, i.e. (7.6) holds. As above, $H$ embeds in 
$\Omega_{\varepsilon}$ and so also in $M$. Thus, one may view the 
metrics $g_{\varepsilon}|_{H}$ as a sequence of metrics on $H$.

 The next result proves that the metrics $g_{\varepsilon}$ are rigid.
\begin{proposition} \label{p 7.3.}
Under the assumptions of Theorem 7.2, the metrics 
$g_{\varepsilon}|_{H}$ are constant, up to isometry, i.e.
\begin{equation} \label{e7.7}
(\Omega_{\varepsilon}, g_{\varepsilon}) = (H, g_{o}).
\end{equation}
\end{proposition}

\noindent
{\bf Proof:}
  This is a simple application of Theorem 2.11 and the maximum 
principle. If (7.7) does not hold, then the metrics $g_{\varepsilon}$ 
vary somewhere on $H$, say at base points $x_{\varepsilon} \rightarrow 
x \in H$. Now the limit potential $\nu$ on $H_{x}$, satisfies (2.65), 
i.e.
\begin{equation} \label{e7.8}
\Delta \nu = \lambda \nu - d,
\end{equation}
where $d \geq 0$ and $\lambda = \lim \frac{1}{4}\rho^{2}\sigma T = 
\frac{1}{4}|\sigma(M)| \lim \rho^{2}$, for $\rho = 
\rho(x_{\varepsilon})$, c.f. (2.72). Since $\rho(x_{\varepsilon})$ is 
bounded away from 0 and $\infty$, observe that $\lambda > 0$. Further, 
as noted following Theorem 2.11, the potential 
$\nu_{\delta(y_{\varepsilon})} = (u-1) / \delta(y_{\varepsilon})$ 
converges to the limit $\nu$, in that there is no constant or affine 
indeterminacy. Thus, $\nu \leq 0$. This implies that $\nu$ is also 
uniformly bounded below, i.e. there is a constant $K < \infty$ such that
\begin{equation} \label{e7.9}
\nu \geq -K > - \infty,
\end{equation}
c.f. the remark following (2.69). Hence, if $p_{j}$ is a minimizing 
sequence for $\nu$ on $H_{x}$, i.e. $\nu(p_{j}) \rightarrow \inf \nu$, 
then one must have $\Delta \nu(p_{j}) \geq - \mu_{j}$, for some 
sequence $\mu_{j} \rightarrow 0$. It then follows from (7.8) that $d = 
0$ and also $\lambda = 0$, giving a contradiction. Hence (7.7) must 
hold.

  This completes the proof of Proposition 7.3 and hence also Theorem 
7.2.

{\endproof}

$\mathbf{7.2. \sigma (M) < 0.}$

 In this section, still assuming $\sigma (M) < $ 0, we analyse the case 
where the degeneration hypothesis (4.2), opposite to (4.1) or (7.1) 
holds. The main result is that Theorem 7.2 also holds in this context.

\begin{theorem} \label{t 7.4.}
  Suppose $\sigma (M) < $ 0, $M$ is irreducible and $\sigma$-tame, and 
the sequence $(\Omega_{\varepsilon_{i}}, g_{\varepsilon_{i}})$, 
satisfies the degeneration hypothesis (4.1). Then again the conclusions 
of Theorem 7.2 hold.
\end{theorem}

\noindent
{\bf Proof:}
 In view of Theorem 5.12, one may assume that on the sequence 
$(\Omega_{\varepsilon}, g_{\varepsilon})$, $\varepsilon = 
\varepsilon_{i}$, one of the non-collapse hypotheses NCA I or NCA II 
fails. The strategy is then to reduce to the proof to that of Theorem 
7.2, using the results of \S 6. 

 First, recall that the degeneration hypothesis and Theorem 4.4 
guarantee the existence of admissible level sets $L_{o}$ as in (4.34), 
with corresponding superlevel sets $U^{o}$. Such levels are not unique, 
i.e. different admissible or preferred base points may lie on distinct 
levels $L_{o}.$ (Similarly, recall that degeneration at some point 
$x_{\varepsilon}\in L_{o}$ does not imply degeneration everywhere on 
$L_{o}).$ However, it is always the case, by (4.4), that $|u(L_{o})- 1| 
\geq  \bar \varepsilon^{2\mu}$. Thus for $L^{1}$ and $U^{1}$ as in 
(7.2), the domain $U^{1}$ is contained in all superlevel sets $U^{o}$ 
of admissible levels $L_{o}$.

  Now let $t^{1} = dist_{g_{\varepsilon}}(L^{1}, \cdot)$ and set
\begin{equation} \label{e7.10}
U^{2} = \{x\in U^{1}: t^{1}(x) \geq  \tau_{o}\}, 
\end{equation}
where $\tau_{o}$ is a small parameter. Then we have the following 
analogue of the non-degeneration hypothesis (7.1), Lemma 7.1 and (7.5) 
on $U^{2}$.

\begin{lemma} \label{l 7.5.}
  The domain $U^{2} = U^{2}(\tau_{o})$ is non-empty and there is a 
constant $\rho_{o} = \rho_{o}(\tau_{o}) > $ 0 s.t.
\begin{equation} \label{e7.11}
\rho (x_{\varepsilon}) \geq  \rho_{o}, 
\end{equation}
$\forall x_{\varepsilon}\in U^{2}$. Further, the metrics 
$g_{\varepsilon}$ collapse everywhere on the boundary $L^{2} = \partial 
U^{2}$ as $\varepsilon  \rightarrow $ 0, while
\begin{equation} \label{e7.12}
vol_{g_{\varepsilon}}U^{2} \geq  1 -  \tau_{1}, 
\end{equation}
where $\tau_{1} \rightarrow $ 0 as $\tau_{o} \rightarrow $ 0, on any 
sequence $\varepsilon  \rightarrow $ 0.
\end{lemma}

\noindent
{\bf Proof:}
 The estimate (7.11) is an immediate consequence of the definition 
(7.10) and Proposition 1.1. The main point is to prove that $U^{2}$ is 
non-empty and that $g_{\varepsilon}$ collapses everywhere on $L^{2}.$ 
From this, (7.12) follows easily.

 Thus, let $x_{\varepsilon}$ be any sequence of base points on $L^{1}.$ 
Suppose first that $\rho (x_{\varepsilon}) \geq  \rho_{1},$ for some 
$\rho_{1} > $ 0 (in a subsequence). If $(\Omega_{\varepsilon}, 
g_{\varepsilon}, x_{\varepsilon})$ does not collapse, then it 
(sub)-converges to a limit hyperbolic manifold $(H, g_{o}, x)$, 
containing at least the full geodesic ball $B_{x}(\rho_{1}).$ In 
particular then $vol_{g_{\varepsilon}}B_{x_{\varepsilon}}(\rho_{1}) 
\geq  v_{1},$ for some $v_{1} > $ 0. However, a definite part of the 
volume of $B_{x_{\varepsilon}}(\rho_{1})$ is contained in $U_{1} = 
\Omega_{\varepsilon}\setminus U^{1},$ which contradicts (7.3). Hence 
$(\Omega_{\varepsilon}, g_{\varepsilon}, x_{\varepsilon})$ collapses as 
$\varepsilon  \rightarrow $ 0. Theorem 6.5 then implies first that 
$U^{2} = U^{2}(\tau_{o})$ is non-empty, for a fixed $\tau_{o}$ small, 
and further that $(\Omega_{\varepsilon}, g_{\varepsilon})$ collapses at 
all points on $L^{2}$ within bounded distance to $x_{\varepsilon}.$ An 
examination of the proof of Lemma 7.1, c.f. (7.4ff), shows that the 
same conclusion holds if $\rho (x_{\varepsilon}) \rightarrow $ 0 
sufficiently slowly as $\varepsilon  \rightarrow $ 0, for instance 
$\rho (x_{\varepsilon}) \geq  \bar \varepsilon^{\mu /10}.$

 Next, consider base points $x_{\varepsilon}$ on $L^{1}$ where $\rho 
(x_{\varepsilon}) \rightarrow $ 0, e.g. $\rho (x_{\varepsilon}) \leq  
\bar \varepsilon^{\mu /10}$. Apriori all points on $L^{1}$ may have 
this property. As mentioned following the proof of Theorem 4.4, the 
admissible and preferred base points are 1-dense, in fact $\eta$-dense, 
within the region 
\begin{equation} \label{e7.13}
S_{\varepsilon} = \{ \bar \varepsilon^{2\mu}\leq  |u-1| \leq  2/(\ln 
\bar \varepsilon)^{2}\} \cap \{\rho  \leq  \chi_{\varepsilon} \}, \ 
{\rm where} \ \chi_{\varepsilon} \rightarrow  0 \ {\rm as} \ 
\varepsilon  \rightarrow  0.
\end{equation}
To be definite, choose $\chi_{\varepsilon} = \bar \varepsilon^{\mu 
/10}$. Clearly, $x_{\varepsilon} \in S_{\varepsilon}$. Choose $\eta$ to 
be a fixed number with $\eta < \tau_{o}/3.$

 Now as mentioned above, one of the two non-collapse hypotheses fails. 
Suppose first that NCA I fails. This means that for {\it  all}  
admissible or preferred base points $y_{\varepsilon},$ there is 
collapse at some $p_{\varepsilon}\in B_{y_{\varepsilon}}(\eta )\cap 
L_{o},$ on the scale of $\rho (p_{\varepsilon}),$ where $L_{o}$ is the 
$u$-level through $y_{\varepsilon}$. The collection of balls 
$B_{y_{\varepsilon}}(\eta)$, for $y_{\varepsilon}$ admissible, covers 
$S_{\varepsilon}$. Now Lemma 6.3 and the analysis following it implies 
that the collapse at such $p_{\varepsilon}$ propagates in cones to 
collapse on $L^{1}$. Together with the collapse result above when 
$\rho(x_{\varepsilon}) \geq \bar \varepsilon^{\mu/10}$, it follows that 
$L^{1}$ collapses on an $\eta$-dense set. Again, Lemma 6.3 or Theorem 
6.5 then imply that $L^{2} = L^{2}(\tau_{o})$ collapses everywhere and 
$U^{2}$ is non-empty.

 Suppose instead that NCA I holds, but NCA II fails. When NCA I holds, 
there are some admissible or preferred base points $y_{\varepsilon}$ 
which satisfy (5.1) and possibly some that don't. For those base points 
$y_{\varepsilon}$ that don't satisfy (5.1), the same conclusion as 
above holds, (i.e. $U^{2}$ is non-empty and $L^{2}$ collapses within 
bounded distance to $y_{\varepsilon}).$ Thus, one only needs to 
consider those base points $y_{\varepsilon}$ where (5.1) holds. For any 
such base point, {\it  all}  distinguished base points 
$q_{\varepsilon}$ within $g_{\varepsilon}$-distance $\eta$ to 
$y_{\varepsilon}$ collapse, possibly on a scale large compared with 
$\rho (q_{\varepsilon}).$ Since $\rho (q_{\varepsilon}) \rightarrow $ 
0, this means there are points $q_{\varepsilon}' ,$ with 
$dist_{g_{\varepsilon}}(q_{\varepsilon}, q_{\varepsilon}' ) \rightarrow 
$ 0, $\rho (q_{\varepsilon}' ) \rightarrow $ 0, such that $\omega 
(q_{\varepsilon}' ) <<  \rho (q_{\varepsilon}' ).$ As above, Lemma 6.3 
and Theorem 6.5 imply that $L^{2}$ is non-empty and collapses within 
bounded distance to $q_{\varepsilon}' .$ By the same $\eta$-density 
arguments as before, it follows that in all cases, $L^{2}$ collapses 
everywhere w.r.t. $g_{\varepsilon}.$

 Finally, to prove (7.12), all limits based at points in $L^{2}$ are 
hyperbolic manifolds $H$ with free isometric $S^{1}$ or $T^{2}$ action, 
with $\partial H$ consisting of a locally uniformly bounded number of 
$S^{1}$ orbits. In particular, $\partial H$ has measure 0. Thus, by the 
continuity of the volume as $\varepsilon \rightarrow 0$, (7.12) follows 
immediately from Lemma 7.1.
{\endproof}

 Lemma 7.5 proves that one has non-degeneration on $U^{2},$ with 
$vol_{g_{\varepsilon}}U^{2} \geq  1-\tau ,$ where $\tau $ may be made 
small by choosing $\tau_{o}$ small, (for all $\varepsilon $ 
sufficiently small). Further, all base points on $L^{2} = \partial 
U^{2}$ collapse as $\varepsilon  \rightarrow $ 0. The rest of the proof 
of Theorem 7.4 then proceeds exactly as in the proof of Theorem 7.2, by 
choosing $\tau $ smaller than the Margulis constant.

{\endproof}

\smallskip

  Theorems 7.2 and 7.4 now lead easily to the proof of the Main Theorem 
in case $\sigma(M) < 0$.

\begin{theorem} \label{t 7.6.}
  Suppose $\sigma (M) < $ 0, and $M$ is irreducible and $\sigma$-tame. 
Then $M$ has a decomposition, unique up to isotopy, as
\begin{equation} \label{e7.14}
M = H \cup_{\mathcal{T}}G, 
\end{equation}
where $H$ is a finite collection of complete, connected, hyperbolic 
manifolds with scalar curvature $\sigma (M)$ and total volume 1, and 
$G$ is a finite collection of connected graph manifolds. The union is 
along a finite collection of tori $\mathcal{T}  = \partial H = \partial 
G,$ each incompressible in M. The Sigma constant $\sigma (M)$ is given 
by
\begin{equation} \label{e7.15}
\sigma (M) = - (6vol_{g_{-1}}H)^{2/3}, 
\end{equation}
where $vol_{g_{-1}}H$ is the total volume of $H$ in the metric of 
constant curvature $- 1.$
\end{theorem}

\noindent
{\bf Proof:}
 Theorems 7.2 and 7.4 prove the existence of the hyperbolic manifold 
$H$ embedded in $M$ satisfying (7.7). In particular 
$\Omega_{\varepsilon}$ embeds in $M$. As mentioned in \S1, the 
complement of $\Omega_{\varepsilon}$ in $M$ is then a finite union $G$ 
of connected graph manifolds, collapsed under a minimizing sequence for 
$I_{\varepsilon}^{~-}$ for any $\varepsilon > 0$. This gives the 
decomposition (7.14). The remaining parts of Theorem 7.6, namely that 
$\mathcal{T} = \partial H = \partial G$ is incompressible in $M$, that 
the decomposition (7.14) is unique up to isotopy, and that (7.15) 
holds, then follow from [1, Thm.0.2]; equivalently, the family 
$(\Omega_{\varepsilon}, g_{\varepsilon})$ is tame, (since it is 
constant), and hence Theorem 7.6 follows from [1], as mentioned in \S 0.

{\endproof}

$\mathbf{7.3. \sigma (M) = 0.}$

 In this section, the cases of non-degeneration (4.2) and degeneration 
on $L_{o}$ are handled together. The following result proves the Main 
Theorem in case $\sigma (M) = 0$.

\begin{theorem} \label{t 7.7.}
  Suppose $\sigma (M) = 0$, and $M$ is irreducible and $\sigma$-tame. 
Then $M$ is a graph manifold and either $\Omega_{\varepsilon} = 
\emptyset  ,$ or $\Omega_{\varepsilon} = M$ and $g_{\varepsilon}$ is a 
sequence of flat metrics on $M$.
\end{theorem}

 Theorem 7.7 will be proved in a sequence of steps below, but logically 
is proved by contradiction. Thus, we assume that 
$(\Omega_{\varepsilon}, g_{\varepsilon})$, $\varepsilon = 
\varepsilon_{i}$, is a non-empty sequence of minimizers for 
$I_{\varepsilon}^{~-}$ as in \S 1, with $g_{\varepsilon}$ non-flat for 
all $\varepsilon  > $ 0. As noted in \S 1, this already implies that 
$M$ is not a graph manifold.

 Suppose first the non-degeneration hypothesis (7.1) holds on 
$(\Omega_{\varepsilon}, g_{\varepsilon}),$ so that there is 
non-degeneration everywhere on $U^{o} = \{x\in\Omega_{\varepsilon}: 
u(x) \geq  1/(\ln \bar \varepsilon)^{2}\}$. All based limits in this 
region are flat manifolds $(F, g_{o}, x)$, $x = \lim x_{\varepsilon}, 
x_{\varepsilon}\in U^{o}.$ Such limits are necessarily non-compact; in 
fact by Lemma 4.9, $dist_{g_{o}}(\partial F, \cdot )$ is unbounded. 
Thus, $(F, g_{o})$ has infinite volume and so the sequence 
$(\Omega_{\varepsilon}, g_{\varepsilon}, x_{\varepsilon})$ {\it must} 
collapse. (This is in strong contrast to the situation when $\sigma (M) 
<  0$). Hence all such limits have a free isometric $S^{1}$ or $T^{2}$ 
action. Further, by Proposition 6.2 or Theorem 6.5, the function $t^{o} 
= dist_{g_{\varepsilon}}(L^{o}, \cdot)$ is unbounded, i.e. assumes 
arbitrarily large values as $\varepsilon \rightarrow 0$. 

 Next, suppose instead that the degeneration hypothesis (4.1) holds. 
Then as in the proof of Theorem 7.4, one of the non-collapse 
assumptions must fail. In this case, let $L^{1}$ and $U^{1}$ be as in 
(7.2). The proof of all of Lemma 7.5 except the volume statement (7.1) 
carries over without any change and implies that one has 
non-degeneration on $U^{2},$ i.e. (7.11) holds, and there is collapse 
everywhere on $L^{2} = \partial U^{2}$, for $U^{2}$ as in (7.10). 
Again, all limits $(F, g_{o}, x)$ of $(\Omega_{\varepsilon}, 
g_{\varepsilon}, x_{\varepsilon})$ for $x_{\varepsilon}\in U^{2}$ are 
flat manifolds, collapsed and so with free isometric $S^{1}$ or $T^{2}$ 
action. As before, the distance function $t^{1}$ from (7.10) is 
unbounded as $\varepsilon  \rightarrow 0$.

 To unify these two situations, let $\bar L = L^{o}, \bar U = U^{o}$ as 
above in case (7.1) holds, while let $\bar L = L^{2}$, $\bar U = U^{2}$ 
as above in case (7.1) fails. Hence, all limits $(F, g_{o}, x)$ of 
$(\Omega_{\varepsilon}, g_{\varepsilon}, x_{\varepsilon}), 
x_{\varepsilon}\in\bar U$ are flat manifolds, with free isometric 
$S^{1}$ or $T^{2}$ action, with $\bar t = dist_{g_{\varepsilon}}(\bar 
L, \cdot)$ unbounded as $\varepsilon \rightarrow 0$.

\smallskip

 As in the previous sections, it is then natural to consider the 
rescaled metrics
\begin{equation} \label{e7.16}
g_{\varepsilon}'  = \rho (x_{\varepsilon})^{-2}\cdot  g_{\varepsilon}. 
\end{equation}
Of course $\rho (x_{\varepsilon}) \geq  \rho_{o} > $ 0 within $\bar U.$ 
In the region where $\rho (x_{\varepsilon}) \rightarrow  \infty $ as 
$\varepsilon  \rightarrow $ 0, (i.e. where $\bar t(x_{\varepsilon}) 
\rightarrow  \infty $ as $\varepsilon  \rightarrow $ 0), the metrics 
$g_{\varepsilon}$ are then being ``blown-down''. As before $\rho' 
(x_{\varepsilon}) =$ 1, and any sequence $(\Omega_{\varepsilon}, 
g_{\varepsilon}' , x_{\varepsilon})$ has a subsequence converging 
smoothly to a maximal limit $(N_{x}, g', x)$. Since the unscaled 
sequence $(\Omega_{\varepsilon}, g_{\varepsilon}, x_{\varepsilon})$ is 
collapsing, the rescaled sequence (7.16) is collapsing even more, so 
that all limits $(N_{x}, g', x)$ necessarily have a free isometric 
$S^{1}$ or $T^{2}$ action.

\smallskip

 The metrics $g_{\varepsilon}' $ satisfy the Euler-Lagrange equations 
(2.24). Here $\bar \alpha = \varepsilon /T\rho^{2}, \rho  = \rho 
(x_{\varepsilon})$ so that obviously $\bar \alpha \rightarrow $ 0. 
Thus, as before in \S 2, the terms in (2.24) containing $\bar 
\alpha\nabla\mathcal{Z}^{2}$ and $\bar \alpha|z|^{2}$ converge to 0. On 
the other hand, the constant term $\bar c_{\varepsilon}' $ in (2.24) is 
given by $\bar c_{\varepsilon}'  = \rho^{2}\cdot \bar c_{\varepsilon},$ 
so that $\bar c_{\varepsilon}'  >>  \bar c_{\varepsilon}$ as 
$\varepsilon  \rightarrow $ 0. Recall from (1.6) that 
\begin{equation} \label{e7.17}
\bar c_{\varepsilon} = \frac{1}{12}\bar \sigma + \frac{1}{6}\bar 
\varepsilon\mathcal{Z}^{2} \rightarrow  0, \ \ {\rm as} \ \  
\varepsilon  \rightarrow  0. 
\end{equation}
It follows that the limit metric $g' $ is a solution of the equations
\begin{equation} \label{e7.18}
L^{*}u = (-\frac{1}{4}su + \bar c_{\infty})g, 
\end{equation}
\begin{equation} \label{e7.19}
2\Delta u + \frac{1}{4}su = - 3\bar c_{\infty}. 
\end{equation}
Here $s$, ($= s_{g'}$), is the scalar curvature of the limit metric 
$g'$, and the prime has been dropped from the notation. The constant 
$\bar c_{\infty}$ is given by $\bar c_{\infty} = \lim \bar 
c_{\varepsilon}'$. Apriori, this could be infinite, but the following 
simple Lemma shows it to be finite.

\begin{lemma} \label{l 7.8.}
  There is a constant $K <  \infty $ such that, for any sequence of 
base points $x_{\varepsilon}$ in $\bar U,$
\begin{equation} \label{e7.20}
0 \leq  \bar c_{\infty} \leq  K. 
\end{equation}
\end{lemma}

\noindent
{\bf Proof:}
 To see this, return to (2.24) in the scale (7.16) and note that, by 
definition, 0 $\leq  u \leq $ 1. Let $\eta $ be a bounded smooth cutoff 
function supported in $B_{x_{\varepsilon}}' (1).$ Then multiply the 
trace equation (2.24) by $\eta $ and integrate by parts over 
$B_{x_{\varepsilon}}' (1).$ Since $\bar \alpha \rightarrow $ 0, while 
$u$ and $\Delta\eta $ are bounded, it follows that $\int\eta\bar 
c_{\varepsilon}' $ is also bounded as $\varepsilon  \rightarrow $ 0, 
independent of $x_{\varepsilon}.$ Hence $\bar c_{\varepsilon}' $ is 
also uniformly bounded, which gives (7.20).
{\endproof}

 As noted in \S 2.1, Proposition A below shows that all maximal limits 
$(N_{x}, g', x)$ of $(\Omega_{\varepsilon}, g_{\varepsilon}' , 
x_{\varepsilon}), x_{\varepsilon}\in\bar U,$ are of constant curvature, 
either flat or hyperbolic, i.e. of constant negative curvature. As 
previously, denote such limits as $(F, g_{o}', x)$. The limit $F$ is 
flat when $\bar c_{\infty} =$ 0, hyperbolic when $\bar c_{\infty} > $ 
0. Recall, c.f. Remark 2.5(ii), that Theorem 2.11 holds for such 
limits, as does Proposition 6.2.

\smallskip

 Now, for any given $\varepsilon  > $ 0, $\bar c_{\varepsilon}' = 
\rho^{2} \bar c_{\varepsilon} >  0$, and hence Lemma 7.8 implies that 
$\rho $ has a uniform upper bound on $\bar U$, depending on 
$\varepsilon$, i.e. for $\varepsilon$ sufficienty small,
\begin{equation} \label{e7.21}
\rho_{max} = \sup_{\bar U} \rho (x_{\varepsilon}) \leq  2(K/\bar 
c_{\varepsilon})^{1/2} <  \infty . 
\end{equation}
On the other hand, $\rho_{max} \rightarrow  \infty $ as $\varepsilon  
\rightarrow 0$, since $\bar t$ is unbounded and $\sigma(M) = 0$. The 
next result characterizes the geometry of limits at points 
$x_{\varepsilon}\in \bar U$ (almost) realizing $\rho_{max}.$ 

\begin{lemma} \label{l 7.9.}
  If $x_{\varepsilon}$ satisfies $\rho (x_{\varepsilon})/\rho_{max} 
\rightarrow $ 1, then the limit $(F, g_{o}', x)$ of 
$(\Omega_{\varepsilon}, g_{\varepsilon}' , x_{\varepsilon})$ cannot be 
flat.
\end{lemma}

\noindent
{\bf Proof:} 
If the limit $(F, g_{o}', x)$ is flat, then Proposition 6.2 holds and 
so the distance function $t'$ to $\partial F$ assumes arbitrarily large 
values on $F$. Hence, so does $\rho'$, and so there exist 
$y_{\varepsilon} \in \Omega_{\varepsilon}$ such that 
$\rho(y_{\varepsilon}) >> \rho(x_{\varepsilon})$. If there exist such 
$y_{\varepsilon} \in \bar U$, then one has a contradiction. On the 
other hand, this can only not be the case if $\Delta \nu_{x} \leq -\mu 
< 0$, as in (6.10), which leads to the same contradiction as following 
(6.11).
{\endproof}

 Thus, blow-downs based at points $x_{\varepsilon}$ (almost) realizing 
$\rho_{max}$ must be hyperbolic. Note that the scalar curvature 
$s_{\varepsilon}'  = \rho^{2}s_{\varepsilon}$ of $g_{\varepsilon}'$, 
$\rho = \rho(x_{\varepsilon})$, satisfies $s_{\varepsilon}'  
\rightarrow  const < 0$ on such limits. Of course, since 
$s_{\varepsilon} = -\sigma T u$ and $u \rightarrow 1$, one has 
\begin{equation} \label{e7.22}
\rho^{2} \sigma T \rightarrow const > 0.
\end{equation}

  There are now several ways to complete the proof of Theorem 7.7, but 
perhaps the simplest is the following. For $x_{\varepsilon}$ as above, 
the limit $(F, g_{o}', x)$ is hyperbolic. Proposition 6.4 (or Theorem 
6.5) describes the structure of the limit. By Lemma 4.9, the function 
$t = dist_{g_{o}'}(\partial F, \cdot)$ is unbounded on $F$. Hence, by 
choosing new base points $y_{j} \in F$ if necessary, with $t(y_{j}) 
\rightarrow \infty$, and approximating them by points $y_{\varepsilon} 
= y_{\varepsilon_{j}} \in \Omega_{\varepsilon}$, one obtains a new 
hyperbolic limit $(F, g_{o}', y)$ with $\partial F = \emptyset$, i.e. 
$(F, y)$ is complete, without boundary. (Note that $\rho'$ is constant 
on hyperbolic limits and does not increase with $t$, as is the case 
with flat limits).

  This brings one exactly to the situation of Proposition 7.3. Namely, 
the limit potential $\nu$ on $F$ satisfies (7.8)-(7.9), and the 
potentials $\nu_{\delta(y_{\varepsilon})} = (u-1) / 
\delta(y_{\varepsilon})$ converge to the limit $\nu$. The same proof as 
in Proposition 7.3 then gives a contradiction to (7.22). This completes 
the proof of Theorem 7.7.

{\endproof}

\section*{Appendix A.}
\setcounter{equation}{0}
\begin{appendix}
\setcounter{section}{1}

 In this Appendix, we prove the analogue of Proposition 1.2 when
\begin{equation} \label{eA.1}
u(x_{\varepsilon}) \rightarrow  1 \ \ {\rm and} \ \ \rho 
(x_{\varepsilon}) \rightarrow  \infty . 
\end{equation}
Of course, this situation is only possible when $\sigma (M) =$ 0 and 
the metrics $(\Omega_{\varepsilon}, g_{\varepsilon})$ are collapsing at 
$x_{\varepsilon}.$ For $x_{\varepsilon}$ satisfying (A.1), let
\begin{equation} \label{eA.2}
g_{\varepsilon}'  = \rho (x_{\varepsilon})^{-2}\cdot  g_{\varepsilon}, 
\end{equation}
and let $(F, g', x)$ be the maximal limit of a convergent subsequence 
of $(\Omega_{\varepsilon}, g_{\varepsilon}' , x_{\varepsilon}),$ 
unwrapped in covering spaces to obtain convergence. From \S 7.3, $g' $ 
is a solution of the equations (7.18)-(7.19).

\medskip
\noindent
{\bf Proposition A.}
  {\it Let $(F, g_{o}', x)$ be a maximal limit constructed above, 
satisfying the equations (7.18)- (7.19). If $\bar c_{\infty} =$ 0, then 
$(F, g_{o})$ is flat, with $u \equiv $ 1. If $\bar c_{\infty} > $ 0, 
then $(F, g_{o})$ is hyperbolic, i.e. of constant negative sectional 
curvature, again with $u \equiv $ 1.}

\medskip
\noindent
{\bf Proof:}
 From (7.19) and (2.3), (2.6), at every point of $(F, g_{o})$ one has
\begin{equation} \label{eA.3}
 -\tfrac{1}{4}su - 3\bar c_{\infty} \leq  0. 
\end{equation}
Further by the definition of $u$, $- su \geq  0$ on $F$. Hence, if 
$\bar c_{\infty} =$ 0, then $su \equiv 0$ on $F$ and the equations 
(7.18)-(7.19) are the static vacuum equations (1.25). The maximum 
principle then implies that the limit is flat, with $u \equiv $ 1.

 Thus suppose $\bar c_{\infty} > $ 0. The proof in this case is a 
rather long computation. Suppose first that at the base point $x$ where 
$u(x) = 1 = \sup u$,
\begin{equation} \label{eA.4}
 -\tfrac{1}{4}s(x) - 3\bar c_{\infty} = 0, 
\end{equation}
or equivalently, $\Delta u(x) = 0$, by (7.19) We then claim that 
essentially the same proof as Proposition 1.2 then shows that $g'$ is 
hyperbolic and $u \equiv 1$. To see this, recall that by definition on 
the sequence $(\Omega_{\varepsilon}, g_{\varepsilon}), u = - 
(s^{-})/\sigma T = - (s^{-})\rho^{2}/\sigma T\rho^{2}, \rho  = \rho 
(x_{\varepsilon}).$ The term $s^{-}\rho^{2}$ is the non-positive part 
of the scalar curvature $s_{\varepsilon}' $ of $g_{\varepsilon}' $ in 
(A.2). Passing to the limit, it follows that the trace equation (7.19) 
may be rewritten in the form
\begin{equation} \label{eA.5}
2\Delta s + \tfrac{1}{4}s^{2} - 3c'  = 0, 
\end{equation}
on the limit $(F, g_{o})$ where $c'  = \lim\rho^{4}\sigma 
c_{\varepsilon}.$ Similarly, the equation (7.18) becomes
\begin{equation} \label{eA.6}
L^{*}s + (\tfrac{1}{4}s^{2}+c' )g = 0. 
\end{equation}
Write (A.5) in the form
\begin{equation} \label{eA.7}
2\Delta s + \tfrac{1}{4}(s- (12c' )^{1/2})(s+(12c' )^{1/2}) = 0. 
\end{equation}
Let $w = s+(12c' )^{1/2}$ and note that (A.3) implies that $w \geq 0$ 
while (A.4) implies that $w(x) = 0$. Thus (A.7) is of the form
\begin{equation} \label{eA.8}
2\Delta w + fw = 0, 
\end{equation}
with $f = \frac{1}{4}(s- (12c' )^{1/2}) \leq $ 0. This equation is 
exactly analogous to (1.26) and the same argument as following (1.26) 
implies that $w \equiv $ 0, so that (A.4) holds everywhere on $F$. 
Hence $u \equiv 1$ and (7.18) or (A.6) implies that $g'$ is hyperbolic.

\smallskip

 Thus, it remains to prove that (A.4) holds at $x$. To do this, observe 
first that (A.6) implies that $D_{o}^{2}s = sz$, where $D_{o}^{2}s = 
D^{2}s -  \frac{\Delta s}{3}\cdot  g$ is the trace-free Hessian, and 
hence
\begin{equation} \label{eA.9}
\delta dD_{o}^{2}s = \delta d(sz). 
\end{equation}
Here $d$ is the exterior derivative on vector valued 1-forms defined by 
the metric, and $\delta $ is the adjoint of $d$, c.f. [6,Ch.1]. All the 
computations to follow are at the base point $x$ where $u(x) = 1$. 
Since $\sup u = 1$, one has $ds = 0$ at $x$ and an elementary 
computation then shows that $\delta d sz = - (\Delta s)z + s\delta d 
z$. Consequently, $tr(\delta d sz) = 0$ at $x$, so that
\begin{equation} \label{eA.10}
tr(\delta dD_{o}^{2}s) = 0. 
\end{equation}
at $x$. The rest of the proof is a computation of $tr(\delta 
dD_{o}^{2}s).$ First, $dD^{2}s(X,Y,Z) = R(X,Y,ds,Z)$, which vanishes at 
$x$. Hence (A.6) and the definition of $L^{*}$ in (1.5) imply that $dr 
= \delta R = 0$ at $x$, where the first equality is the Bianchi 
identity. A straightforward computation then gives
$$\delta dD^{2}s = -R \circ (D^{2}s), $$
where $ R \circ $ is the action of the curvature on symmetric 
2-tensors, c.f [6, Ch.1]. Consequently,
$$tr \delta dD^{2}s = - tr \bigl(R \circ (D^{2}s)\bigr) = - < r, 
D^{2}s> , $$
at $x$. Also, $\delta d(\Delta s\cdot  g) = - (\Delta\Delta s)\cdot  g$ 
and so it follows from (A.10) that
\begin{equation} \label{eA.11}
 -< r, D^{2}s> + \Delta\Delta s = 0, 
\end{equation}
at $x$. Expanding $r$ as $r = z + \frac{s}{3}g$ and using the equations 
(A.5)-(A.6), this implies that
$$- s|z|^{2} -  \tfrac{1}{3}s\Delta s -  \tfrac{1}{4}s\Delta s  = 0, $$
so that $s|z|^{2} = -\frac{7}{12}s\Delta s.$ However, $s <  0$, $\Delta 
s \geq $ 0 at $x$, and so $\Delta s = 0$ and $z = 0$, at $x$. This 
proves (A.4), as required.
{\endproof}

\end{appendix}

\section*{Appendix B.}
\setcounter{equation}{0}
\begin{appendix}
\setcounter{section}{2}

 Let $u(x_{\varepsilon}) \rightarrow $ 1 and let $(F, g_{o}', x)$ be a 
maximal limit of $(\Omega_{\varepsilon}, g_{\varepsilon}' , 
x_{\varepsilon}),$ as in \S 2. In this Appendix, we prove that, for 
$\delta  = \delta (x_{\varepsilon})$,
\begin{equation} \label{eB.1}
\bar \alpha\frac{\nabla\mathcal{Z}^{2}}{\delta} \rightarrow  0,
\end{equation}
uniformly in $L^{2}$ on compact subsets of $F$, proving the claims in 
(2.25) and (2.43). We also indicate that the covariant derivatives 
$\nabla^{k}z/\delta , k \geq $ 0, remain bounded locally in $L^{2}$ on 
compact subsets of $F$. Here and below, all metric quantities are 
w.r.t. $g_{\varepsilon}' $ or $g_{o}' .$

 We first prove (B.1) at the base point $x_{\varepsilon}$ and within 
$B_{x_{\varepsilon}}' (\frac{1}{2}).$ Given this, the proof that (B.1) 
holds near any $y\in F,$ with $\delta  = \delta (x_{\varepsilon}),$ 
follows easily and is given at the end.

\smallskip

 To begin, Proposition 1.4 implies that $\bar \alpha \rightarrow $ 0 
while Proposition 1.1 implies that for any $k \geq $ 0, $\nabla^{k}z 
\rightarrow $ 0 in $L^{2}$ on compact subsets of $F$. Also, by 
definition, $z/\delta $ is uniformly bounded in $L^{2}(B), B = 
B_{x_{\varepsilon}}' (\frac{1}{2}).$ Hence 
\begin{equation} \label{eB.2}
\bar \alpha\frac{|z|^{2}}{\delta} \rightarrow  0 \ \ {\rm in} \ \  
L^{2}(B).
\end{equation}

 Suppose first that $\delta  \geq  c\cdot \rho^{2},$ for some $c > $ 0, 
i.e. Case (i) or (ii) of (2.12) holds. Then the same argument as in 
(2.36) implies that
\begin{equation} \label{eB.3}
\Delta\nu_{\delta} \rightarrow  0,
\end{equation}
strongly in $L^{2}(B).$ By [6, Ch.4H] or [2,(3.8)],
\begin{equation} \label{eB.4}
\nabla\mathcal{Z}^{2} = D^{*}Dz + \tfrac{1}{3}D^{2}s -  2R\circ z + 
\tfrac{1}{2}(|z|^{2}-\tfrac{1}{3}\Delta s)g,
\end{equation}
so that the Euler-Lagrange equation (2.24) in the scale 
$g_{\varepsilon}' $ gives
\begin{equation} \label{eB.5}
\frac{\bar \alpha}{\delta}D^{*}Dz + u\frac{z}{\delta} = 
D^{2}\nu_{\delta} + o(1)
\end{equation}
in $L^{2}(B).$ Now integration by parts shows that 
\begin{equation} \label{eB.6}
\frac{\bar \alpha}{\delta}D^{*}Dz \rightarrow  0 \ \ {\rm in} \ \  
L^{-2,2},
\end{equation}
i.e. it converges to 0 when paired with any sequence of symmetric 
2-tensors which are uniformly bounded in $L^{2,2},$ (since $\bar \alpha 
\rightarrow $ 0). Thus, from (B.5),
\begin{equation} \label{eB.7}
D^{2}\nu_{\delta} = 0(1) \ \ {\rm in} \ \ L^{-2,2}.
\end{equation}
It then follows from (B.3) and elliptic regularity theory that 
$D^{2}\nu_{\delta}$ is bounded in $L^{2}(B)$ and converges strongly in 
$L^{2}(B)$ to a limit $D^{2}\nu  \equiv  h$. (Here and in similar 
arguments below, one must actually go to slightly smaller balls $B'  
\subset  B$; since this is of no consequence in these arguments, this 
will be left aside). Returning to (B.5), this means that $\bar \alpha 
D^{*}Dz / \delta  + uz /\delta$ converges strongly in $L^{2}$ to $h$.  
In particular, since $z/\delta $ is bounded in $L^{2}, \bar \alpha 
D^{*}Dz / \delta $ is bounded in $L^{2},$ and thus converges weakly to 
0 in $L^{2}.$

 To prove that $\bar \alpha D^{*}Dz / \delta $ converges strongly to 0 
in $L^{2},$ i.e. (B.1), since $h$ is fixed, one has
$$<\bar \alpha D^{*}Dz/ \delta , \eta^{2}h>  \rightarrow  0, $$
where $<  , > $ denotes the $L^{2}$ inner product w.r.t. 
$g_{\varepsilon}' $ and $\eta $ is a fixed cutoff function supported in 
$B$. Letting $||\cdot ||$ denote the $L^{2}$ norm, it follows that $< 
u(z / \delta), \eta^{2}h> \rightarrow ||\eta h||^{2}$ and so the 
Cauchy-Schwarz inequality gives
\begin{equation} \label{eB.8}
||\eta u\frac{z}{\delta}||^{2} \rightarrow  K \geq  ||\eta h||^{2}.
\end{equation}
On the other hand, $||(\eta \bar \alpha D^{*}Dz / \delta )+(\eta uz / 
\delta )||^{2} \rightarrow  ||\eta h||^{2}.$ Expanding out the norm on 
the left, the cross term has the form
\begin{equation} \label{eB.9}
\frac{\bar \alpha}{\delta^{2}}\int< D^{*}Dz, u\eta^{2}z> \  = 
\frac{\bar \alpha}{\delta^{2}}\int< Dz, D(u\eta^{2}z)> \  = \frac{\bar 
\alpha}{\delta^{2}}\int u\eta^{2}|Dz|^{2}
\end{equation}
$$+\frac{1}{\delta^{2}}\int<\bar \alpha^{1/2}\eta Dz, \bar 
\alpha^{1/2}\eta du\otimes z)>  + \frac{2}{\delta^{2}}\int<\bar 
\alpha^{1/2}\eta Dz, \bar \alpha^{1/2}ud\eta\otimes z)> . $$
The last two terms are estimated as follows. Since $z/\delta $ is 
bounded in $L^{2}$ and $u \rightarrow 1$, $du \rightarrow 0$ smoothly, 
the H\"{o}lder and Young inequalities $(ab \leq  
\frac{1}{2}(a^{2}+b^{2}))$ imply that the $4^{\rm th}$ term in (B.9) is 
the sum of two terms, one small compared with the $3^{\rm rd}$ term and 
the other of order $\bar \alpha \rightarrow 0$. The same holds w.r.t. 
the last $5^{\rm th}$ term, since $d\eta$ is bounded.

 Thus the cross term is positive, or goes to 0 as $\varepsilon 
\rightarrow  0$. It follows that $||\eta \bar \alpha D^{*}Dz / 
\delta||^{2} + ||\eta uz/\delta||^{2} \leq  ||\eta h||^{2}$ in the 
limit and hence (B.1) now follows from (B.8).

\smallskip

 The argument in case $\delta  <<  \rho^{2}$ is similar, although 
further work is required since (B.3) may not hold. Since $z/\delta $ is 
uniformly bounded in $L^{2}(B)$, $div(z/\delta)$ is uniformly bounded 
in $L^{-1,2}(B).$ By the Bianchi identity, $div(z/\delta ) = 
\frac{1}{6}d(s/\delta ),$ where $s = s' $ is the scalar curvature of 
$g_{\varepsilon}' .$ Hence there are constants $m_{\varepsilon}$ such 
that $(s- m_{\varepsilon})/\delta $ remains uniformly bounded in 
$L^{2}(B);$ the constant $m_{\varepsilon}/\delta $ is the mean value of 
$s/\delta $ on $B$, and may be unbounded as $\varepsilon \rightarrow 
0$. It follows that
\begin{equation} \label{eB.10}
(\frac{s- m_{\varepsilon}}{\delta}) \rightarrow  \zeta ,
\end{equation}
weakly in $L^{2}(B)$ as $\varepsilon  \rightarrow $ 0, for some 
$\zeta\in L^{2}(B).$ Returning to the trace equation (1.31), it follows 
from (B.2) and (B.10) that there are constants $\hat c_{\varepsilon},$ 
possibly unbounded as $\varepsilon  \rightarrow $ 0, such that
\begin{equation} \label{eB.11}
\Delta\nu_{\delta} = \zeta_{\delta} + \hat c_{\varepsilon},
\end{equation}
where $\zeta_{\delta}$ is bounded in $L^{2}(B)$ and converges to $\zeta 
$ weakly in $L^{2}(B).$ Suppose the limit $(F, g_{o}', x)$ is flat. 
Then as in (2.53), the function $\phi  = \frac{1}{6}\hat 
c_{\varepsilon}r^{2}$ satisfies $\Delta\phi  = \hat c_{\varepsilon}$ on 
the limit $F$. Define $\hat \nu_{\delta} = \nu_{\delta} -  
\frac{1}{6}\hat c_{\varepsilon}r^{2};$ it then follows that $\Delta\hat 
\nu_{\delta}$ is uniformly bounded in $L^{2}(B)$ and
\begin{equation} \label{eB.12}
\Delta\hat \nu_{\delta} \rightarrow  \zeta ,
\end{equation}
weakly in $L^{2}(B).$ The same argument holds for hyperbolic limits, 
using (2.57) in place of (2.53).

 Now (B.6) holds as before and hence (B.7) holds, with $\hat 
\nu_{\delta}$ in place of $\nu_{\delta},$ or with $D^{2}_{o}$ in place 
of $D^{2}.$ Together with (B.12), it follows again from elliptic 
regularity that $D^{2}\hat \nu_{\delta}$ is uniformly bounded in 
$L^{2}(B)$ which in turn implies that $\bar 
\alpha\frac{\nabla\mathcal{Z}^{2}}{\delta}$ is uniformly bounded in 
$L^{2}(B).$

\smallskip

 To prove (B.1), by the same arguments as above, it suffices to prove 
that the convergence in (B.12) is in the strong $L^{2}$ topology. To 
prove this, return to the trace equation (2.69), i.e.
$$\Delta\nu_{\delta} = \frac{1}{2}\sigma T\rho^{2}\nu_{\delta} + 
c_{\varepsilon} + o(1) = \zeta_{\delta} + \hat c_{\varepsilon} + 
o(1),$$ 
where the constants may be unbounded as $\varepsilon  \rightarrow $ 0 
and $o(1) \rightarrow  0$ in $L^{2}$ as $\varepsilon  \rightarrow $ 0. 
In particular, $\zeta_{\delta} = \frac{1}{2}\sigma 
T\rho^{2}\nu_{\delta},$ mod additive constants. The coefficient $\sigma 
T\rho^{2}$ is uniformly bounded as $\varepsilon  \rightarrow $ 0, since 
$-\sigma T\rho^{2}u$ is the scalar curvature of $g_{\varepsilon}' .$ It 
follows that $\Delta\zeta_{\delta} = \frac{1}{2}\sigma 
T\rho^{2}\zeta_{\delta} + \frac{1}{2}\sigma T\rho^{2}\hat 
c_{\varepsilon}.$ If one pairs this equation with a cutoff function 
$\eta $ and uses the self-adjointness of $\Delta ,$ together with the 
fact that $\zeta_{\delta}$ is weakly bounded in $L^{2},$ it follows 
that $\frac{1}{2}\sigma T\rho^{2}\hat c_{\varepsilon}$ is bounded as 
$\varepsilon  \rightarrow $ 0. (This is the same proof as Lemma 7.8). 
This implies that $\Delta\zeta_{\delta}$ is bounded in $L^{2}$ and 
hence, by elliptic regularity, $\zeta_{\delta}$ converges strongly in 
$L^{2}$ to its limit $\zeta .$

\medskip

 This completes the proof of (B.1) at $x_{\varepsilon}$ within $B$. To 
prove that (B.1) holds on all of $(F, g_{o}', x)$ choose any $y\in F$ 
and let $y_{\varepsilon} \rightarrow y$. Choosing $y_{\varepsilon}$ to 
be new center points, the validity of (B.1) at $y_{\varepsilon}$ 
implies that $\bar \alpha\nabla\mathcal{Z}^{2}/\delta (y_{\varepsilon}) 
\rightarrow $ 0 in $L^{2}(B_{y_{\varepsilon}}' (\frac{1}{2})),$ where 
$B_{y_{\varepsilon}}' (\frac{1}{2})$ the the $\frac{1}{2}$-ball in the 
scale $\rho' (y_{\varepsilon}) =$ 1. Suppose these balls 
$B_{x_{\varepsilon}}' (\frac{1}{2}), B_{y_{\varepsilon}}' 
(\frac{1}{2})$ have non-empty intersection. If $x_{\varepsilon}$ is 
allowable, then Theorem 2.11, i.e. (2.64)-(2.65) holds for each on the 
intersection region. These equations and elliptic regularity imply that 
$\delta (y_{\varepsilon})/\delta (x_{\varepsilon}) \sim $ 1, i.e. the 
ratio is bounded away from 0 and $\infty $ as $\varepsilon  \rightarrow 
$ 0. Hence, a covering argument shows that (B.1) holds on all of $F$, 
provided it holds on $B_{x_{\varepsilon}}' (\frac{1}{2}).$ In fact, the 
arguments above following (B.10) give the same conclusions even if 
$x_{\varepsilon}$ is not allowable.

 Finally, suppose the base points $x_{\varepsilon}$ are allowable, so 
that, after adding an affine term to $\nu_{\delta}$, $\nu_{\delta}$ 
converges to a limit function $\nu .$ It then follows that one obtains 
$C^{\infty}$ convergence of $z/\delta $ to its limit. This will not be 
detailed here, (since its not actually needed), but the proof is 
exactly the same as that given in [2, Thm.4.2, Rmk.4.3].

\end{appendix}

\bibliographystyle{plain}

\bigskip
\begin{center}
December, 2001
\end{center}

\medskip

\noindent

\address
\noindent
{Department of Mathematics\\
S.U.N.Y. at Stony Brook\\
Stony Brook, N.Y. 11794-3651\\}
\email{anderson@math.sunysb.edu}

\end{document}